\documentclass[reqno]{amsart}

\usepackage{mathptmx,amsmath,amssymb,amsthm,graphicx,tikz,mathtools}
\usepackage{multirow,hyphenat,enumitem}

\usetikzlibrary{decorations.pathreplacing,decorations.markings,patterns}

\usepackage[colorlinks=true,linkcolor=red!70!black,citecolor=green!70!black,
    urlcolor=magenta!70!black,backref]{hyperref}

\setlist[enumerate,1]{font=\bfseries,label=\arabic*.}

\providecommand{\affiliation}{\address}
\providecommand{\acknowledgments}{\section*{Acknowledgments}}

\colorlet{darkblue}{blue!70!black}
\colorlet{darkred}{red!70!black}
\colorlet{darkgreen}{green!70!black}
\colorlet{lightyellow}{yellow!50!white}

\newcommand{\arctanh}{\operatorname{arctanh}}
\newcommand{\Candle}{\operatorname{Candle}}
\newcommand{\chord}{\operatorname{chord}}
\newcommand{\diam}{\operatorname{diam}}
\newcommand{\Int}{\operatorname{Int}}
\newcommand{\LCD}{\operatorname{LCD}}
\newcommand{\Ric}{\operatorname{Ric}}
\newcommand{\RRic}{\operatorname{\sqrt{R}ic}}
\newcommand{\Tr}{\operatorname{Tr}}
\newcommand{\Vis}{\operatorname{Vis}}
\renewcommand{\div}{\operatorname{div}}
\newcommand{\sech}{\operatorname{sech}}

\newcommand{\E}{\mathbb{E}}
\renewcommand{\H}{\mathbb{H}}
\newcommand{\Q}{\mathbb{Q}}
\newcommand{\R}{\mathbb{R}}
\renewcommand{\S}{\mathbb{S}}
\newcommand{\Z}{\mathbb{Z}}
\newcommand{\CH}{\mathbb{C}\mathrm{H}}

\newcommand{\lab}{{\mathrm{lab}}}
\newcommand{\hmu}{\hat{\mu}}
\newcommand{\tN}{\tilde{N}}
\newcommand{\dd}{\mathrm{d}}
\newcommand{\ishort}{{\mathrm{short}}}
\newcommand{\ilong}{{\mathrm{long}}}
\newcommand{\sdd}{\:\dd}
\newcommand{\del}{\partial}
\newcommand{\dbyd}[2]{\frac{\del #1}{\del #2}}
\newcommand{\eps}{\epsilon}
\renewcommand{\setminus}{\smallsetminus}

\newcommand{\dOmega}{{\del\Omega}}
\newcommand{\dW}{{\del W}}
\newcommand{\dM}{{\del M}}

\newcommand{\hF}{\hat{F}}

\newcommand{\braket}[1]{\langle {#1} \rangle}
\renewcommand{\le}{\leqslant}
\renewcommand{\ge}{\geqslant}
\newcommand{\longto}{\longrightarrow}
\newcommand{\st}{\mid}
\newcommand{\defeq}{\stackrel{\mathrm{def}}{=}}

\newtheorem{theorem}{Theorem}[section]
\newtheorem{conjecture}[theorem]{Conjecture}
\newtheorem{question}{Question}[section]
\newtheorem{proposition}[theorem]{Proposition}
\newtheorem{corollary}[theorem]{Corollary}
\newtheorem{lemma}[theorem]{Lemma}
\newtheorem{linprog}[theorem]{LP Problem}
\theoremstyle{definition}
\newtheorem*{definition}{Definition}
\theoremstyle{remark}
\newtheorem*{remark}{Remark}

\newcommand{\ie}{\emph{i.e.}}
\newcommand{\eg}{\emph{e.g.}}
\newcommand{\Eg}{\emph{E.g.}}
\newcommand{\Ie}{\emph{I.e.}}

\newenvironment{eq}[1]{\begin{equation}\label{#1}}
    {\end{equation}\ignorespacesafterend}

\newcommand{\back}{\hspace{-.5em}}
\newcommand{\eatline}{\vspace{-\baselineskip}}

\begin{document}

\title[Le Petit Prince]{The Cartan-Hadamard conjecture and The Little Prince}
\author{Beno\^{\i}t R. Kloeckner}
\email{benoit.kloeckner@u-pec.fr}
\thanks{Supported by ANR grant ``GMT'' JCJC - SIMI 1 - ANR 2011 JS01 011 01}
\affiliation{Universit\'e Paris-Est, Laboratoire d'Analyse et de
    Mat\'ematiques Appliqu\'ees (UMR 8050), UPEM, UPEC, CNRS, F-94010,
    Cr\'eteil, France}
\author{Greg Kuperberg}
\email{greg@math.ucdavis.edu}
\thanks{Supported by NSF grant CCF \#1013079 and CCF \#1319245.}
\affiliation{Department of Mathematics, University of California, Davis}

\begin{abstract}
The generalized Cartan-Hadamard conjecture says that if $\Omega$ is a domain
with fixed volume in a complete, simply connected Riemannian $n$-manifold
$M$ with sectional curvature $K \le \kappa \le 0$, then $\dOmega$ has
the least possible boundary volume when $\Omega$ is a round $n$-ball with
constant curvature $K=\kappa$.  The case $n=2$ and $\kappa=0$ is an old
result of Weil.  We give a unified proof of this conjecture in dimensions
$n=2$ and $n=4$ when $\kappa=0$, and a special case of the conjecture
for $\kappa < 0$ and a version for $\kappa > 0$.  Our argument uses a new
interpretation, based on optical transport, optimal transport, and linear
programming, of Croke's proof for $n=4$ and $\kappa=0$.  The generalization
to $n=4$ and $\kappa \ne 0$ is a new result.  As Croke implicitly did, we
relax the curvature condition $K \le \kappa$ to a weaker candle condition
$\Candle(\kappa)$ or $\LCD(\kappa)$.

We also find counterexamples to a na\"{\i}ve version of the
Cartan{\hyp}Hadamard conjecture:  For every $\eps > 0$, there is a Riemannian
$\Omega \cong B^3$ with $(1-\eps)$-pinched negative curvature, and with
$|\dOmega|$ bounded by a function of $\eps$ and $|\Omega|$ arbitrarily large.

We begin with a pointwise isoperimetric problem called ``the problem
of the Little Prince.'' Its proof becomes part of the more general method.
\end{abstract}

\maketitle

\tableofcontents

\section{Introduction}
\label{s:intro}

In this article, we will prove new, sharp isoperimetric inequalities for a
manifold with boundary $\Omega$, or a domain in a manifold.  Before turning
to motivation and context, we state a special case of one of our main
results (Theorem~\ref{th:positive}).

\begin{theorem}
Let $\Omega$ be a compact Riemannian $n$-manifold with boundary, with $n
\in \{2,4\}$. Suppose that $\Omega$ has unique geodesics, has sectional
curvature bounded above by $+1$, and that the volume of $\Omega$ is at
most half the volume of the sphere $\S^n$ of constant curvature $1$.
Then the volume of $\dOmega$ is at least the volume of the boundary of a
spherical cap in $\S^n$ with the same volume as $\Omega$.
\label{th:first} \end{theorem}

Here and in the sequel, we say that a manifold (possibly with boundary)
has \emph{unique geodesics} when every pair of points is connected by at
most one Riemannian geodesic.  (More precisely, at most one connecting
curve $\gamma$ with $\nabla_{\gamma'} \gamma' = 0$.  We do not consider
locally shortest curves that hug the boundary to be geodesics.)

\subsection{The generalized Cartan-Hadamard conjecture}

An isoperimetric inequality has the form
\begin{eq}{e:isop}
|\dOmega| \ge I(|\Omega|)
\end{eq}
where $I$ is some function.   (We use $|\cdot|$ to denote volume and $|\del
\cdot|$ to denote boundary volume or perimeter; see Section~\ref{s:basic}.)
The largest function $I = I_M$ such that \eqref{e:isop} holds for all domains
of a Riemannian $n$-manifold $M$ is called the \emph{isoperimetric profile}
of $M$.

Besides the intrinsic appeal of the isoperimetric profile and isoperimetric
inequalities generally, they imply other important comparisons.  For example,
they yield estimates on the first eigenvalue $\lambda_1(\Omega)$ of the
Laplace operator by the Faber-Krahn argument \cite{Chavel:book}.  As a
second example, the first author has shown \cite{Kloeckner:anisometry} that
they imply a lower bound on a certain isometric defect of a continuous map
$\phi:M \to N$ between Riemannian manifolds.  Both of these applications
also yield sharp inequalities when the isoperimetric optimum is a metric
ball, which will be the case for the main results in this article.

The isoperimetric profile is unknown for most manifolds.  The main case in
which it is known is when $M$ is a complete, simply connected manifold with
constant curvature.   Let $X_{n,\kappa}$ be this manifold in dimension $n$
with curvature $\kappa$, and let $I_{n,k}$ be its isoperimetric profile.
In other words, $X_{n,\kappa} = \sqrt{\kappa}\S^n$ is a sphere of radius
$\sqrt{\kappa}$ when $\kappa > 0$; $X_{n,0} = \E^n$ is a Euclidean space;
and $X_{n,\kappa} = \sqrt{-\kappa}\H^n$ is a rescaled hyperbolic space when
$\kappa < 0$.  Then a metric ball $B_{n,\kappa}(r)$ has the least boundary
volume among domains of a given volume.  Thus the profile is given by
\[I_{n,\kappa}(|B_{n,\kappa}(r)|) = |\del B_{n,\kappa}(r)|.\]
Moreover, the volume $|B_{n,\kappa}(r)|$ and its boundary volume $|\del
B_{n,\kappa}(r)|$ are easily computable.

Instead of calculating the isoperimetric profile of a given manifold,
we can look for a sharp isoperimetric inequality in a class of manifolds.
Since $I_{n,\kappa}(V)$ decreases as a function of $\kappa$ for each fixed
$V$, it is natural to consider manifolds whose sectional curvature bounded
above by some $\kappa$.  This motivates the following well-known conjecture.

\begin{conjecture}[Generalized Cartan-Hadamard Conjecture] If $M$ is a
complete, simply connected $n$-manifold with sectional curvature $K$ bounded
above by some $\kappa\le 0$, then every domain $\Omega \subseteq M$ satisfies
\begin{eq}{e:sharp}
|\dOmega | \ge I_{n,\kappa}(|\Omega|).
\end{eq}
\eatline \label{c:main} \end{conjecture}

(If $M$ is not simply connected, then there are many counterexamples.
For example, we can let $M$ be a closed, hyperbolic manifold and let
$\Omega \subseteq M$ be the complement of a small ball.)

The history of Conjecture~\ref{c:main} is as follows
\cite{Osserman:isoperimetric,Druet:curved,Berger:panoramic}.  In 1926,
Weil \cite{Weil:negative} established Conjecture~\ref{c:main} when $n=2$
and $\kappa = 0$ for Riemannian disks $\Omega$, without assuming an ambient
manifold $M$, thus answering a question of Paul L\'evy.  Weil's result
was established independently by Beckenbach and Rad\'o \cite{BR:negative},
who are sometimes credited with the result.  When $n=2$, the case of disks
implies the result for other topologies of $\Omega$ in the presence of $M$.
It was first established by Bol \cite{Bol:inequalities} when $n=2$ and
$\kappa \ne 0$.  Rather later, Conjecture~\ref{c:main} was mentioned
by Aubin \cite{Aubin:sobolev} and Burago-Zalgaller for $\kappa \le 0$
\cite{BZ:inequalities}, and by Gromov \cite{Gromov:struct,Gromov:structen}.
The case $\kappa=0$ is called the Cartan-Hadamard conjecture, because
a complete, simply connected manifold with $K \le 0$ is called a
Cartan-Hadamard manifold.

Soon afterward, Croke proved Conjecture~\ref{c:main} in dimension $n=4$
with $\kappa = 0$ \cite{Croke:sharp}. Kleiner \cite{Kleiner:isoperimetric}
proved Conjecture~\ref{c:main} in dimension $n=3$, for all $\kappa \le 0$,
by a completely different method. (See also Ritor\'e \cite{Ritore:optimal}.)
Morgan and Johnson \cite{MJ:sharp} established Conjecture~\ref{c:main}
for \emph{small} domains (see also Druet \cite{Druet:sharp} where the
curvature hypothesis is on \emph{scalar} curvature); however their argument
does not yield any explicit size condition.

Actually, Croke does not assume an ambient Cartan-Hadamard manifold $M$,
only the more general hypothesis that $\Omega$ has unique geodesics.
We believe that the hypotheses of Conjecture~\ref{c:main} are negotiable,
and it has some generalization to $\kappa > 0$.  But the conjecture is not
as flexible as one might think; in particular, Conjecture~\ref{c:main}
is false for Riemannian 3-balls.  (See Theorem~\ref{th:ball} below and
Section~\ref{s:topology}.)  With this in mind, we propose the following.

\begin{conjecture} If $\Omega$ is a manifold with boundary with unique
geodesics, if its sectional curvature is bounded above by some $\kappa >
0$, and if $|\Omega|\le |X_{n,\kappa}|/2$, then
$|\dOmega | \ge I_{n,\kappa}(|\Omega|)$.
\label{c:positive} \end{conjecture}

The volume restriction in Conjecture~\ref{c:positive} is justified for
two reasons.   First, the comparison ball in $X_{n,\kappa}$ only has
unique geodesics when $|\Omega| < |X_{n,\kappa}|/2$.  Second, Croke
\cite{Croke:some} proved a curvature-free inequality, using only
the condition of unique geodesics, that implies a sharp extension
of Conjecture~\ref{c:positive} when $|\Omega| \ge |X_{n,\kappa}|/2$
(Theorem~\ref{th:croke2}).

Of course, one can extend Conjecture~\ref{c:positive} to negative
curvature bounds (and then the volume condition is vacuous).  The resulting
statement is strictly stronger than Conjecture~\ref{c:main}, since every
domain in a Cartan-Hada\-mard manifold has unique geodesics, but there are
unique-geodesic manifolds that cannot embed in a Cartan-Hadamard manifold
of the same dimension (Figure~\ref{f:cone}).

Another type of generalization of Conjecture~\ref{c:main} is one that
assumes a bound on some other type of curvature.   For example,
Gromov \cite[Rem. 6.28$\frac12$]{Gromov:struct} suggests that
Conjecture~\ref{c:main} still holds when $K \le \kappa$ is replaced by
\begin{eq}{e:mixed}
K \le 0, \quad \Ric \le (n-1)^2 \kappa g.
\end{eq}
In fact, Gromov's formulation is ambiguous: He considers $\kappa=-1$ and
writes $\operatorname{Ricci} \le -(n-1)$, which could mean either $\Ric
\le -(n-1)^2 g$ or $\Ric \le -(n-1) g$. The latter inequality is false for
complex hyperbolic spaces $\CH^n$.  The former is similar to our root-Ricci
curvature condition; see below.

Meanwhile Croke \cite{Croke:sharp} only uses a non-local condition that
we call $\Candle(0)$ rather than the curvature condition $K \le 0$; we
state this as Theorem~\ref{th:croke}.

Our previous work \cite{K:gunther} subsumes both of these two
generalizations. More precisely, most of our results will be stated
in terms of two volume comparison conditions, $\Candle(\kappa)$ and
$\LCD(\kappa)$; see Section~\ref{s:candles} for their definitions. One
can interpret our two main results below (in weakened form) without
referring to Section~\ref{s:candles} by replacing $\Candle(\kappa)$ and
$\LCD(\kappa)$ by $K \le \kappa$, since $K \le \kappa \implies \LCD(\kappa)$
is G\"unther's inequality \cite{Gunther:volume,BC:book}, while $\LCD(\kappa)
\implies \Candle(\kappa)$ is elementary.  When $\kappa\le 0$, one can
also replace the $\Candle(\kappa)$ and $\LCD(\kappa)$ hypotheses by the
mixed curvature bounds \eqref{e:mixed}. In \cite{K:gunther} we introduced
a general curvature bound on what we call the \emph{root-Ricci curvature}
$\RRic$, which is more general than both $K\le \kappa$ and \eqref{e:mixed},
and we proved that this bound implies $LCD(\kappa)$ and $\Candle(\kappa)$.

\subsection{Main results}
\label{s:results}

For simplicity, we will consider isoperimetric inequalities only for compact,
smooth Riemannian manifolds $\Omega$ with smooth boundary $\dOmega$;
or for compact, smooth domains $\Omega$ in Riemannian manifolds $M$.
Our constructions will directly establish inequalities for all such $\Omega$.
We therefore don't have to assume a minimizer or prove that one exists.
Our results automatically extend to any limit of smooth objects in a
topology in which volume and boundary volume vary continuously, \eg, to
domains with piecewise smooth boundary.  Note that our uniqueness result,
Theorem~\ref{th:equality}, does not automatically generalize to a limit
of smooth objects; but its proof might well generalize to some limits of
this type.

Our two strongest results are in the next two subsections.   They both
include Croke's theorem in dimension $n=4$ as a special case.  Each theorem
has a volume restriction that we can take to be vacuous when $\kappa = 0$.

\subsubsection{The positive case}

\begin{theorem} Let $\Omega$ be a compact Riemannian manifold with boundary,
of dimension $n \in \{2,4\}$.  Suppose that $\Omega$ has unique geodesics
and is $\Candle(\kappa)$ with $\kappa \ge 0$ (\eg, $K \le \kappa$), and that
$|\Omega|\le |X_{n,\kappa}|/2$.  Then $|\dOmega|\ge I_{n,\kappa}(|\Omega|)$.
\label{th:positive} \end{theorem}

This is our fully general version of Theorem~\ref{th:first}.
As mentioned, Theorem~\ref{th:croke2} provides an optimal
extension of Theorem~\ref{th:positive} to the case $|\Omega| \ge
|X_{n,\kappa}|/2$. Observe that the volume condition is vacuous when
$\kappa=0$, so that Theorem \ref{th:positive} implies Croke's theorem
\ref{th:croke}.

\subsubsection{The negative case}

When $\kappa$ is negative and $n=4$, we only get a partial result.
(But see Section~\ref{s:questions}.)  To state it, we let $r_{n,\kappa}(V)$
be the radius of a ball of volume $V$ in $X_{n,\kappa}$.  We define
$\chord(\Omega)$ to be the length of the longest geodesic in $\Omega$;
we have the elementary inequality
\[ \chord(\Omega) \le \diam(\Omega). \]

\begin{theorem}
Let $M$ be a Cartan-Hadamard manifold of dimension $n \in \{2,4\}$ which
is $\LCD(\kappa)$ with $\kappa \le 0$ (\eg, $K \le \kappa$).  Let $\Omega$
be a domain in $M$, and if $n=4$, suppose that
\begin{eq}{e:small}
\tanh(\chord(\Omega)\sqrt{-\kappa})\tanh(r_{n,\kappa}(|\Omega|)\sqrt{-\kappa})
    \le \frac12.
\end{eq}
Then $|\dOmega| \ge I_{n,\kappa}(|\Omega|)$.
\label{th:negative} \end{theorem}

Actually, Theorem~\ref{th:negative} only needs $M$ to be convex
with unique geodesics rather than Cartan-Hadamard.  However, we
do not know whether that is a more general hypothesis for $\Omega$.
(See Section~\ref{s:topology}.) Observe that \eqref{e:small} is vacuous
when $\kappa=0$, and thus Theorem \ref{th:negative} also implies Croke's
theorem \ref{th:croke}.

The smallness condition \eqref{e:small} means that Theorem~\ref{th:negative}
is only a partial solution to Conjecture~\ref{c:main} when $n=4$.  Note that
since $\tanh(x) < 1$ for all $x$, it suffices that \emph{either} the chord
length \emph{or} the volume of $\Omega$ is small.  \Ie, it suffices that
\[ \sqrt{-\kappa} \min\left(\chord(\Omega),r_{n,\kappa}(|\Omega|)\right)
    \le \arctanh(\frac12) = \frac{\log(3)}{2}. \]
If we think of Conjecture~\ref{c:main} as parametrized by dimension,
volume, and the curvature bound $\kappa$, then Theorem~\ref{th:negative}
is a complete solution for a range of values of the parameters.

\subsubsection{Pointwise illumination}

We prove a pointwise inequality which, in dimension $2$, generalizes
Weil's isoperimetric inequality \cite{Weil:negative}.  We state it in terms
of \emph{illumination} of the boundary of a domain $\Omega$ by light sources 
lying in $\Omega$, defined rigorously in Section~\ref{s:prince}.

\begin{theorem} Let $\Omega$ be a compact Riemannian $n$-manifold with
boundary, with unique geodesics, and which is $\Candle(0)$; and let $p
\in \dOmega$.  If we fix the volume $|\Omega|$, then the illumination
of $p$ by a uniform light source in $\Omega$ is maximized when $\Omega$
is Euclidean and is given by the polar relation
\begin{eq}{e:polar}
r \le k\cos(\theta)^{1/(n-1)}
\end{eq}
for some constant $k$, with $p$ at the origin.  In particular, in dimension
$n=2$, the optimum $\Omega$ is a round disk.
\label{th:point} \end{theorem}

Theorem~\ref{th:point} generalizes the elementary Proposition~\ref{p:prince},
the problem of the Little Prince, which was part of the inspiration for
the present work.

When $n=2$, Theorem~\ref{th:point} shows that a Euclidean, round disk
maximizes illumination simultaneously at all points of its boundary, and
therefore maximizes the \emph{average} illumination over the boundary.  But,
as a consequence of the divergence theorem, the \emph{total} illumination
over the boundary is proportional to $|\Omega|$. A Euclidean, round disk must
therefore minimize $|\dOmega|$, which is precisely Weil's theorem.

\subsubsection{Equality cases}

We also characterize the equality cases in Theorems \ref{th:positive}
and \ref{th:negative}, with a moderate weakening when $\kappa = 0$.

\begin{theorem} Suppose that $\Omega$ is optimal in Theorem~\ref{th:positive}
or \ref{th:negative}, again with $n\in\{2,4\}$.  When $\kappa=0$, suppose
further that $\Omega$ is $\RRic$ class $0$.  Then $\Omega$ is isometric
to a metric ball in $X_{n,\kappa}$.
\label{th:equality} \end{theorem}

Again, see Section~\ref{s:candles} for the definition of root-Ricci
curvature $\RRic$.  In particular, $\RRic$ class 0 implies $\Candle(0)$,
but it does not imply $K \le 0$ when $n > 2$.

We will prove Theorem~\ref{th:equality} in Section~\ref{s:uniq};
see also Section~\ref{s:questions}.

\subsubsection{Relative inequalities and multiple images}

Choe \cite{Choe:relative2,Choe:double} generalizes Weil's and Croke's
theorems in dimensions 2 and 4 to a domain $\Omega \subseteq M$ which
is outside of a convex domain $C$, which is allowed to share part of its
boundary with $\del C$; he then minimizes the boundary volume $|\dOmega
\setminus \del C|$.  The optimum in both cases is half of a Euclidean ball.

Choe's method in dimension 4 is to consider reflecting geodesics that
reflect from $\del C$ like light rays.  (This dynamic is also called
\emph{billiards}, but we use optics as our principal metaphor.)  Such an
$\Omega$ cannot have unique reflecting geodesics; rather two points in
$\Omega$ are connected by at most two geodesics.  We generalize Choe's result
by bounding the number of connecting geodesics by any positive integer.

\begin{theorem} Let $\Omega$ be a compact $n$-manifold with boundary
with $n=2$ or $4$, let $\kappa \ge 0$, and let $W \subset \dOmega$ be a
(possibly empty) $(n-1)$-dimensional submanifold.  Suppose that $\Omega$
is $\Candle(\kappa)$ for geodesics that reflect from $W$ as a mirror, and
suppose that every pair of points in $\Omega$ can be linked by at most $m$
(possibly reflecting) geodesics.  Suppose also that
\[ |\Omega| \le \frac{|X_{n,\kappa}|}{2m}. \]
Then
\begin{eq}{e:multiple}
|\dOmega \setminus \dW| \ge \frac{I_{n,\kappa}(m |\Omega|)}{m}.
\end{eq}
\eatline \label{th:multiple} \end{theorem}

Note that G\"unther's inequality generalizes to this case (Proposition
\ref{p:mgunther}):  If $\Omega$ satisfies $K \le \kappa$, and if the
mirror region $W$ is locally concave, then $\Omega$ is $\LCD(\kappa)$
and therefore $\Candle(\kappa)$ for reflecting geodesics.

Theorem \ref{th:multiple} is sharp, as can be seen from various examples.
Let $G$ be a finite group that acts on the ball $B_{n,\kappa}(r)$
by isometries.  Then the orbifold quotient $\Omega = B_{n,\kappa}(r)/G$
matches the bound of Theorem~\ref{th:multiple}, if we take the reflection
walls to be mirrors and if we take $m = |G|$.   Although $\Omega$ has
lower-dimensional strata where it fails to be a smooth manifold, we can
remove thin neighborhoods of them and smooth all ridges to make a manifold
with nearly the same volume and boundary volume.

We could state a version of Theorem~\ref{th:multiple} for $\kappa < 0$ using
the $\LCD(\kappa)$ condition, but it would be much more restricted because
we would require an ambient $M$ in which every two points are connected
by \emph{exactly} $m$ geodesics. We do not know any interesting example
of such an $M$.  (\Eg, if the boundary of $M$ is totally geodesic, then
this case is equivalent to simply doubling $M$ and $\Omega$ across $\dM$.)

\subsubsection{Counterexamples in dimension 3}

We find counterexamples to justify the hypotheses of an ambient
Cartan-Hadamard manifold and unique geodesics in Conjectures \ref{c:main}
and \ref{c:positive}. One might like to replace these geometric hypotheses
by purely topological ones.   However, we show that even if $\Omega$
is diffeomorphic to a ball, this does not imply any isoperimetric inequality.

\begin{theorem} For every $\eps > 0$, there is a Riemannian 3-manifold
$\Omega \cong B^3$ with $(1-\eps)$-pinched negative curvature and
with arbitrarily large volume $|\Omega|$ and bounded surface
area $|\dOmega|$ (depending only on $\eps)$.
\label{th:ball} \end{theorem}

Recall that a Riemannian manifold has \emph{\emph{$\delta$}-pinched negative
curvature} when its sectional curvature $K$ satisfies $-1 \le K \le -\delta$
everywhere.

While the manifold $\Omega$ we construct in Theorem~\ref{th:ball} has
trivial topology, its geometry is decidedly non-trivial.  Most of its volume
consists of a truncated hyperbolic knot complement $S_3 \setminus J$ with
constant curvature $K = \kappa \approx -1$.   Such an $\Omega$ has closed
geodesics, which strongly contradicts the property of unique geodesics.
Informally, we can say that $\Omega$ is ``a 3-ball that wants to be a
hyperbolic knot complement".

Theorem~\ref{th:ball} was inspired by Joel Hass' construction of a
negatively curved ball with totally concave boundary \cite{Hass:concave}.
Both constructions yield counterexamples to a Riemannian extension problem
considered by Pigola and Veronelli \cite{PV:extension}.  In both cases,
the ball $\Omega$ has closed geodesics; if $\Omega$ could extend to a
complete manifold $W$ that satisfies $K\le -1+\eps$ or even $K\le 0$, then
its univeral cover $M = \tilde{W}$ would be a Cartan-Hadamard manifold with
closed geodesics, contradicting the Cartan-Hadamard theorem.  The ultimate
purpose of either construction also obstructs a Cartan-Hadamard extension.
In Hass' case, a compact $\Omega$ in a Cartan-Hadamard manifold cannot
have totally concave boundary.   In our case, by Kleiner's isoperimetric
inequality \cite{Kleiner:isoperimetric}, $\Omega$ cannot have arbitrarily
large volume and bounded surface area.

It is not hard to convert the result of Theorem~\ref{th:ball} to a complete
refutation of any possible isoperimetric relation for negatively curved
3-balls.

\begin{corollary} For each $V, A > 0$, there is a Riemannian 3-ball $\Omega$
with $K \le -1$ and with $|\Omega| = V$ and $|\dOmega| = A$.
\label{c:ball} \end{corollary}

We sketch the proof of Corollary~\ref{c:ball}:  Starting with $|\Omega|
\gg V$ and $|\dOmega|$ bounded, we can rescale $\Omega$ to make $|\Omega|
= V-\eps$ and $|\dOmega| < A$.  We can then increase $|\dOmega|$ while
increasing $|\Omega|$ by an arbitrarily small constant by adding a long,
thin finger to $\Omega$.  Pinched negative curvature is an interesting extra
property.  We do not know whether one can achieve $|\dOmega| \to 0$ with
$|\Omega|$ bounded below, and with $(1-\eps)$-pinched negative curvature.

\subsection{The linear programming model}
\label{s:model}

Our method to prove Theorems~\ref{th:positive} and \ref{th:negative}
(and indirectly Theorem~\ref{th:point}) is a reinterpretation and
generalization of Croke's argument, based on optical transport, optimal
transport, and linear programming.

We simplify our manifold $\Omega$ to a measure $\mu_\Omega$ on the set
of triples $(\ell,\alpha,\beta)$, where $\ell$ is length of a complete
geodesic $\gamma \subseteq \Omega$ and $\alpha$ and $\beta$ are its boundary
angles.  Thus $\mu_\Omega$ is always a measure on the set $\R_{\ge 0} \times
[0,\pi/2)^2$, regardless of the geometry or even the dimension of $\Omega$.
We then establish a set of linear constraints on $\mu_\Omega$, by combining
Theorem \ref{th:santalo} (more precisely equations~\eqref{e:etendue} and
\eqref{e:margin}) with Lemmas \ref{l:core}, \ref{l:extend}, and \ref{l:cbk}.
The result is the basic LP Problem \ref{lp:basic} and an extension
\ref{lp:extend}.  The constraints of the model depend on the volume $V =
|\Omega|$ and the boundary volume $A = |\dOmega|$, among other parameters.

Given such a linear programming model, we can ask for which pairs $(V,A)$
the model is \emph{feasible}; \ie, does there exist a measure $\mu$ that
satisfies the constraints?  On the one hand, this is a vastly simpler
problem than the original Conjecture~\ref{c:main}, an optimization over all
possible domains $\Omega$.  On the other hand, the isoperimetric problem,
minimizing $A$ for any fixed $V$, becomes an interesting question in its
own right in the linear model.

Regarding the first point, finite linear programming is entirely algorithmic:
It can be solved in practice, and provably in polynomial time in general.
Our linear programming models are infinite, which is more complicated
and should technically be called convex programming. Nonetheless,
each model has the special structure of optimal transport problems,
with finitely many extra parameters.  Optimal transport is even nicer
than general linear programming.   All of our models are algorithmic in
principle.  In fact, our proofs of optimality in the two most difficult
cases are computer-assisted using Sage~\cite{Sage}.

Regarding the second point, our model is successful in two different
ways:  First, even though it is a relaxation, it sometimes yields a sharp
isoperimetric inequality, \ie, Theorems \ref{th:positive}, \ref{th:negative},
and \ref{th:multiple}.  Second, our models subsume several previously
published isoperimetric inequalities.  We mention six significant ones.
Note that the first four, Theorems \ref{th:weil}-\ref{th:choe}, are special
cases of Theorems \ref{th:positive}, \ref{th:negative}, and \ref{th:multiple}
as mentioned in Section \ref{s:results}.  The other two results are separate,
but they also hold in our linear programming models.

\begin{theorem}[Variation of Weil \cite{Weil:negative} and Bol
\cite{Bol:inequalities}] Let $\Omega$ be a compact Riemannian surface with
curvature $K \le \kappa \ge 0$ with unique geodesics, and suppose that
$\kappa|\Omega| \le 2\pi$.  Then for fixed area $|\Omega|$, the
perimeter $|\dOmega|$ is minimized when $|\Omega|$ has constant curvature
$K = \kappa$ and is a geodesic ball.
\label{th:weil} \end{theorem}

\begin{theorem}[Variation of Bol \cite{Bol:inequalities}] Suppose that
$\Omega \subseteq M$ is a domain in a Cartan-Hadamard surface $M$ that
satisfies $K \le \kappa \le 0$.  Then for fixed area $|\Omega|$, the
perimeter $|\dOmega|$ is minimized when $|\Omega|$ has constant curvature
$K = \kappa$ and is a geodesic ball.
\label{th:bol} \end{theorem}

\begin{theorem}[Croke \cite{Croke:sharp}] If $\Omega$ is a compact
4-manifold with boundary, with unique geodesics, and which is $\Candle(0)$,
then for each fixed volume $|\Omega|$, the boundary volume $|\dOmega|$
is minimized when $\Omega$ is a Euclidean geodesic ball.
\label{th:croke} \end{theorem}

\begin{theorem}[Choe \cite{Choe:relative2,Choe:double}] Let $M$
be a Cartan-Hadamard manifold of dimension $n \in \{2,4\}$,
and let $\Omega \subseteq M$ be a domain whose interior is disjoint
from a convex domain $C \subseteq M$. Then
\[ |\dOmega \setminus \del C| \ge \frac{I_{n,0}(2|\Omega|)}2. \]
\eatline \label{th:choe} \end{theorem}

\begin{theorem}[Croke \cite{Croke:some}] If $\Omega$ is an $n$-manifold with
boundary with unique geodesics, then $|\dOmega| \ge |\del Y_{n,\rho}|$ where
$Y_{n,\rho}$ is a hemisphere with constant curvature $\rho$ and $\rho$ is
chosen so that $|\Omega| = |Y_{n,\rho}|$.
\label{th:croke2} \end{theorem}

Note that when $|\Omega| \ge |X_{n,\kappa}|/2$, we obtain $\rho \le \kappa$,
so that Croke's inequality extends Theorem~\ref{th:positive}, as promised.
See the end of Section~\ref{s:croke2} for further remarks about this result.

\begin{theorem}[Yau \cite{Yau:isoperimetric}] Let $M$ be a Cartan-Hadamard
$n$-manifold which is $\LCD(\kappa)$ with $\kappa < 0$.  Then every domain
$\Omega \subseteq M$ satisfies
\[ |\dOmega| \ge (n-1)\sqrt{-\kappa}|\Omega|. \]
\eatline \label{th:yau} \end{theorem}

Finally, we state the result that our models subsume all of these bounds.

\begin{theorem} Let $\mu$ be a measure on $\mathbb{R}^+ \times [0,\pi/2)^2$
that satisfies LP Problem \ref{lp:basic}, with formal dimension
$n$, formal curvature bound $\kappa$, formal volume $V(\mu)$ (defined
by \eqref{e:lpvol}), and formal boundary volume $A(\mu)$ (defined by
\eqref{e:lpmargin}). Then $\mu$ satisfies Theorem \ref{th:positive} and
therefore \ref{th:weil}.  If $\mu$ satisfies LP Problem \ref{lp:extend},
then it satisfies Theorems~\ref{th:yau} and \ref{th:negative}, and
therefore \ref{th:bol}.  If $\mu$ satisfies LP Problem \ref{lp:croke2},
then it satisfies \ref{th:croke2}.  If $\mu$ satisfies the LP model
\ref{lp:multiple}, then it satisfies Theorem~\ref{th:multiple} and therefore
\ref{th:choe}.
\label{th:subsume} \end{theorem}

We will prove some cases of Theorem~\ref{th:subsume} in the course
of proving our other results; the remaining cases will be done in
Section~\ref{s:oldwine}.

Our linear programming models are similar to the important Delsarte linear
programming method in the theory of error{\hyp}correcting codes and sphere
packings \cite{Delsarte:linear,CS:splag,CE:new}.  Delsarte's original
result was that many previously known bounds for error{\hyp}correcting codes
are subsumed by a linear programming model.  But his model also implies
new bounds, including sharp bounds.  For example, consider the kissing
number problem for a sphere in $n$ Euclidean dimensions \cite{CS:splag}.
The geometric maximum is of course an integer, but in a linear programming
model this may no longer be true.   Nonetheless, Odlyzko and Sloane
\cite{OS:unit} established a sharp geometric bound in the Delsarte model,
which happens to be an integer and the correct one, in dimensions 2, 8,
and 24.  (The bounds are, respectively, 6, 240, and 196,560 kissing spheres.)
The basic Delsarte bound for the sphere kissing problem is quite strong in
other dimensions, but it is not usually an integer and not usually sharp
even if rounded down to an integer.

Another interesting common feature of the Delsarte method and ours is
that they are both sets of linear constraints satisfied by a two-point
correlation function, \ie, a measure derived from taking pairs of points
in the geometry.

\subsection{Other results}

\subsubsection{Croke in all dimensions}

There is a natural version of Croke's theorem in all dimensions.  This is
a generalized, sharp isoperimetric inequality in which the volume of
$\Omega$ is replaced by some other functional when the dimension $n \ne 4$.
This result might not really be new; we state it here to further illustrate
of our linear programming model.

If $\Omega$ is a manifold with boundary and unique geodesics, then the
space $G$ of geodesic chords in $\Omega$ carries a natural measure $\mu_G$,
called Liouville measure or \'etendue (Section~\ref{s:geodesics}).

\begin{theorem} Let $\Omega$ be a compact manifold with boundary of dimension
$n \ge 4$, with unique geodesics, and with non-positive sectional curvature.
Let
\[ L(\Omega) = \int_{G} \ell(\gamma)^{n-3} \sdd\mu_G(\gamma)\]
If $B_{n,0}(r)$ is the round, Euclidean ball such that
\[ L(\Omega) = L(B_{n,0}(r)), \]
then
\[ |\dOmega| \ge |\del B_{n,0}(r)|. \]
\eatline \label{th:alldim} \end{theorem}

By Theorem~\ref{th:santalo} (Santal\'o's equality),
\[ \omega_{n-1}|\Omega| = \int_G \ell(\gamma) \sdd\mu_G(\gamma). \]
(Here $\omega_n = |X_{n,1}|$ is the $n$-sphere volume; see
Section~\ref{s:basic}.)  Thus Theorem~\ref{th:alldim} is Croke's Theorem
if $n=4$.   The theorem is plainly a sharp isoperimetric bound for the
boundary volume $|\dOmega|$ in all cases given the value of $L(\Omega)$,
which happens to be proportional to the volume $|\Omega|$ only when $n = 4$.

Similar results are possible with a curvature bound $K < \kappa$,
only with more complicated integrands $F(\ell)$ over the space $G$.

\subsubsection{Non-sharp bounds and future work}
\label{s:nonsharp}

We mention three cases in which the methods of this paper yield improved
non-sharp results.

First, when $n=3$ and $\kappa=0$, Problem~\ref{lp:basic} yields a non-sharp
version of Kleiner's theorem under the weaker hypotheses of $\Candle(0)$ and
unique geodesics.  Croke \cite{Croke:sharp} established the isoperimetric
inequality in this case up to a factor of $\sqrt[3]{36/32} = 1.040\ldots$.
Meanwhile Theorem~\ref{th:point} implies the same isoperimetric inequality
up to a factor of $\sqrt[3]{27/25} = 1.026\ldots$.  The wrinkle is that
Croke's proof uses only \eqref{e:lpcroke}, while Theorem~\ref{th:point}
uses only \eqref{e:lpprince}.  A combined linear programming problem should
produce a superior if still non-sharp bound.

Second, it is a well-known conjecture that a metric ball is the unique
optimum to the isoperimetric problem for domains in the complex hyperbolic
plane $\CH^2$.  (The same conjecture is proposed for any non-positively
curved symmetric space of rank 1.)  If we normalize the metric on $\CH^2$
so that it is $(-4,-1)$-pinched (Section~\ref{s:basic}), then $\CH^2$
is $\LCD(-16/9)$.  Then Theorem~\ref{th:negative} is, to our knowledge,
better than what was previously established for moderately small volumes.
Even so, this is a crude bound; we could do even better with a version of
Problem~\ref{lp:extend} that uses the specific candle function of $\CH^2$.

Third, even for domains in Cartan-Hadamard manifolds with $K\le-1$ (or more
generally $\LCD(-1)$), we can relax the smallness condition \eqref{e:small}
in Theorem~\ref{th:negative} simply by increasing the curvature bound
$\kappa$ from $\kappa = -1$.  This is still a good bound for a range
of volumes until it is eventually surpassed by Theorem~\ref{th:yau}.
This is also a crude bound that can surely be improved, given that both
Theorem~\ref{th:negative} and Theorem~\ref{th:yau} hold in the same linear
programming model, Problem \ref{lp:extend}.

\acknowledgments

The authors would like to thank Sylvain Gallot, Joel Hass, Misha Kapovich,
and Qinglan Xia for useful discussions about Riemannian geometry and
optimal transport.  The authors would also like to thank an anonymous
referee for detailed corrections and comments.

\section{Conventions}

\subsection{Basic conventions}
\label{s:basic}

If $f:\R_{\ge 0} \to \R$ is an integrable function, we let
\[ f^{(-1)}(x) \defeq \int_0^x f(t) \sdd t \]
be its antiderivative that vanishes at $0$, and then by induction its
$n$th antiderivative $f^{(-n)}$.   This is in keeping with the standard
notation that $f^{(n)}$ is the $n$th derivative of $f$ for $n > 0$.

If $M$ is a Riemannian manifold, we let $\nu_M$ denote the Riemannian
measure on $M$.  As usual, $TM$ is the tangent bundle of $M$, while we use
$UM$ to denote the unit tangent bundle.  Also, if $\Omega$ is a manifold
with boundary $\dOmega$, then we let
\[ U^+ \dOmega \defeq \big\{u=(p,v) \st p\in \dOmega, v\in U_p\Omega
    \mbox{ inward pointing}\big\}. \]

We say that $M$ has \emph{$(\delta_1,\delta_2)$-pinched curvature} if its
sectional curvature $K$ satisfies $\delta_1 \le K \le \delta_2$ everywhere.
To paraphrase, we may say that $M$ is pinched, its metric is pinched, etc.

We let $|M|$ be the volume of $M$:
\[ |M| \defeq \int_M \dd\nu_M. \]
We let
\[ \omega_n = |X_{n,1}| = \frac{2\pi^{(n+1)/2}}{\Gamma(\frac{n+1}2)} \]
be the volume of the unit $n$-dimensional sphere $X_{n,1} = S^n \subseteq
\R^{n+1}$.

\subsection{Candles}
\label{s:candles}

Our main results are stated in terms of conditions $\Candle(\kappa)$ and
$\LCD(\kappa)$ that follow from the sectional curvature condition $K \le
\kappa$ by G\"unther's comparison theorem \cite{Gunther:volume,BC:book}.
These conditions are non-local, but in previous work \cite{K:gunther},
we showed that they follow from another local condition, more general than
$K \le \kappa$ that we called $\RRic$ class $(\rho,\kappa)$.  The original
motivation is that Croke's theorem only needs that the manifold $D$ is
$\Candle(0)$, and even then only for pairs of boundary points.  Informally,
a Riemannian manifold $M$ is $\Candle(\kappa)$ if a candle at any given
distance $r$ from an observer is dimmer than it would be at distance $r$
in a geometry of constant curvature $\kappa$.

More rigorously, let $M$ be a Riemannian manifold and let $\gamma=\gamma_u$
be a geodesic in $M$ that begins at $p = \gamma(0)$ with initial velocity
$u\in U_pM$.  Then the \emph{candle function} $j_M(\gamma,r)$ of $M$ is
by definition the normalized Jacobian of the exponential map
\[ u \mapsto \gamma_u(r) = \exp_p(ru), \]
given by the equation
\[ \dd\nu_M(\gamma_u(r)) = j_M(\gamma_u,r) \sdd\nu_{U_pM}(u) \sdd r \]
for $r > 0$, where $\nu_M$ is the Riemannian volume on $M$ and $\nu_{U_pM}$
is the Riemannian measure on the round unit sphere $U_pM$.  More generally,
if $a < b$, we define
\[ j_M(\gamma,a,b) = j_M(\gamma_a,b-a), \]
where $\gamma_a$ is the same geodesic as $\gamma$ but with parameter
shifted by $a$.  We also define
\[ j_M(\gamma,b,a) = j_M(\overline{\gamma_b},b-a), \]
where $\overline{\gamma_b}$ is the same geodesic as $\gamma$, but reversed
and based at $\gamma(b)$. (But see Corollary~\ref{c:iseeyou}.)

The candle function of the constant-curvature geometry $X_{n,\kappa}$ is
independent of the geodesic.  We denote it by $s_{n,\kappa}(r)$; it is
given by the following explicit formulas:
\begin{eq}{e:ccandle}
s_{n,\kappa}(r) = \begin{cases} {\displaystyle
    \Big(\frac{\sin(r\sqrt{\kappa})} {\sqrt{\kappa}}\Big)^{n-1}}
    & \mbox{if $\kappa > 0$, $\displaystyle r\le\frac{\pi}{\sqrt{\kappa}}$} \\
    r^{n-1} & \mbox{if $\kappa = 0$} \\
    {\displaystyle \Big(\frac{\sinh(r\sqrt{-\kappa})}
    {\sqrt{-\kappa}}\Big)^{n-1}} & \mbox{if $\kappa < 0$.} \end{cases}
\end{eq}
We will also need the extension $s_{n,\kappa}(r) = 0$ when
$\kappa > 0$ and $r \ge \pi/\sqrt{\kappa}$.

\begin{definition} An $n$-manifold $M$ is $\Candle(\kappa)$ if
\[ j_M(\gamma,r) \ge s_{n,\kappa}(r) \]
for all $\gamma$ and $r$.  It is $\LCD(\kappa)$, for \emph{logarithmic candle
derivative}, if
\[ \log(j_M(\gamma,r))' \ge \log(s_{n,\kappa}(r))' \]
for all $\gamma$ and $r$.  (Here the derivative is with respect to $r$.)
The $\LCD(\kappa)$ condition implies the $\Candle(\kappa)$ condition by
integration.  If $\kappa > 0$, then these conditions are only required
up to the focal distance $\pi/\sqrt{\kappa}$ in the comparison geometry.
\end{definition}

To illustrate how $\Candle(\kappa)$ is more general than $K \le \kappa$,
we mention \emph{root-Ricci curvature} \cite{K:gunther}.  Suppose that
$M$ is a manifold such that $K\le 0$ and let $\kappa < 0$.   For any unit
tangent vector $u \in U_pM$ with $p \in M$, we define
\[ \RRic(u) \defeq \Tr(\sqrt{-R(\cdot,u,\cdot,u)}). \]
Here $R(u,v,w,x)$ is the Riemann curvature tensor expressed as a tetralinear
form, and the square root is the positive square root of a positive
semidefinite matrix or operator.  We say that $M$ \emph{is of $\RRic$
class $\kappa$} if $K\le 0$ and
\[ \RRic(u) \ge (n-1)\sqrt{-\kappa}. \]
Then
\[ K\le\kappa \implies \RRic \textrm{ class } \kappa \implies \LCD(\kappa)
    \implies \Candle(\kappa) \implies \Ric \le (n-1) \kappa g. \]
The second implication, from $\RRic$ to $\LCD$, is the main result of
\cite{K:gunther}.  (We also established a version of the result that
applies for any $\kappa \in \R$.  This version uses a generalized $\RRic$
class $(\rho,\kappa)$ condition that also requires $K \le \rho$ for a
constant $\rho > \max(\kappa,0)$.)  All implications are strict when $n>2$.
By contrast in dimension 2, the last condition trivially equals the
first one, so all of the conditions are equivalent.

We conclude with two examples of $4$-manifolds of $\RRic$ class $-1$,
and which are therefore $\LCD(-1)$, but that do not satisfy $K\le -1$:
\begin{itemize}
\item The complex hyperbolic plane, normalized to be
$(-\frac94,-\frac9{16})$-pinched.
\item The product of two simply connected surfaces that each satisfy
$K < -9$.
\end{itemize}
Actually, the most important regime where $\Candle(\kappa)$ is weaker
than $K \le \kappa$ is at short distances.  Since 
\[ j_M(\gamma,r) = r^{n-1} - \frac{\Ric(\gamma'(0),\gamma'(0))}6 r^{n+1}
    + O(r^{n+2}) \]
in dimension $n$, we can write informally that
\[ \Candle(\kappa) \stackrel{\approx}{\iff} \Ric \le (n-1) \kappa g \]
as $\diam(M) \to 0$.

\section{The Little Prince and other stories}
\label{s:prince}

\begin{figure}[htb]
\begin{center}
\frame{\includegraphics[height=3in]{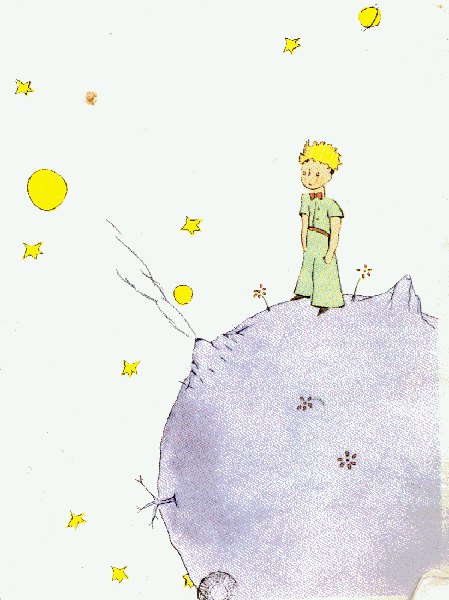}}
\end{center}
\caption{The Little Prince on his not-very-big planet, actually an asteroid.}
\label{f:prince} \end{figure}

\subsection{The problem of the Little Prince}

As Saint-Exup\'ery related to inhabitants of our planet, the
Little Prince lives on his own planet, also known as asteroid B-612
(Figure~\ref{f:prince}).  Since this planet is not very big, its
gravitational pull is small and its habitation is precarious.  The question
arises as to what shape it should be to maximize the normal component of
gravity for the Little Prince, assuming that the planet has a fixed mass,
and a uniform, fixed mass density.  Let $\Omega$ be the shape of the
planet.  The divergence theorem tells us that the average normal gravity
is proportional to $|\Omega|/|\dOmega|$, so maximizing the average would
be exactly the isoperimetric problem.  Suppose instead that the Little
Prince has a favorite location, and does not mind less gravity elsewhere.
(After all, in the illustrations he usually stands on top of the planet.)

We cannot be sure of the dimension of the Little Prince or his planet.
The illustrations are 2-dimensional, but the Prince visits the Sahara Desert
which suggests that he is 3-dimensional.  In any case higher-dimensional
universes, which are a fashionable topic in physics these days, would each
presumably have their own Little Prince.  So we assume that the Little
Prince is $n$-dimensional for some $n \ge 2$.  We first assume Newtonian
gravity and therefore a Euclidean planet; recall that in $n$ dimensions,
a divergenceless central gravitational force is proportional to $r^{1-n}$.

\begin{proposition}[Little Prince Problem] Let $\Omega$ be the shape of a
planet in $n$ Euclidean dimensions with a pointwise gravitational force
proportional to $r^{1-n}$.  Suppose that the planet has a fixed volume
$|\Omega|$ and a uniform, fixed mass density, and let $p \in \dOmega$.
Then the total normal gravitational force $F(\Omega,p)$ at $p$ is maximized
when $\Omega$ is bounded by the surface $r = k \cos(\theta)^{1/(n-1)}$
for some constant $k$, in spherical coordinates centered at $p$.
\label{p:prince} \end{proposition}

The problem of the Little Prince in 3 dimensions is sometimes used as
an undergraduate physics exercise \cite{McDonald:gravity}.  It has also
been previously used to prove the isoperimetric inequality in 2 dimensions
\cite{HHM:monthly}.  However, our further goal is inequalities for curved
spaces such as Theorem~\ref{th:point}.

\begin{proof} For convenience, we assume that the gravitational constant
and the mass density of the planet are both 1.  Given $x \in \Omega$,
let $r = r(x)$ and $\theta = \theta(x)$ be the radius and first angle
in spherical coordinates with the point $p$ at the origin, and such
that the normal component of gravity is in the direction $\theta = 0$.
Then the total gravitational effect of a volume element $\dd x$ at $x$
is $\cos(\theta)r^{1-n}\sdd x$, so the total gravitational force is
\[ F(\Omega,p) = \int_{x \in \Omega} \cos(\theta)r^{1-n} \sdd x. \]
In general, if $f(x)$ is a continuous function and we want to choose a
region $\Omega$ with fixed volume to maximize
\[ \int_\Omega f(x) \sdd x, \]
then by the ``bathtub principle", $\Omega$ should be bounded by a level
curve of $f$, \ie,
\[ \Omega = f^{-1}([k,\infty)) \]
for some constant $k$.  Our $f$ is not continuous at the origin, but the
principle still applies.  Thus $\Omega$ is bounded by a surface of the form
\[ r = k \cos(\theta)^{1/(n-1)}. \qedhere \]
\end{proof}

As explained above in words, the integral over $\dOmega$ of the normal
component of gravity is proportional to $|\Omega|$ by the divergence theorem.
More rigorously:  We switch to a vector expression for gravitational force
and we do not assume that $p=0$.  Then
\[ F(\Omega,p) = \int_{\Omega} (x-p)|x-p|^{-n} \sdd x. \]
Since for each fixed $x \in \Int(\Omega)$, the vector field $p\mapsto
(x-p)|x-p|^{-n}$ is divergenceless except at its singularity, we have
\[ \int_{\dOmega} \braket{-w(p),x-p}|x-p|^{-n} \sdd p = \omega_{n-1} \]
where $w(p)$ is the outward unit normal vector at $p$. Thus
\[ \int_{\dOmega} \braket{-w(p),F(\Omega,p)} \sdd p
  = \omega_{n-1} |\Omega| \]
by switching integrals.  Then
\begin{eq}{e:isoprince}
\omega_{n-1} |\Omega| \le |\dOmega| F_{\max},
\end{eq}
where $F_{\max}$ is the upper bound established by
Proposition~\ref{p:prince}.

In particular, when $n=2$, the optimum $\Omega$ is the polar plot of $r
= k\cos(\theta)$, which is a round circle.  In this case
\[ \braket{-w(p),F(\Omega,p)} = F_{\max} \]
at all points simultaneously.  Thus when $n=2$, equation~\eqref{e:isoprince}
is exactly the sharp isoperimetric inequality \eqref{e:sharp}.

\subsection{Illumination and Theorem~\ref{th:point}}
\label{s:illum}

Proposition~\ref{p:prince} is close to a special case of
Theorem~\ref{th:point}.  To make it an actual special case, we slightly
change its mathematics and its interpretation, but we will retain the sharp
isoperimetric corollary using the divergence theorem.  Instead of the shape
of a planet, we suppose that $\Omega$ is the shape of a uniformly lit room,
and we let $I(\Omega,p)$ be the total intensity of light at a point on the
wall $p \in \dOmega$.  More rigorously, if $\Vis(\Omega,p)$ is the subset of
$\Omega$ which is visible from $p$ (assuming that the walls are opaque, but 
allowing geodesics to be continued when they meet the boundary tangentially 
so that $\Vis(\Omega,p)$ is closed), then
\[ I(\Omega,p) = \int_{\Vis(\Omega,p)} \back\back
    \braket{-w(p),x-p}|x-p|^{-n} dx. \]
We still have
\[ \int_{\Vis(\dOmega,x)} \back\back
    \braket{-w(p),x-p}|x-p|^{-n} \sdd p = \omega_{n-1} \]
and we can still exchange integrals.  Moreover,
\[ I(\Omega,p) = F(\Omega,p) \]
when $\Omega$ is convex.  Thus, this variation of Proposition~\ref{p:prince}
is also true and also implies \eqref{e:sharp}.

We now consider the case when $\Omega$ is a curved Riemannian manifold,
that is, Theorem~\ref{th:point}.  The proof is a simplified version of the
proof of Theorems \ref{th:positive} and \ref{th:negative}.  Before giving the
proof, we give a rigorous definition of illumination in the curved setting.
(The definition agrees with the natural geometric assumption that light
rays travel along geodesics.)

Let $\Omega$ be a compact Riemannian $n$-manifold with boundary and unique
geodesics. We define a Riemannian analogue of $p\mapsto -(x-p)|x-p|^{-n}$,
changing sign here to match the illumination interpretation.  Namely,
for each fixed $x \in \Int(\Omega)$, we define a tangent vector field $v_x$
as follows.  If $y \in \Vis(\Omega,x)$, then we let $\gamma$ be the geodesic
with $\gamma(0) = x$ and $\gamma(r) = y$, and then let
\[ v_x(y) = \frac{\gamma'(r)}{j_\Omega(\gamma,r)}. \]
If $y \notin \Vis(\Omega,x)$, then we let $v_x(y) = 0$.  The motivation,
as above, is that this formula describes the radiation from a point source
of light at $x$ to the rest of $\Omega$.

We claim that $\div v_x = \omega_{n-1}\delta_x$ in a distributional sense,
where $\delta_x$ is the Dirac measure at $x$, so that we can then use $v_x$
in the divergence theorem.   It is routine to check that this holds at $x$
itself and at any point $y$ where $v_x$ is continuous.  The only delicate
case is when $y \in \del \Vis(\Omega,x) \setminus \dOmega$.  The vector
field $v_x$ is not continuous at these points; however, it is parallel to
$\del \Vis(\Omega,x)$ and thus does not have any singular divergence.

We fix a point $p \in \dOmega$ and again let $w(p)$ be the outward unit
normal vector to $\dOmega$ at $p$. Then the illumination at $p$ is defined by
\[I(\Omega,p) = \int_{\Vis(\Omega,p)}
    \back\back\back \braket{w(p), v_x(p)} \sdd \nu_\Omega(x). \]

\begin{proof}[Proof of Theorem~\ref{th:point}] First, we express
$I(\Omega,p)$ as an integral over $U = U^+_p\dOmega$, the unit inward tangent
vectors at $p$.  Given $u\in U$, let $\ell(u)$ be the length of the maximal
geodesic segment defined by $u$ and let $\alpha(u)$ be the angle of $u$
with the inward normal $-w(p)$.  Then, in polar coordinates we get
\begin{align*}
I(\Omega,p) &= \int_U \int_0^{\ell(u)} \back \cos(\alpha(u))
    \sdd t \sdd \nu_U(u) \\
    &= \int_U \ell(u) \cos(\alpha(u)) \sdd \nu_U(u).
\end{align*}
The first equality expresses the fact that the norm $||v_x(p)||$ is the
reciprocal of the Jacobian of the exponential map from $p$.   In other words,
it is based on an optical symmetry principle (Corollary~\ref{c:iseeyou}):
If two identical candles are at $x$ and $p$, then each one looks exactly
as bright from the position of the other one.

Second, the $\Candle(0)$ hypothesis tells us that
\[ |\Omega| \ge |\Vis(\Omega,p)| = \int_{\Vis(\Omega,p)} \back \back \dd
    \nu_\Omega(x) \ge \int_U \int_0^{\ell(u)} t^{n-1} \sdd t \sdd \nu_U(u), \]
so that
\begin{eq}{e:pprince}
|\Omega|\ge \int_U \frac{\ell(u)^n}{n} \sdd \nu_U(u).
\end{eq}

Third, we apply the linear programming philosophy that will be important
in the rest of the paper.

All of our integrands depend only on $\ell$ and $\alpha$.  Thus we can
summarize all available information by projecting the measure $\nu_U$
to a measure
\[ \sigma_\Omega = (\ell,\alpha)_*(\dd\nu_U) \]
on the space of pairs 
\[ (\ell,\alpha) \in \R_{\ge 0} \times [0,\frac{\pi}2). \]
Then we want to maximize 
\begin{eq}{e:illumI}
I = \int_{\ell,\alpha} \ell \cos(\alpha) \sdd\sigma_\Omega
\end{eq}
subject to the constraint
\begin{eq}{e:illumV}
\int_{\ell,\alpha} \frac{\ell^n}n \sdd\sigma_\Omega \le V.
\end{eq}
We have one other linear piece of information: If we project volume on
the hemisphere $U$ into the angle coordinate $\alpha \in [0,\frac{\pi}2)$,
then the result is
\begin{eq}{e:illumal}
\alpha_*(\dd\sigma_\Omega) = \alpha_*(\dd \nu_U)
    = \omega_{n-2}\sin(\alpha)^{n-2}\sdd\alpha,
\end{eq}
since the latitude on $U$ at angle $\alpha$ is an $(n-2)$-sphere
with radius $\sin(\alpha)$.

We temporarily ignore geometry and maximize \eqref{e:illumI} for an abstract
positive measure $\sigma = \sigma_\Omega$ that satisfies \eqref{e:illumV}
and \eqref{e:illumal}.  To do this, choose $a>0$, and let
\begin{eq}{e:illumax}
f(\alpha) = \sup_{\ell>0} \left(\ell \cos(\alpha) - \frac{a\ell^n}n \right).
\end{eq}
We obtain
\begin{align}
0 &\le \int_{\ell,\alpha} \left(f(\alpha) + \frac{a\ell^n}n
    - \ell\cos(\alpha)\right) \sdd\sigma(\ell,\alpha) \nonumber \\
    & \le \int_0^{\pi/2} \back f(\alpha) \omega_{n-2}\sin(\alpha)^{n-2}
    \sdd\alpha + aV - I.
\label{e:illumend} \end{align}
The integral on the right side of \eqref{e:illumend} is a function of
$a$ only.  Finally \eqref{e:illumend} is an upper bound on $I$, one
that achieves equality if \eqref{e:illumV} is an equality and $\sigma =
\sigma_\Omega$ is supported on the locus
\[ \cos(\alpha) = a\ell^{n-1}, \]
because that is the maximand of \eqref{e:illumax}.  The first condition
tells us that $\Omega$ is Euclidean and visible from $p$.  The second gives
us the polar plot \eqref{e:polar} if we take $k = a^{-1/(n-1)}$.
\end{proof}

\begin{remark} It is illuminating to give an alternate Croke-style end to
the proof of Theorem~\ref{th:point}.  Namely, H\"older's inequality says that
\begin{align*}
I &= \int_{\ell,\alpha} \back \ell \cos(\alpha) \sdd\sigma_\Omega  \\
    &\le \Big(\int_{\ell,\alpha} \back \ell^n \sdd\sigma_\Omega\Big)^{\frac1n}
    \Big(\int_{\ell,\alpha} \back \cos(\alpha)^{\frac{n}{n-1}}
    \sdd\sigma_\Omega \Big)^{\frac{n-1}n} \\
    &\le (nV)^{\frac1n} \Big(\int_0^{\frac{\pi}{2}} \back
    \cos(\alpha)^{\frac{n}{n-1}} \omega_{n-2}\sin(\alpha)^{n-2}
    \sdd\alpha \Big)^{\frac{n-1}n}.
\end{align*}
The last expression depends only on $V$ and $n$, while the inequality is
an equality if \eqref{e:illumV} is an equality, and if
\[ \ell^n \propto \cos(\alpha)^{\frac{n}{n-1}}.\]
The first condition again tells us that $\Omega$ is Euclidean and visible
from $p$; the second one gives us the same promised shape \eqref{e:polar}.

The Croke-style argument looks simpler than our proof of
Theorem~\ref{th:point}, but what was elegance becomes misleading for our
purposes.  For one reason, our use of the auxiliary $f(\alpha)$ amounts
to a proof of this special case of H\"older's inequality.  Thus our
argument is not really different; it is just another way to describe the
linear optimization.  For another, we will see more complicated linear
programming problems in the full generality of Theorems~\ref{th:positive}
and \ref{th:negative} that do not reduce to H\"older's inequality.
\end{remark}

\section{Topology and geodesics}
\label{s:topology}

In this section we will analyze the effect of topology and geodesics on
isoperimetric inequalities.

Weil and Bol established the sharp isoperimetric inequality \eqref{e:sharp}
for Riemannian disks $\Omega$ with curvature $K \le \kappa$, without
assuming an ambient manifold $M$, and for any $\kappa \in \R$.  The cases
$\kappa \le 0$ of the Weil and Bol theorems is equivalent to the $n=2$
case of Conjecture~\ref{c:main} \cite{Druet:curved}.

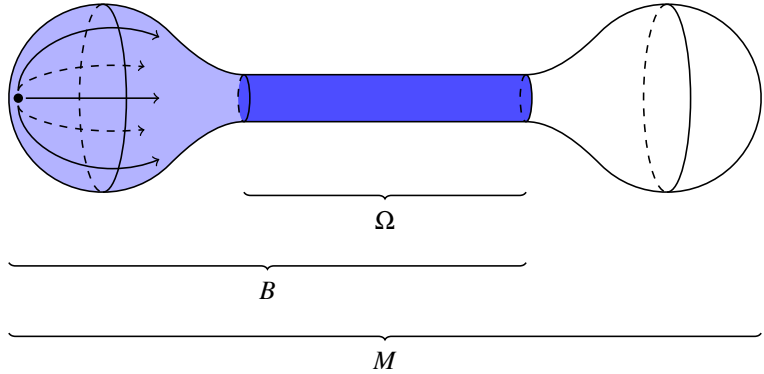
\begin{figure*}[htb]
\begin{center} \begin{tikzpicture}[semithick,scale=1.25]
\fill[color=blue!30!white] (-1.5,.25) .. controls (-1.8,.25) and 
    (-2.093,.507) .. (-2.293,.707) arc (45:315:1) .. controls (-2.093,-.507)
    and (-1.8,-.25) .. (-1.5,-.25) arc (270:90:.0625 and .25) -- cycle;
\fill[color=blue!70!white] (-1.5,.25) arc (90:270:.0625 and .25) -- (1.5,-.25)
   arc (-90:90:.0625 and .25) -- cycle;
\draw (-2.293,.707) arc (45:315:1) .. controls (-2.093,-.507) and
    (-1.8,-.25) .. (-1.5,-.25) -- (1.5,-.25) .. controls (1.8,-.25) and
    (2.093,-.507) .. (2.293,-.707) arc (-135:135:1) .. controls (2.093,.507)
    and (1.8,.25) .. (1.5,.25) -- (-1.5,.25) .. controls (-1.8,.25) and
    (-2.093,.507) .. (-2.293,.707);
\draw[dashed] (-3,1) arc (90:270:.25 and 1);
\draw (-3,-1) arc (-90:90:.25 and 1);
\draw[dashed] (-1.5,.25) arc (90:270:.0625 and .25);
\draw (-1.5,-.25) arc (-90:90:.0625 and .25);
\draw[dashed] (1.5,.25) arc (90:270:.0625 and .25);
\draw (1.5,-.25) arc (-90:90:.0625 and .25);
\draw[dashed] (3,1) arc (90:270:.25 and 1);
\draw (3,-1) arc (-90:90:.25 and 1);
\draw[decorate,decoration=brace] (1.5,-1) -- (-1.5,-1);
\draw (0,-1.1) node[anchor=north] {$\Omega$};
\draw[decorate,decoration=brace] (1.5,-1.75) -- (-4,-1.75);
\draw (-1.25,-1.85) node[anchor=north] {$B$};
\draw[decorate,decoration=brace] (4,-2.5) -- (-4,-2.5);
\draw (0,-2.6) node[anchor=north] {$M$};
\fill (-3.9,0) circle (.05);
\draw[shorten <= .1cm,->] (-3.9,0) arc (180:60:1 and 0);
\draw[shorten <= .1cm,->] (-3.9,0) arc (180:60:1 and .75);
\draw[shorten <= .1cm,->] (-3.9,0) arc (180:300:1 and .75);
\draw[dashed,shorten <= .1cm,shorten >= .2cm,->]
    (-3.9,0) arc (180:60:1 and .35);
\draw[dashed,shorten <= .1cm,shorten >= .2cm,->]
    (-3.9,0) arc (180:300:1 and .35);
\end{tikzpicture} \end{center}
\caption{A counterexample $\Omega \subseteq B \subseteq M$ to Aubin's
    conjecture with $\kappa > 0$, in which both $|\Omega|$ and $|\dOmega|$
    are unrestricted.}
\label{f:barbell} \end{figure*}

The case $\kappa > 0$ is more delicate, even in 2 dimensions.  Aubin
\cite{Aubin:sobolev} assumed that $B$ is a Riemannian ball with $K \le
\kappa$ and then that $\Omega \subseteq B$; but this formulation does
not work.  Even if $B$ is a metric ball with an injective exponential map,
and even if in addition $B \subset M$ and $M$ is complete and simply
connected with the same $K \le \kappa$, there may be no control over
the size of $\dOmega$.  We can let $M$ be a ``barbell" consisting of two
large, nearly round 2-spheres connected by a rod (Figure~\ref{f:barbell}).
Then $\Omega$ can be just the rod, while $B$ is $\Omega$ union one end of
the barbell.  $B$ is also a metric ball with an injective exponential map.
Then $\Omega$ is an annulus $S^1 \times I$ in which both the meridian
$S^1$ and the longitude $I$ can have any length.  Thus both $|\Omega|$
and $|\dOmega|$ can have any value. Morgan and Johnson \cite{MJ:sharp}
made the same point and used it to justify their small-volume hypothesis;
of course, they require an upper bound on the volume that depends on the
geometry of the ambient manifold.

Theorem~\ref{th:ball} says that Weil's theorem fails completely for
negatively curved Riemannian 3-balls.  Our proof is similar to Hass's
construction \cite{Hass:concave} of a negatively curved 3-ball with
concave boundary.

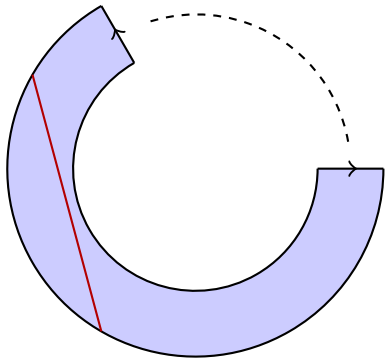
\begin{figure}[htb]
\begin{center} \begin{tikzpicture}[thick,scale=1.25,decoration={markings,
    mark=at position 0.6 with {\arrow[scale=1.5,semithick]{>}}}]
\fill[blue!20!white] (120:2) arc (120:360:2) -- (0:1.3) arc
    (360:120:1.3) -- cycle;
\draw (120:2) arc (120:360:2);
\draw (120:1.3) arc (120:360:1.3);
\draw[dashed] (10:1.65) arc (10:110:1.64);
\draw[postaction={decorate}] (120:1.3)--(120:2);
\draw[postaction={decorate}] (0:1.3)--(0:2);
\draw[darkred] (150:2) -- (240:2);
\end{tikzpicture} \end{center}
\caption{A conical, locally Euclidean annulus that has unique geodesics
    but does not embed in a Cartan-Hadamard surface, with a geodesic
    indicated in red.  (We glue together the edges marked with arrows.)}
\label{f:cone} \end{figure}

If $\Omega$ is a smooth domain in a Cartan-Hadamard manifold as in
Conjecture~\ref{c:main}, then it has unique geodesics, but unique geodesics
is a strictly weaker hypothesis even in 2 dimensions.  For example,
if $\Omega$ is a thin, locally Euclidean annulus with an angle deficit
(Figure~\ref{f:cone}), then it has unique geodesics, but by the Gauss-Bonnet
formula its inner circle cannot be filled by a non-positively-curved disk.
Theorem~\ref{th:ball} tells us that we need some geometric condition on
a manifold $\Omega$ to obtain an isoperimetric inequality, because even
the strictest topological condition, that $\Omega$ be diffeomorphic to a
ball, is not enough.  One natural condition is that $\Omega$ has unique
geodesics. (See also Section~\ref{s:multiple} for a generalization.)

\begin{question} If $M$ is a Cartan-Hadamard manifold and $\Omega$ minimizes
$|\dOmega|$ for some fixed value of $|\Omega|$, then is it convex?  Is it
a topological ball?
\end{question}

Since we first proposed this question in an earlier version of the present article, Hass \cite{Hass:regions} proved that an isoperimetric minimizer
$\Omega$ in a Cartan-Hadamard manifold need not be connected, which is thus
a negative answer to both parts of the question.  He also shows that
in two dimensions, each connected component is a convex disk.  He then
gives partial evidence that nonconvex, connected minimizers exist in
three dimensions.  It may still be interesting to ask what restrictions
there are on the topology of $\Omega$, for simplicity given that it is a
manifold with boundary.

In two dimensions, if $\Omega$ is a non-positively curved disk, then it
has unique geodesics.  (Proof: If a disk does not have unique geodesics,
then it contains a geodesic ``digon''.  By the Gauss-Bonnet theorem, a
geodesic digon cannot have non-positive curvature.)  Thus the $\kappa =
0$, $n=2$ case of Theorem~\ref{th:positive} implies Weil's theorem; but,
as explained in the previous paragraph, it is more general.

In higher dimensions, there are non-positively curved smooth balls with
closed geodesics.  Hass's construction has closed geodesics, and so does
our construction in Theorem~\ref{th:ball}.

\subsection{Proof of Theorem~\ref{th:ball}}

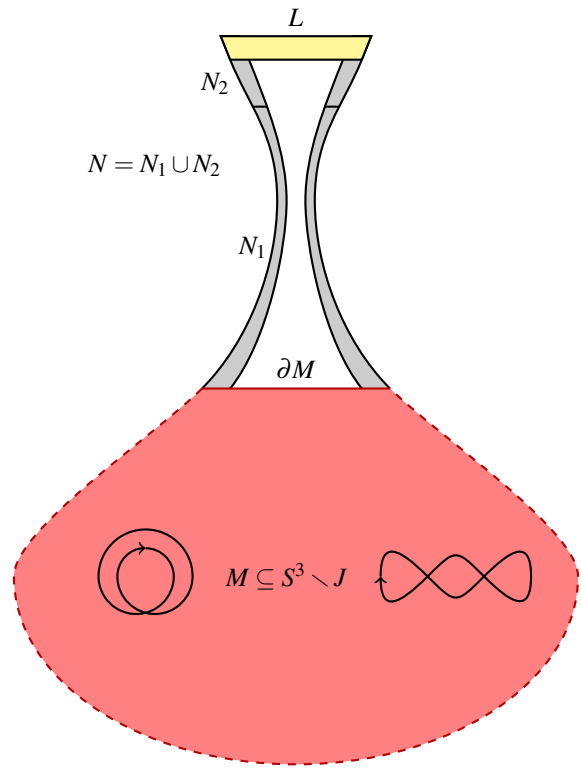
\begin{figure}[htb]
\begin{center} \begin{tikzpicture}[thick,scale=1.25]
\draw[fill=black!20!white] (-1,2)
    .. controls (-.5,2.5) and (-.2,3.3) .. (-.2,4)
    .. controls (-.2,4.7) and (-.5,5) .. (-.7,5.5)
    -- (-.8,5.75) -- (-.6,5.75) -- (-.5,5.5)
    .. controls (-.3,5) and (-.1,4.5) .. (-.1,4) 
    .. controls (-.1,3.5) and (-.3,2.5) .. (-.7,2);
\draw[fill=black!20!white] (1,2)
    .. controls (.5,2.5) and (.2,3.3) .. (.2,4)
    .. controls (.2,4.7) and (.5,5) .. (.7,5.5)
    -- (.8,5.75) -- (.6,5.75) -- (.5,5.5)
    .. controls (.3,5) and (.1,4.5) .. (.1,4) 
    .. controls (.1,3.5) and (.3,2.5) .. (.7,2);
\draw (-.47,5) -- (-.3,5) (.47,5) -- (.3,5);
\draw[fill=lightyellow] (.8,5.75) -- (.7,5.5) -- (-.7,5.5)
    -- (-.8,5.75) -- cycle;
\path[dashed,draw=darkred,fill=red!50!white] (-1,2) .. controls
    (-1.5,1.5) and (-3,.5) .. (-3,0) arc (180:360:3 and 2) .. controls
    (3,.5) and (1.5,1.5) .. (1,2);
\draw[darkred] (1,2) -- (-1,2);
\draw[anchor=east] (-.6,5.25) node {$N_2$};
\draw (-.2,3.5) node[anchor=east] {$N_1$};
\draw (-1.5,4.375) node {$N = N_1 \cup N_2$};
\draw (-0,2) node[anchor=south] {$\del M$};
\draw (-.1,0) node {$M \subseteq S^3 \setminus J$};
\draw[anchor=south] (0,5.75) node {$L$};
\draw[->] (-1.6,.3) arc (90:0:.3) arc (360:180:.4) arc (180:0:.5)
    arc (360:180:.4) arc (180:90:.3);
\draw[->] (.9,0) .. controls (.9,.4) and (1.1,.3) .. (1.4,0) ..
    controls (1.7,-.3) and (1.7,-.3) .. (2,0) .. controls (2.3,.3) and
    (2.5,.4) .. (2.5,0) .. controls (2.5,-.4) and (2.3,-.3) .. (2,0) ..
    controls (1.7,.3) and (1.7,.3) .. (1.4,0) .. controls (1.1,-.3)
    and (.9,-.4) .. (.9,0);
\end{tikzpicture} \end{center}
\caption{A diagram of $\Omega = M \cup N \cup L$, schematically like
    a decanter.  It consists of a truncated hyperbolic knot complement
    $M$, plus a 2-handle $N \cup L$ that consists of a 3-ball lid $L$
    and a thickened annulus neck $N = N_1 \cup N_2$.}
\label{f:decanter} \end{figure}

In our proof of Theorem~\ref{th:ball}, we will construct $\Omega$
to be $(-1-\eps,-1+\eps)$-pinched, so that $\sqrt{1+\eps}\Omega$ is
$(-1,-1+2\eps)$-pinched.  We can then change $\eps$ to match the stated
conclusion of Theorem~\ref{th:ball}.

Our construction is shown schematically in Figure~\ref{f:decanter}.
We make $\Omega$ as a union of three pieces $M$, $N$, and $L$.  $M$ is a
truncated hyperbolic knot complement,  $L$ is a ``lid" which is a horospheric
pseudocylinder, and $N$ is a connecting neck which is a thickened annulus.
Both $M$ and $L$ have constant curvature $K = -1$, while the neck $N$ has
a $(-1-\eps,-1+\eps)$-pinched metric that interpolates between the metrics
on $\del M$ and $\del L$ and meets each one along an annulus.  Although the
two ends of $N$ are both horospheric annuli, they are mismatched in two
ways:  First, $\del M$ and the bottom of $\del L$ are both concave, so
their extrinsic curvature must be interpolated.  Second, the annulus $\del
M \cap N$ is \emph{vertical} (isometric to a cylinder $S^1 \times I$ with
the product metric) while $\del L \cap N$ is \emph{horizontal} (isometric
to an annulus in $\E^2$).  In order to achieve both interpolations, we
further divide $N = N_1 \cup N_2$ into two thickened annuli $N_1$ and $N_2$.

To construct $N_1$ and $N_2$, which will be the most technical
part of the proof, we review some facts about warped products
\cite{BO:negative,AB:warped}.  Recall that if $B$ and $F$ are two Riemannian
manifolds and $h:B \to \R_+$ is a smooth function, we can define a Riemannian
metric on $B \times F$ by the formula
\[ \dd s_{B \times F}^2(p,q) = \dd s^2_B(p) + f(p)^2 \dd s^2_F(q) \]
for $(p,q) \in B \times F$.  The manifold $B \times F$ with this metric
is denoted $B \times_f F$ and is called a \emph{warped product}, while
the function $f$ is a \emph{warping function}.  In this paper, we will
only need warped products of the form $I \times_f M$, where the base $I$
is an interval.

\begin{lemma} Let $I$ be an interval, let $f:I \to \R_+$
be a warping function, let $F$ be a Riemannian manifold, 
and let $W = I \times_f F$.  Then:
\begin{enumerate}
\item If $F$ is locally Euclidean and $f(t) = e^t$, then $W$
has constant curvature $K = -1$.
\item If $F$ has constant curvature $K = -1$ and $f(t) = \cosh(t)$,
 then $W$ also has constant curvature $-1$.
\item If $F$ is 1-dimensional, then the intrinsic curvature of $W$
is given by
\[ K_W(t,x) = -\frac{f''(t)}{f(t)} \]
for $t \in I$ and $x \in F$.
\item If $F$ is $(-1-\eps,-1+\eps)$-pinched and $f(t) = \cosh(t)$, then $W$
is also $(-1-\eps,-1+\eps)$-pinched.
\end{enumerate}
\label{l:warp} \end{lemma}

In the proof of Lemma~\ref{l:warp}, and later in the proof of
Theorem~\ref{th:ball}, we will make use of the standard upper
half-space model for hyperbolic space:
\begin{eq}{e:hmetric}
\dd s_{\H^{n+1}}^2 = \frac{\dd s_{\E^n}^2 + \dd z^2}{z^2}
    = \frac{\dd x_1^2 + \ldots + \dd x_n^2 + \dd z^2}{z^2}.
\end{eq}
Recall that in this model, $x_k \in \R$ for every $k$ and $z > 0$.

\begin{proof} Cases 1, 2, and 3 all follow quickly from the conditions
\eqref{e:cba} and \eqref{e:cbb} below.  However, as these are well-known
facts in differential geometry, we also give separate calculations.
Case 1 is confirmed by a standard metric model of $\H^{n+1}$,
\[ \dd s_{\H^{n+1}}^2 = e^{2t}\dd s_{\E^n}^2 + \dd t^2, \]
which is obtained from \eqref{e:hmetric} by the change of variables $z
= e^{-t}$.  Case 2 is confirmed by a standard metric model of $\H^{n+1}$,
\[ \dd s^2 = \cosh(t)^2\dd s_{\H^n}^2 + \dd t^2, \]
which may also be obtained from \eqref{e:hmetric} by the change of variables
\[ (x_n,z) = (y\tanh(t),y\sech(t)), \]
to obtain the metric
\[ \dd s_{\H^{n+1}}^2 = \cosh(t)^2\frac{\dd x_1^2 + \ldots + \dd x_{n-1}^2 +
    \dd y^2}{y^2} + \dd t^2. \]

Case 3 follows from the Jacobi field equation \eqref{e:jacobi}, considering
that the $I$ fibers in a warped product $I \times_f F$ are geodesic curves.

Case 4 is a special case of a result of Alexander-Bishop \cite[Prop
2.2.]{AB:warped}.  Reducing to the case of a one-dimensional base $B =
I$, this proposition says that:
\begin{enumerate}
\item $W$ has curvature bounded above by $K_0$ if
\begin{eq}{e:cba}
f'' \ge -K_0f \qquad \mbox{and} \qquad K_F \le K_0f^2 + f'^2.
\end{eq}
\item $W$ has curvature bounded below by $K_0$ if
\begin{eq}{e:cbb}
f'' \le -K_0f \qquad \mbox{and} \qquad K_F \ge K_0f^2 + f'^2.
\end{eq}
\end{enumerate}
(The proposition states ``if and only if" and requires that $B$ be
complete, but their proof makes clear that the ``if" direction does not
require completeness.)  If we take $K_0 = -1+\eps$ for the upper
bound and $K_0 = -1-\eps$ for the lower bound, and if $f(t) = \cosh(t)$,
we obtain the requirements
\begin{gather*}
\cosh(t)(1-\eps) \le \cosh(t) \le \cosh(t)(1+\eps) \\
-1-\eps\cosh(t)^2 \le K_F \le -1+\eps\cosh(t)^2.
\end{gather*}
All of these inequalities hold immediately.
\end{proof}

\begin{lemma} Let $\eps, \delta, a, b > 0$.  Then there
exists $c > 2$ and a smooth $f:[0,c] \to \R_+$ such that the manifold
\[ N_1 = [-\delta,\delta] \times \R/2\pi\Z \times [0,c] \]
with coordinates $(\rho,\theta,h)$ and metric
\[ \dd s^2 = \dd \rho^2 + \cosh(\rho)^2(f(h)^2\sdd\theta^2 + \dd h^2) \]
is $(-1-\eps,-1+\eps)$-pinched, and such that
\[ f(h) = \begin{cases} a e^{-h} & h \le 1\\
    b e^{h-c} & h \ge c-1 \end{cases} \]
and
\[ |\del N_1| = O_\eps\big(\cosh(\delta)^2(a+b)\big). \]
\eatline \label{l:n1} \end{lemma}

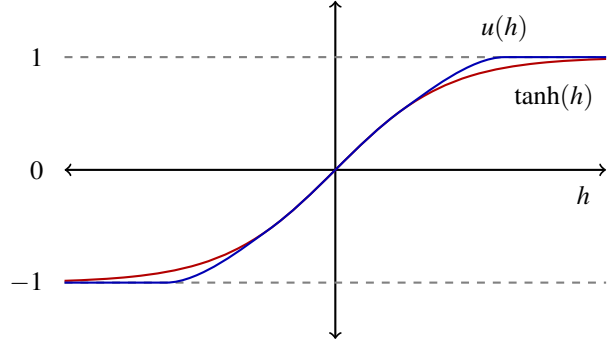
\begin{figure}[htb]
\begin{tikzpicture}[thick,scale=1.5]
\draw[<->] (0,-1.5) -- (0,1.5);
\draw[<->] (-2.4,0) -- (2.4,0);
\draw[dashed,gray] (-2.4,-1) -- (2.4,-1);
\draw[dashed,gray] (-2.4,1) -- (2.4,1);
\draw[darkred] (0,0) .. controls (.2,.2) and (.4,.395) .. (.6,.537)
    .. controls (.8,.679) and (1,.773) .. (1.2,.834)
    .. controls (1.4,.895) and (1.6,.926) .. (1.8,.947)
    .. controls (2,.968) and (2.2,.977) .. (2.4,.984);
\draw[darkred] (0,0) .. controls (-.2,-.2) and (-.4,-.395) .. (-.6,-.537)
    .. controls (-.8,-.679) and (-1,-.773) .. (-1.2,-.834) 
    .. controls (-1.4,-.895) and (-1.6,-.926) .. (-1.8,-.947)
    .. controls (-2,-.968) and (-2.2,-.977) .. (-2.4,-.984);
\draw[darkblue] (0,0) .. controls (.2,.2) and (.4,.395) .. (.6,.537)
    .. controls (.8,.679) and (1.2,1) .. (1.5,1) -- (2.4,1);
\draw[darkblue] (0,0) .. controls (-.2,-.2) and (-.4,-.395) .. (-.6,-.537)
    .. controls (-.8,-.679) and (-1.2,-1) .. (-1.5,-1) -- (-2.4,-1);
\draw[anchor=north] (2.2,-.05) node {$h$};
\draw[anchor=north west] (1.5,.85) node {$\tanh(h)$};
\draw[anchor=south] (1.5,1.05) node {$u(h)$};
\draw[anchor=east] (-2.5,-1) node {$-1$} (-2.5,0) node {$0$} (-2.5,1) node {$1$};
\end{tikzpicture}
\caption{The function $u(h)$, an approximate Riccati solution that
    transitions from $-1$ to $1$ and which partly agrees with the exact
    solution $\tanh(h)$}
\label{f:upinch} \end{figure}

\begin{proof} Following Lemma~\ref{l:warp}, $N_1$ is
$(-1-\eps,-1+\eps)$-pinched if and only if
\[ (1-\eps)f(h) \le f''(h) \le (1+\eps)f(h). \]
This relation holds if and only if the logarithmic derivative
$u(h) = f'(h)/f(h)$ approximately satisfies a Riccati equation:
\begin{eq}{e:upinch} 1-\eps \le u'(h) + u(h)^2 \le 1+\eps. \end{eq}
We first construct $u(h)$, for convenience for all $h \in \R$.  Actually,
we will shift $h$ by a constant, which has no effect on \eqref{e:upinch}.
Let $u(h)$ be a smooth function such that:
\begin{enumerate}
\item $u(h) = \tanh(h)$ when $|\tanh(h)| \le 1-\frac{\eps}2$.
\item $0 \le u'(h) \le \eps$ when $|\tanh(h)| > 1-\frac{\eps}2$. 
\item $u(h)$ increases with $h$ until it reaches 1 at $h=c_0$ and then
stays constant.
\item $u(-h) = -u(h)$ for all $h$.
\end{enumerate}
(See Figure~\ref{f:upinch}.)  When $|\tanh(h)| \le 1-\frac{\eps}2$, we have
\[u'(h) + u(h)^2 = 1. \]
For other values of $h$, we have
\[ 1-\eps < u(h)^2 \le 1, \qquad 0 \le u'(h) \le \eps, \]
so that \eqref{e:upinch} is satisfied in all cases.

Choose $c_1 \le -c_0-1$ and $c_2 \ge c_0+1$ such that also 
\[ \int_{c_1}^{c_2} u(t) \sdd t = \log(b) -\log(a), \]
and let $c = c_2 - c_1$.  Then
\[ f(h) = a\exp\bigg(\int_{c_1}^{h+c_1} \back \back u(t) \sdd t\bigg)\]
has all of the required properties on the interval $[0,c]$.

We can estimate $|\del N_1|$ by first considering the area of the level
surface $\rho = 0$, which is $\int_0^c 2\pi f(h)\sdd h$.  Splitting the
integral of $f$ into
\[ \int_0^c f(h)\sdd h = \int_0^{-c_0-c_1} \back f(h) \sdd h
    + \int_{-c_0-c_1}^{c_0-c_1} \back f(h)\sdd h
    + \int_{c_0-c_1}^c \back f(h) \sdd h, \]
we observe that
\begin{align*}
\int_0^{-c_0-c_1} \back \back f(h) \sdd h
    &= \int_0^{-c_0-c_1} \back\back ae^{-h} \sdd h < a \\
\int_{c_0-c_1}^c \back \back f(h) \sdd h
    &= \int_{c_0-c_1}^c \back\back be^{h-c} \sdd h < b \\
\int_{-c_0-c_1}^{c_0-c_1} \back\back f(h)\sdd h
    &< 2c_0\max(a,b) = O_\eps(a+b).
\end{align*}
The estimate for $|\del N_1|$ follows quickly.
\end{proof}

\begin{lemma} Let $\eps > 0$ and let $f:[0,3] \to
[0,1]$ be a smooth function such that $f(h) = 0$ for $h \le 1$
and $f(h) = 1$ for $h \ge 2$.  Then there exists $b>0$ such that
the manifold
\[ N_2 = [-1,1] \times \R/2\pi\Z \times [0,3] \]
with coordinates $(\rho,\theta,h)$ and with the metric
\begin{eq}{e:n2interp}
\dd s^2 = e^{2h}\big(\dd \rho^2 + (b+f(h)\rho)^2\sdd\theta^2\big) + \dd h^2,
\end{eq}
is $(-1-\eps,-1+\eps)$-pinched.
\label{l:n2} \end{lemma}

Note that if $f(h)$ is locally constant near $h=h_0$, then
Lemma~\ref{l:warp} tells us that the metric \eqref{e:n2interp} has constant
negative curvature $K=-1$ at $h_0$.  In particular, $K=-1$ for every $h \in
[0,1]\cup[2,3]$.

\begin{proof} For clarity, we work in the universal cover $\tN_2$,
so that $\theta \in \R$.  Without yet choosing $b$, we apply the change
of variables $\theta = \alpha/b$ to obtain the metric
\begin{eq}{e:cmetric} \dd s^2 = e^{2h}\bigg(\dd \rho^2 +
    \big(1+\frac{f(h)\rho}b\big)^2\sdd\alpha^2\bigg) + \dd h^2. \end{eq}
Since $h$ and $\rho$ are bounded independently of the choice of $b$, and
since $f(h)$ is a fixed, smooth function, the metric \eqref{e:cmetric}
converges uniformly in the $C^\infty$ topology to the metric
\begin{eq}{e:vmetric} \dd s^2 = e^{2h}(\dd \rho^2 + \sdd\alpha^2) + \dd h^2
\end{eq}
as $b \to \infty$.  Recall that the curvature tensor $R_{ijkl}$ of a
manifold $M$ with a metric $g_{ij}$ has a polynomial formula in terms of the
derivatives $g_{ij,k}$ and $g_{ij,kl}$, and the matrix inverse $g^{ij}$.
In our case, the limiting metric \eqref{e:vmetric} has constant curvature
$K=-1$ by Lemma~\ref{l:warp}.  Both $g_{ij}$ and $g^{ij}$ are uniformly
bounded, since they are independent of the non-compact coordinate $\alpha$.
It follows that the sectional curvature of $\tN_2$ converges uniformly to
$K=-1$ as $b \to \infty$. This is equivalent to the conclusion that $\tN_2$
or $N_2$ is $(-1-\eps,-1+\eps)$-pinched when $b$ is sufficiently large.
\end{proof}

\begin{proof}[Proof of Theorem~\ref{th:ball}] 
Let $J \subseteq S^3$ be a hyperbolic knot and give its complement
$S^2\setminus J$ its complete hyperbolic metric with curvature $K = -1$.
We can choose $J$ so that $|S^3 \setminus J|$ is arbitrarily high by a
theorem of Adams \cite{Adams:handbook}.  A tubular neighborhood of $J$
is metrized as a parabolic cusp, and this cusp can be truncated to obtain
a manifold $M$ with a horospheric torus boundary $\del M$.  Let $\gamma
\subset \del M$ be a geodesic meridian circle.  Note that we can truncate
$|S^3 \setminus J|$ as far out along the cusp as we like.  We choose the
truncation so that
\[ |M| \ge |S^3 \setminus J|-1, \qquad |\del M| \le 1,
    \qquad |\gamma| \le 1. \]

If we attach a 2-handle $D^2 \times I$ to $M$ along $\gamma$, then the result
\[ \Omega = M \cup (D^2 \times I) \]
is diffeomorphic to a ball $B^3$.  (Recall that an $n$-dimensional $k$-handle
is a $B^k \times B^{n-k}$ which is to be attached along $(\del B^k) \times
B^{n-k}$ to some other manifold.)  We want to give the handle $D^2 \times I$
a $(-1-\eps,-1+\eps)$-pinched metric that extends smoothly to the metric on $M$.
We also want to bound $|\del (D^2 \times  I)|$ by a constant, independent of
the choice of $M$ (but depending on $\eps$).

We construct the 2-handle $D^2 \times I$ as the union of a thickened annulus
$N$, the ``neck"; and a 3-ball $L$, the ``lid".  The connecting neck $N$
is divided into two stages, $N_1$ and $N_2$.  Each $N_k$ is of the form
$I_\ishort \times S^1 \times I_\ilong$, where $S^1 \times I_\ilong$ is
thus a long annulus.   More precisely (as in Figure~\ref{f:decanter}),
we will glue $M$ to $N_1$, $N_1$ to $N_2$, and $N_2$ to $L$.  It will be
convenient for each pair to overlap with positive volume.

We construct $N_1$ and $N_2$ using Lemmas~\ref{l:n1} and \ref{l:n2},
which both use the coordinates $(\rho,\theta,h)$.  To distinguish them,
we change the variables names to $\rho_1$ and $h_1$ in Lemma~\ref{l:n1}
and to $\rho_2$ and $h_2$ in Lemma~\ref{l:n2}; the coordinate $\theta$
will be the same.

We choose the constant $a$ in Lemma~\ref{l:n1} so that $|\gamma| =
2\pi a$, and we choose the constant $b$ to match its value provided in
Lemma~\ref{l:n2}.   We parametrize $\gamma$ by $\theta \in \R/2\pi\Z$ so that
$a\theta$ represents the length along $\gamma$ from some starting point.
Since $\gamma$ is horocyclic, it has a neighborhood in $M$ with coordinates
$(\rho_1,\theta,h_1)$ with the metric
\[ \dd s^2 = \cosh(\rho_1)^2(a^2e^{-2h_1}\dd \theta^2
    + \dd h_1^2) + \dd\rho_1^2, \]
and where $\gamma$ itself in these coordinates is $\gamma(\theta) =
(0,\theta,0)$ and $(0,\theta,h_1)$ is at a distance of $h_1$ from $\del M$.
Choose $\delta$ so that the region $\rho_1 \in [-\delta,\delta]$ and
$h_1 \ge -\delta$ is an embedded neighborhood of $\gamma$.  We also want
$\sinh(\delta) \le 1$, for a reason that we will discuss later.
We use these coordinates to glue $M$ to $N_1$, with its metric provided
by Lemma~\ref{l:n1}.

\begin{figure}
\begin{tikzpicture}[thick,scale=1.5]
\useasboundingbox (-3,-.5) -- (3,3.5);
\draw[<->] (-2.5,0) -- (2.5,0);
\draw[->] (0,0) -- (0,3.5);
\fill[blue,opacity=.3] (1,.368) -- (1,1) -- (-1,1) -- (-1,.368);
\draw (1,.368) -- (1,1) -- (-1,1) -- (-1,.368);
\draw[dashed] (1,.368) -- (-1,.368);
\fill[blue,opacity=.3] (1.631,2.175) -- (.6,.8) arc
    (53.1:126.9:1) -- (-1.631,2.175) arc (126.9:53.1:2.718); 
\draw (1.631,2.175) -- (.6,.8) arc (53.1:126.9:1) -- (-1.631,2.175);
\draw[dashed] (-1.631,2.175) arc (126.9:53.1:2.718);
\draw (-.7,1.8) node {$N_1$};
\draw (-.75,.6) node {$N_2$};
\draw[anchor=south west] (0,.368) node {$1/e$}; \fill (0,.368) circle (.05); 
\draw[anchor=south west] (0,1) node {$1$}; \fill (0,1) circle (.05); 
\draw[anchor=south west] (0,2.718) node {$e$}; \fill (0,2.718) circle (.05); 
\draw[anchor=north] (-1,-.05) node {$-1$}; \fill (-1,0) circle (.05);
\draw[anchor=north] (1,-.05) node {$1$}; \fill (1,0) circle (.05);
\draw[anchor=south] (2.2,.05) node {$x_1$};
\draw[anchor=east] (-.05,3.2) node {$z_1$};
\draw[->] (2,1.6) -- (1.3,1.6);
\draw[anchor=west] (2,1.6) node
    {$\displaystyle \frac{z_1}{x_1} = \frac1{\sinh(\delta)}$};
\end{tikzpicture}
\caption{Gluing $N_1$ to $N_2$ in the half-space model in coordinates $x_1$
and $z_2$.  If $\sinh(\delta) \le 1$, then $N_1$ inserts into the end of $N_2$
as shown.}
\label{f:n1n2} \end{figure}
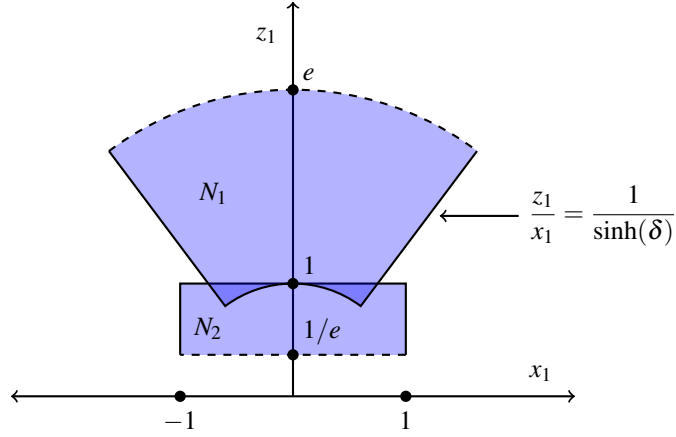

To glue $N_1$ to $N_2$, we introduce the coordinates $(x_1,y_1,z_1)$ with
$x_1,z_1 \in \R$ and $\theta \in \R/2\pi\Z$, and with the hyperbolic metric
\[ \dd s^2 = \frac{\dd x_1^2 + \dd y_1^2 + \dd z_1^2}{z_1^2}. \]
As in the proof of Lemma~\ref{l:warp}, we can identify these coordinates
with both $(\rho_1,\theta,h_1)$ and $(\rho_2,\theta,h_2)$ using
the equations
\begin{align*}
(x_1,y_1,z_1) &= (e^{c-h_1}\tanh(\rho_1),b\theta,e^{c-h_1}\sech(\rho_1)) \\
(x_1,y_1,z_1) &= (\rho_2,b\theta,e^{-h_2}).
\end{align*}
Both changes of variables preserve the defined metrics.  We can also solve
for $(\rho_2,h_2)$ in terms of $(\rho_1,h_1)$ to obtain
\[ (\rho_2,h_2) = (e^{c-h_1}\tanh(\rho_1),h_1-c-\log(\sech(\rho_1))). \]
If $\sinh(\delta) \le 1$, the constant curvature end of
$N_1$ stays within the constant curvature end of $N_2$, and is inserted
between its corners, as in Figure~\ref{f:n1n2}.

Finally we define the lid $L$ to be the region
\[ L = \{(x_2,y_2,h_2) \st x_2^2+y_2^2 \le b^2, 2 \le h \le 3\} \]
in $\R^3$ with metric
\[ \dd s^2 = e^{2h_2}(\dd x_2^2 + \dd y_2^2) + \dd h_2^2. \]
We glue $N_2$ to $L$ by changing to polar coordinates,
\[ (x_2,y_2,h_2) = ((\rho_2+b)\cos(\theta),(\rho_2+b)\sin(\theta),h_2), \]
which also converts between the metrics on $L$ and $N_2$.

By construction, $N_1 \cup N_2 \cup L$ is a 2-handle that attaches to $M$
along $\gamma$.  The result is a  $(-1-\eps,-1+\eps)$-pinched 3-ball
\[ \Omega = M \cup N_1 \cup N_2 \cup L \]
with piecewise smooth boundary; we can make the boundary smooth by shaving
it slightly.   Since $M$ itself has arbitrarily large volume, $\Omega$
does too.  It remains to bound the surface area
\[ |\del \Omega| \le |\del M| + |\del N_1| + |\del N_2| + |\del L|. \]
The submanifolds $N_1$ and $L$ only depends on $\eps$.   Meanwhile we have
already specified that $|\del M| \le 1$.  Finally $|\del N_1|$ is estimated
in Lemma~\ref{l:n1} in terms of the constants $a$ and $b$.  The constant $a$
is bounded because $|\gamma| \le 1$, while the constant $b$ only depends on
$\eps$.  Thus $|\del \Omega|$ is bounded by a constant, depending on $\eps$.
\end{proof}

\section{Geodesic integrals}
\label{s:geodesics}

In this section we will study Santal\'o's integral formula
\cite[Sec. 19.4]{Santalo:book} in the formalism of geodesic flow and
symplectic geometry.  See McDuff and Salamon \cite[Sec. 5.4]{MS:intro} for
properties of symplectic quotients.  The formulas we derive are those of
Croke \cite{Croke:sharp}; see also Teufel \cite{Teufel:riemannian}.

\subsection{Symplectic geometry}
\label{s:symplectic}

Let $W$ be an open symplectic $2n$-manifold with a symplectic
form $\omega_W$.  Then $W$ also has a canonical volume form $\mu_W
= \omega_W^{\wedge n}$ which is called the Liouville measure on $W$.
Let $h:W \to \R$ be a Hamiltonian, by definition any smooth function on
$W$, suppose that $0$ is a regular value of $h$, and let $H = h^{-1}(0)$
then be the corresponding smooth level surface.  Then $\omega_W$ converts
the 1-form $\dd h$ to a vector field $\xi$ which is tangent to $H$.
Suppose that every orbit $\gamma$ of $\xi$ only exists for a finite time
interval.  Let $G$ be the set of orbits of $\xi$ on $H$; it is a type of
\emph{symplectic quotient} of $W$.  $G$ is a smooth open manifold except
that it might not be Hausdorff.  

The manifold $G$ is also symplectic with a canonical 2-form $\omega_G$ and
its own Liouville measure $\mu_G$.  $H$ cannot be symplectic since it is
odd{\hyp}dimensional, but it does have a Liouville measure $\mu_H$.  (In fact
$G$ and $\omega_G$ only depend on $H$, and not otherwise on $h$, while
$\mu_H$ depends on the specific choice of $h$.)  Let $(a(\gamma),b(\gamma))$
be the time interval of existence of $\gamma \in G$; here only the difference
\[ \ell(\gamma) = b(\gamma)-a(\gamma) \]
is well-defined by the geometry.  In this general setting, if $f:H \to \R$
is a suitably integrable function, then
\begin{eq}{e:sympint}
\int_H f(x) \sdd\mu_H(x) = \int_{\gamma \in G}
    \int_{a(\gamma)}^{b(\gamma)} f(\gamma(t)) \sdd t \sdd\mu_G(\gamma).
\end{eq}
Or, if $\sigma$ is a measure on $H$, we can consider the push-forward
$(\pi_G)_*(\sigma)$ of $\sigma$ under the projection $\pi_G:H \to G$.
Taking the special case that $f$ is constant on orbits of $\xi$, the
relation \eqref{e:sympint} says that
\[ (\pi_G)_*(\mu_H) = \ell\mu_G. \]

\subsection{The space of geodesics and \'etendue}
\label{s:etendue}

\begin{figure}[htb]
\begin{center} \begin{tikzpicture}[semithick,scale=1.25]
\fill[pattern=horizontal lines,pattern color=darkblue]
    (0,0) .. controls (-.5,0) and (-.5,1.5) .. (-1.25,1.5) arc (90:180:.75)
    arc (180:360:2) arc (0:90:.75) .. controls (.5,1.5) and (.5,0) .. (0,0);
\draw[dashed]
    (0,0) .. controls (-.5,0) and (-.5,1.5) .. (-1.25,1.5) arc (90:180:.75)
    arc (180:360:2) arc (0:90:.75) .. controls (.5,1.5) and (.5,0) .. (0,0);
\draw[->] (2.2,0) -- (3,0);
\draw (4,-1.25) -- (4,0);
\draw (3.8,0) arc (0:20:4.386) node (a) {};
\draw (4.2,0) arc (180:160:4.386) node (b) {};
\fill (3.8,0) circle (.07);
\fill (4.2,0) circle (.07);
\draw[fill=white] (4,-1.25) circle (.07);
\draw[fill=white] (a) circle (.07);
\draw[fill=white] (b) circle (.07);
\draw[fill=white] (4,0) circle (.07);
\end{tikzpicture} \end{center}
\caption{A manifold $M$ in which the space of geodesics is not Hausdorff.
    The horizontal chords make a ``zipper" 1-manifold.}
\label{f:nonhausdorff} \end{figure}
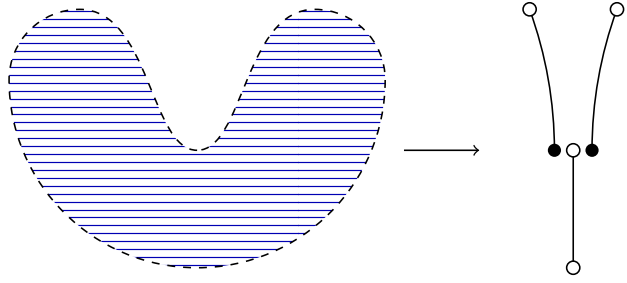

If $M$ is a smooth $n$-manifold, then $W = T^*M$ is canonically a symplectic
manifold.  If $M$ has a Riemannian metric $g$, then $g$ gives us a canonical
identification $TM \cong T^*M$.  It also gives us a Hamiltonian $h:TM \to
\R$ defined as
\[ h(v) = (g(v,v)-1)/2. \]
The level surface $h^{-1}(0)$ is evidently the unit tangent bundle $UM$.
It is less evident, but still routine, that the Hamiltonian flow $\xi$
of $h$ is the geodesic flow on $UM$.  Suppose further that $M$ only has
bounded-time geodesics.  Then the corresponding symplectic quotient $G$
is the space of oriented geodesics on $M$.  The structure on $G$ that
particularly interests us is its Liouville measure $\mu_G$.  The Liouville
measure on $H = UM$ is also important, and happens to equal the Riemannian
measure $\nu_{UM}$.   Even in this special case, $G$ might not be Hausdorff
if the geodesics of $M$ merge or split, as in Figure~\ref{f:nonhausdorff}.

The Liouville measure $\mu_G$ is important in geometric optics
\cite{Smith:optical}, where among other names it is called
\emph{\'etendue}\footnote{In English, not just in French.}.  Lagrange
established that \'etendue is conserved.  Mathematically, this says
exactly that the $(2n-2)$-form $\mu_G$, which is definable on $H$,
descends to $G$.  More explicitly, suppose (in the full generality of
Section~\ref{s:symplectic}) that $K_1,K_2 \subseteq H$ are two transverse
open disks that are identified by the holonomy map
\[ \phi:K_1 \stackrel{\cong}{\longto} K_2 \]
induced by the set of orbits.  Then the Liouville measures $K_1$ and $K_2$
match, \ie, $\phi_*(\mu_G) = \mu_G$.  (As in the proof of Liouville's
theorem, $\phi$ is even a symplectomorphism.)

Now suppose that $\Omega$ is a compact Riemannian manifold with boundary
and with unique geodesics, and let $M$ be the interior of $\Omega$.
Then $G$, the space of oriented geodesics of $\Omega$ or $M$, is canonically
identified in two ways to $U^+ \dOmega$.  We can let $\gamma = \gamma_u$ be
the geodesic generated by $u$, or we can let $\gamma = \overline{\gamma_u}$
be the geodesic with the same image but inverse orientation.
These are both examples of identifying
part of $G$, in this case all of $G$, with transverse submanifolds as in
the previous paragraph.  Let
\[ \sigma_+:U^+\dOmega \stackrel{\cong}{\longto} G, \qquad
    \sigma_-:U^+\dOmega \stackrel{\cong}{\longto} G \]
be the two corresponding identifications.

The maps $\sigma_\pm$ are smooth bijections; when $\dOmega$ is convex,
they are diffeomorphisms.  In general, the inverses $\sigma^{-1}_\pm$
are smooth away from the non-Hausdorff points of $G$. These correspond to
geodesics tangent to $\dOmega$, and they are a set of measure 0 in $G$.
Thus, the composition
\[ \phi = \sigma_-^{-1} \circ \sigma_+:U^+\dOmega \stackrel{\cong}{\longto}
    U^+\dOmega \]
is an involution of $U^+\dOmega$ that preserves the measure $\mu_G$.  (It is
even almost everywhere a local symplectomorphism with respect to $\omega_G$.)
We call the map $\phi$ the \emph{optical transport} of $\Omega$.

We define several types of coordinates on $G$, $U\Omega$ and $U^+\dOmega$.
Let $u = (x,v)$ be the position and vector components of a tangent vector $u
\in U\Omega$, and let $u = (p,v)$ be the same for $u \in U^+\dOmega$.  On $G$
itself, we already have the length function $\ell(\gamma)$.  In addition,
if $\gamma = \gamma_u$ for
\[ u = (p,v) \in U^+\dOmega, \]
let $\alpha(\gamma)$ be the angle between $v$ and the inward normal vector
$w(p)$.  If $\gamma = \overline{\gamma_u}$, then let $\beta(\gamma)$
be that angle instead.

The map $\sigma_+$ relates the Liouville measure $\mu_G$ with
Riemannian measure $\nu_{U^+\dOmega}$.  More loosely, the projection $\pi_G$
from Section~\ref{s:symplectic} relates $\mu_G$ with $\nu_{U\Omega}$.
Then by slight abuse of notation,
\begin{eq}{e:etendue}
\dd\mu_G = \cos(\alpha) \sdd\nu_{U^+\dOmega} = \frac{\dd\nu_{U\Omega}}\ell .
\end{eq}
In words, $\mu_G$ is close to $\nu_{U^+\dOmega}$ but not the same:
If a beam of light is incident to a surface at an angle of $\alpha$,
then its illumination has a factor of $\cos(\alpha)$.  The measure
$(\pi_G)_*(\nu_{U\Omega})$ is also close but not the same, because the
\'etendue of a family of geodesics does not grow with the length of the
geodesics.

Another important comparison of measures relates geodesics to pairs
of points.

\begin{lemma}
Suppose that $p,q \in \Omega$ lie on a geodesic $\gamma$ and that $p
\ne q$.  Let $p = \gamma(a)$ and $q = \gamma(b)$.  Then
\[ \dd\nu_{\Omega \times \Omega}(p,q) =
    j_\Omega(\gamma,a,b) \sdd\mu_G(\gamma)\sdd a \sdd b. \]
In the cases $(p,q)\in \Omega\times \del\Omega$,
$(p,q)\in\del\Omega\times\Omega$ and $(p,q)\in\del\Omega\times\del\Omega$
we respectively have
\begin{align*}
\dd\nu_{\Omega \times \del\Omega}(p,q) &= \frac{j_\Omega(\gamma,a,b)}
    {\cos\beta(\gamma)} \sdd\mu_G(\gamma)\sdd a \\
\dd\nu_{\del\Omega \times \Omega}(p,q) &= \frac{j_\Omega(\gamma,a,b)}
    {\cos\alpha(\gamma)} \sdd\mu_G(\gamma)\sdd b \\
\dd\nu_{\del\Omega \times \del\Omega}(p,q) &= \frac{j_\Omega(\gamma,a,b)}
    {\cos\alpha(\gamma)\cos \beta(\gamma)} \sdd\mu_G(\gamma).
\end{align*}
\eatline \label{l:etendcand} \end{lemma}

\begin{proof}
We start with the case $p,q\in\Omega$.  On the one hand, a localized
version of formula \eqref{e:etendue} is
\[ \dd\mu_G(\gamma)\sdd a = \dd\nu_{U\Omega}(u)
    = \dd\nu_{U_p\Omega}(v) \sdd\nu_\Omega(p) \]
when $u=(p,v)$ and $\gamma$ is the geodesic such that $\gamma(a)=p$
and $\gamma'(a)=v$.  On the other hand, by the definition of the candle
function, we have for all fixed $p$:
\[ \dd\nu_\Omega(q) = j_\Omega(\gamma,a,b) \sdd\nu_{U_p\Omega}(v) \sdd b \]
when $\gamma(a)=p$, $\gamma(b)=q$ and $v=\gamma'(a)$.  Thus,
we obtain a pair of equalities of measures:
\[ \dd\mu_G(\gamma)\sdd a\sdd b = \dd\nu_\Omega(p)\sdd\nu_{U_p\Omega}(v)\sdd b
    = \frac{\dd\nu_\Omega(p)\sdd\nu_\Omega(q)}{j_\Omega(\gamma,a,b)}. \]
The other cases are handled in the same way.  When $p\in\Omega$ and $q\in
\del\Omega$ we have
\begin{align*}
\dd\mu_G(\gamma)\sdd a &= \dd\nu_{U_p\Omega}(v) \sdd\nu_\Omega(p) \\
\dd\nu_{\del\Omega}(q) &= \frac{j_\Omega(\gamma,a,b)}
    {\cos\beta(\gamma)}\sdd\nu_{U_p\Omega}(v)
\end{align*}
where $b$ is entirely determined by $p$ and $v$, instead of being a variable;
when $p\in\del\Omega$ and $q\in \Omega$ then $a$ can be fixed to $0$ by
choosing a suitable parametrization of geodesics, and we have
\begin{align*}
\dd\mu_G(\gamma)\sdd b &= \cos\alpha(\gamma)
    \sdd\nu_{U_p^+\del\Omega}(v) \sdd\nu_{\del\Omega}(p)\sdd b \\
\dd\nu_{\Omega}(q) &=j_\Omega(\gamma,a,b) \sdd \nu_{U_p^+\del\Omega}(v) \sdd b
\end{align*}
and finally when $(p,q)\in\del\Omega\times \del\Omega$ we use
\begin{align*}
\dd\mu_G(\gamma) &= \cos\alpha(\gamma)\sdd\nu_{U_p^+\del\Omega}(v)
    \sdd\nu_{\del\Omega}(p) \\
\dd\nu_{\del\Omega}(q) &=\frac{j_\Omega(\gamma,a,b)}{\cos\beta(\gamma)}
    \sdd\nu_{U_p^+\del\Omega}(v). \qedhere
\end{align*}
\end{proof}

Lemma \ref{l:etendcand} has the important corollary that the candle function
is symmetric.  To generalize from an $\Omega$ with unique geodesics to
an arbitrary $M$, we can let $\Omega$ be a neighborhood of the geodesic
$\gamma$, immersed in $M$.

\begin{corollary}[Folklore {\cite[Lem. 5]{Yau:isoperimetric}}]
In any Riemannian manifold $M$,
\[ j_M(\gamma,a,b) = j_M(\gamma,b,a). \]
\eatline \label{c:iseeyou} \end{corollary}

Combining \eqref{e:sympint} with \eqref{e:etendue} yields Santal\'o's
equality.

\begin{theorem}[Santal\'o {\cite[Sec. 19.4]{Santalo:book}}] If $\Omega$
is as above, and if $f:U\Omega \to \R$ is a continuous function, then
\[ \int_{U\Omega} f(u) \dd\nu_{U\Omega}(u)
    = \int_{U^+\dOmega} \int_0^{\ell(\gamma_u)} \back f(\gamma_u(t))
    \cos(\alpha(u)) \sdd t \sdd\nu_{U^+\dOmega}(u). \]
\eatline \label{th:santalo} \end{theorem}

Finally, we will consider another reduction of the space $G$, the projection
\[ \pi_\lab:G \to \R_{\ge 0} \times [0,\pi/2)^2, \quad
    \pi_\lab(\gamma) = (\ell(\gamma),\alpha(\gamma),\beta(\gamma)). \]
Let
\[ \mu_\Omega = (\pi_\lab)_*(\mu_G) \]
be the push-forward of Liouville measure.  Then $\mu_\Omega$ is a
measure-theoretic reduction of the optical transport map $\phi$, and
is close to a transportation measure in the sense of Monge-Kantorovich.
More precisely, equation \eqref{e:etendue} yields a formula for the $\alpha$
and $\beta$ marginals of $\mu_\Omega$, so we can view $\mu_\Omega$, or
rather its projection to $[0,\pi/2)^2$, as a transportation measure from
one marginal to the other.  The projection onto the $\alpha$ coordinate is
\[ \alpha_*(\mu_\Omega) \defeq \int_{\ell,\beta} \dd\mu_\Omega =
    |\dOmega| \omega_{n-2}\sin(\alpha)^{n-2}\cos(\alpha)\sdd\alpha, \]
where as in \eqref{e:illumal} we use the volume of a latitude sphere on
$U^+_p \dOmega$.  Using the abbreviation
\[ z(\theta) = \frac{\omega_{n-2}\sin(\theta)^{n-1}}{n-1}, \]
we can give a simplified formula for both marginals:
\begin{eq}{e:margin}
\alpha_*(\mu_\Omega) = |\dOmega| \sdd z(\alpha), \qquad
\beta_*(\mu_\Omega) = |\dOmega| \sdd z(\beta).
\end{eq}
A final important property of $\mu_\Omega$ that follows from its construction
is that it is symmetric in $\alpha$ and $\beta$.

\subsection{The core inequalities}
\label{s:core}

In this section, we establish three geometric comparisons that convert
our curvature hypotheses to linear inequalities that can then be used for
linear programming. (Section~\ref{s:etendue} does not use either unique
geodesics or a curvature hypothesis.  Thus, the results there are not
strong enough to establish an isoperimetric inequality.)

\begin{lemma} If $\Omega$ is $\Candle(\kappa)$ and has unique geodesics,
then:
\begin{align}
\int_{\ell,\alpha,\beta} \frac{s_{n,\kappa}(\ell)}{\cos(\alpha)\cos(\beta)}
    \sdd \mu_\Omega & \le |\dOmega|^2 && \mbox{(Croke)} \label{e:croke} \\
\int_{\ell,\alpha,\beta} \frac{s_{n,\kappa}^{(-1)}(\ell)}{\cos(\alpha)}
    \sdd\mu_\Omega & \le |\dOmega| |\Omega|  && \mbox{(Little Prince)}
    \label{e:prince} \\
\int_{\ell,\alpha,\beta} s_{n,\kappa}^{(-2)}(\ell) \sdd\mu_\Omega
    & \le |\Omega|^2 && \mbox{(Teufel)} \label{e:teufel}
\end{align}
If $\Omega$ is convex and has constant curvature $\kappa$, then 
all three inequalities are equalities.
\label{l:core} \end{lemma}

The first case of Lemma~\ref{l:core}, equation \eqref{e:croke}, is due to
Croke \cite{Croke:sharp}.  Equation \eqref{e:prince} generalizes the integral
over $p \in \dOmega$ of equation \eqref{e:pprince} in Theorem~\ref{th:point}.
Finally equation \eqref{e:teufel} generalizes an isoperimetric inequality
of Teufel \cite{Teufel:plane}.  Nonetheless all three inequalities can be
proven in a similar way.

\begin{proof} We define a \emph{partial} map $\tau: \Omega \times \Omega
\to G$ by letting $\tau(p,q)$ be the unique geodesic $\gamma \in G$ that
passes through $p$ and $q$, if it exists.  We define $\tau(p,q)$ only when
$p \ne q$ and only when $\gamma$ is available.  Also, if $\gamma$ exists,
we parametrize it by length starting at the initial endpoint at $0$.

By construction,
\begin{align*}
||\tau_*(\nu_{\dOmega \times \dOmega})|| &\le |\dOmega|^2 \\ 
||\tau_*(\nu_{\dOmega \times \Omega})|| &\le |\dOmega| |\Omega| \\
||\tau_*(\nu_{\Omega \times \Omega})|| &\le |\Omega|^2.
\end{align*}
Note that each inequality is an equality if and only if $\Omega$
is convex.

Using Lemma \ref{l:etendcand}, we can write integrals for each of the
left sides
\begin{align*}
||\tau_*(\nu_{\dOmega \times \dOmega})|| &= \int_G
    \frac{j_\Omega(\gamma,0,\ell)}{\cos(\alpha)\cos(\beta)}\sdd\mu_G(\gamma) \\
||\tau_*(\nu_{\dOmega \times \Omega})|| &= \int_G \int_0^\ell
    \frac{j_\Omega(\gamma,0,r)}{\cos(\alpha)} \sdd r \sdd\mu_G(\gamma) \\
||\tau_*(\nu_{\Omega \times \Omega})|| &= \int_G \int_0^\ell \int_0^t
    j_\Omega(\gamma,r,t) \sdd r \sdd t \sdd\mu_G(\gamma).
\end{align*}
Because $\Omega$ is $\Candle(\kappa)$,
\begin{align*}
j_\Omega(\gamma,0,\ell) &\ge s_{n,\kappa}(\ell), \\
\int_0^\ell j_\Omega(\gamma,0,t) \sdd t &\ge s_{n,\kappa}^{(-1)}(\ell), \\
\int_0^\ell \int_0^t j_\Omega(\gamma,r,t)\sdd r\sdd t
    &\ge s_{n,\kappa}^{(-2)}(\ell),
\end{align*}
and note that each inequality is an equality if $\Omega$ has constant
curvature $\kappa$.  We thus obtain
\begin{align*}
\int_G \frac{s_{n,\kappa}(\ell)}{\cos(\alpha)\cos(\beta)} \sdd\mu_G(\gamma)
    &\le |\dOmega|^2 \\
\int_G \frac{s_{n,\kappa}^{(-1)}(\ell)}{\cos(\alpha)} \sdd\mu_G(\gamma)
    &\le |\dOmega| |\Omega| \\
\int_G s_{n,\kappa}^{(-2)}(\ell) \sdd\mu_G(\gamma) &\le |\Omega|^2.
\end{align*}
Because these integrands only depend on $\ell$, $\alpha$, and $\beta$,
we can now descend from $\mu_G$ to $\mu_\Omega$.
\end{proof}

\subsection{Extended inequalities}
\label{s:extend}

Lemma~\ref{l:core} will yield a linear programming model that is strong
enough to prove Theorem~\ref{th:positive}, but not Theorem~\ref{th:negative}
nor many of the other cases of Theorem~\ref{th:subsume}.  In this section,
we will establish several variations of Lemma~\ref{l:core} using alternate
hypotheses.

The following lemma is the refinement needed for Theorem~\ref{th:negative}
and Theorem~\ref{th:yau}.

\begin{lemma} Suppose that $\Omega$ is a compact domain in an $\LCD(-1)$
Cartan-Hadamard $n$-manifold $M$ and let
\[ \chord(\Omega) \le L \in (0,\infty]. \]
Then
\begin{align}
\int_{\ell,\alpha,\beta} \Big(\frac{s_{n,-1}^{(-1)}(\ell)}{\cos(\alpha)} -
    \frac{(n-1)s_{n,-1}^{(-2)}(\ell)}{\tanh(L)}\Big) \sdd \mu_\Omega
    &\le |\dOmega| |\Omega|- \frac{(n-1)|\Omega|^2}{\tanh(L)} 
\label{e:extend1} \\
\int_{\ell,\alpha,\beta} \Big(\frac{s_{n,-1}(\ell)}{\cos(\alpha)\cos(\beta)} -
    \frac{(n-1)s_{n,-1}^{(-1)}(\ell)}{\tanh(L)\cos(\alpha)}\Big)
    \sdd \mu_\Omega
    &\le |\dOmega|^2-\frac{(n-1)|\dOmega| |\Omega|}{\tanh(L)}.
\label{e:extend2}
\end{align}
If $\omega$ is convex and has constant curvature $\kappa$, then
the inequalities are equalities.
\label{l:extend} \end{lemma}

\begin{proof}[Proof of \eqref{e:extend1}] We abbreviate
\[ s(\ell) \defeq s_{n,-1}(\ell), \]
and we switch $\alpha$ and $\beta$ in the integral.

Let $G$ be the space of geodesics of $\Omega$ and recall the partial map
$\tau:\Omega \times \Omega \to G$ used in the proof of Lemma~\ref{l:core}
and the measures $\nu_{\dOmega \times \Omega}$ and $\nu_{\Omega \times
\Omega}$.  We consider the signed measure
\[ \sigma_{\Omega \times \Omega} \defeq \nu_{\Omega\times\dOmega} 
    - \frac{n-1}{\tanh(L)} \nu_{\Omega\times\Omega}. \]
To be precise, if $(p,q) \in \Omega \times \dOmega$, then $\gamma =
\tau(p,q)$ is the geodesic that passes through $p$ and ends at $q$.  We claim
two things about the pushforward $\tau_*(\sigma_{\Omega \times \Omega})$:
\begin{enumerate}
\item That the net measure omitted by $\tau$ is non-negative:
\[ ||\tau_*(\sigma_{\Omega\times\Omega})||
    \le |\dOmega||\Omega| - \frac{n-1}{\tanh(L)}|\Omega|^2. \]
\item That the measure that is pushed forward is underestimated
by the comparison candle function:
\[ \int_G \Big(\frac{s^{(-1)}(\ell)}{\cos(\beta)} -
    \frac{(n-1)s^{(-2)}(\ell)}{\tanh(L)}\Big) \sdd \mu_G(\gamma) \le 
||\tau_*(\sigma_{\Omega\times\Omega})||. \]
\end{enumerate}
Just as in the proof of Lemma~\ref{l:core}, equation \eqref{e:extend1}
follows from these two claims.

To prove the second claim, let $\gamma \in G$ be a maximal geodesic
of $\Omega$ with unit speed and domain $[0,\ell]$.  We abbreviate
the candle function along $\gamma$:
\[ j(t) \defeq j(\gamma,t), \qquad j(r,t) \defeq j(\gamma,r,t). \]
Since $M$ and therefore $\Omega$ is $\LCD(-1)$, we have
the inequality
\[ \frac{j'(t)}{j(t)} \ge \frac{s'(t)}{s(t)}. \]
We can rephrase this as saying that
\[ \dbyd{j}{t}(r,t) - \frac{s'(t-r)}{s(t-r)}j(r,t)
 = \dbyd{j}{t}(r,t) - \frac{n-1}{\tanh(t-r)}j(r,t) \]
is minimized (with a value of 0) in the $K = -1$ case.   Now
\[ \tanh(t-r) \le \tanh(L), \]
while $\LCD(-1)$ implies $\Candle(-1)$, \ie,
\[ j(r,t) \ge s(t-r). \]
It follows that 
\begin{eq}{e:syau}
\dbyd{j}{t}(r,t) - \frac{(n-1)j(r,t)}{\tanh(L)} \ge s'(t-r)
    - \frac{(n-1)s(t-r)}{\tanh(L)}.
\end{eq}
We can integrate with respect to $r$ and $t$ to obtain:
\begin{align*}
\int_0^\ell \int_r^\ell \Big[\dbyd{j}{t}(r,t) -
    \frac{(n-1)j(r,t)}{\tanh(L)}\Big] \sdd t \sdd r
&= \int_0^\ell j(r,\ell) \sdd r - \frac{n-1}{\tanh(L)}\int_0^\ell
    \int_r^\ell j(r,t) \sdd r \sdd t \\
&\ge s^{(-1)}(\ell) - \frac{(n-1)s^{(-2)}(\ell)}{\tanh(L)}.
\end{align*}
Then, if the terminating angle of $\gamma$ is $\beta$, we can again use
the $\Candle(-1)$ condition to obtain
\[ \int_0^\ell \frac{j(r,\ell)}{\cos(\beta)} \sdd r - \frac{n-1}{\tanh(L)}
    \int_0^\ell \int_0^t j(r,t) \sdd r \sdd t \ge \frac{s^{(-1)}(\ell)}
    {\cos(\beta)} - \frac{(n-1)s^{(-2)}(\ell)}{\tanh(L)}. \]
Since the left side is the fiber integral of $\tau_*(\sigma_{\Omega \times
\Omega})$, as in the proof of Lemma~\ref{l:core}, this establishes the
second claim.

To establish the first claim, for each $p \in \Omega$, we consider the set
$\Omega \setminus \Vis(\Omega,p)$ consisting of points $q \in \Omega$ that
are \emph{not} visible from $p$.   The union of all of these is exactly
the pairs $(p,q)$ where $\tau$ is not defined.  If $\gamma$ is a geodesic
in $M$ emanating from $p$, we can restrict further to its intersection
\[ \gamma \cap (\Omega \setminus \Vis(\Omega,p)), \]
where we extend the geodesic $\gamma$ from $\Omega$ to $M$.  We claim
that the integral of $\sigma_{\Omega \times \Omega}$ on each of these
intersections, with the appropriate Jacobian factor, is non-negative.

To verify this claim, we suppose that the intersection is non-empty, and
we parametrize $\gamma$ at unit speed so that $\gamma(0) = p$.  Let $I$
be the set of times $t$ such that
\[ \gamma(t) \in \Omega \setminus \Vis(\Omega,p) ,\]
let $\{t_k\}$ be the set of right endpoints of $I$ where $\gamma$ leaves
$\Omega$, and for each $k$, let $\beta_k \in [0,\pi/2]$ be the angle that
$\gamma$ exits $\Omega$ at $\gamma(t_k)$.  Let $\ell$ be the rightmost
point of $I$.  Then the infinitesimal portion of $\sigma_{\Omega \times
\Omega}$ on $\gamma(I)$ is
\[ \sum_k \frac{j(t_k)}{\cos(\beta_k)} - \frac{n-1}{\tanh(L)}\int_I j(t)
    \sdd t \ge j(\ell) - \frac{n-1}{\tanh(L)}\int_0^\ell j(t) \sdd t. \]
(In other words, we geometrically simplify to the worst case: $I = [0,\ell]$
and $\beta = 0$.)  The derivative of the right side is now
\begin{eq}{e:yau}
j'(\ell) - \frac{n-1}{\tanh(L)} j(\ell) \ge 0.
\end{eq}
The inequality holds because it is the same as \eqref{e:syau}, except with
the right side simplified to $0$.   This establishes the first claim and
thus \eqref{e:extend1}.

The equality criterion holds for the same reasons as in Lemma~\ref{l:core}.
\end{proof}

\begin{proof}[Proof of \eqref{e:extend2}] The proof has exactly the same
ideas as the proof of \eqref{e:extend1}, only with some changes to the
formulas.  We keep the same abbreviations.   This time we define
\[ \sigma_{\dOmega \times \Omega} \defeq \nu_{\dOmega\times\dOmega} 
    - \frac{n-1}{\tanh(L)}\nu_{\dOmega\times\Omega}, \]
we consider $\tau_*(\sigma_{\dOmega \times \Omega})$, and we claim:
\begin{enumerate}
\item That the net measure omitted by $\tau$ is non-negative:
\[ ||\tau_*(\sigma_{\dOmega\times\Omega})||
    \le |\dOmega|^2 - \frac{n-1}{\tanh(L)} |\dOmega||\Omega|^2. \]
\item That the integral underestimates the pushforward:
\[ \int_G \Big(\frac{s(\ell)}{\cos(\alpha)\cos(\beta)} -
    \frac{(n-1)s^{(-1)}(\ell)}{\cos(\alpha)\tanh(L)}\Big) \sdd\mu_G(\gamma) \\
    \le  ||\tau_*(\sigma_{\dOmega\times\Omega})||. \]
\end{enumerate}

To prove the second claim, we define $\gamma$ and $j$ as before and we
again obtain \eqref{e:syau}.  In this case, we integrate only with respect
to $t \in [0,\ell]$ to obtain
\[ j(0,\ell) - \frac{n-1}{\tanh(L)} \int_0^\ell j(0,t) \sdd t \ge
    s(\ell) - \frac{(n-1)s^{(-1)}(\ell)}{\tanh(L)}. \]
Now divide through by $\cos(\alpha)$, and we use the $\Candle(-1)$ property
to divide the first term $\cos(\beta)$, to obtain
\[ \frac{j(0,\ell)}{\cos(\alpha)\cos(\beta)}
    - \frac{n-1}{\cos(\alpha)\tanh(L)} \int_0^\ell j(0,t) \sdd t
\ge \frac{s(\ell)}{\cos(\alpha)\cos(\beta)}
    - \frac{(n-1)s^{(-1)}(\ell)}{\cos(\alpha)\tanh(L)}. \]
The left side is the fiber integral of $\tau_*(\sigma_{\dOmega\times\Omega})$,
so this establishes the second claim.

The proof of the first claim is identical to the case of \eqref{e:extend1},
except that $p \in \dOmega$, and we divide through by $\cos(\alpha)$.
\end{proof}

Meanwhile Theorem~\ref{th:croke2} requires the following striking inequality
that depends only on the condition of unique geodesics rather than any
bound on curvature.  We omit the proof as the lemma is equivalent
to Lemma 9 of Croke \cite{Croke:some}.

\begin{lemma}[Croke-Berger-Kazdan] If $\Omega$ is a compact Riemannian
manifold with boundary and with unique geodesics, then
\[ \int_{\ell,\alpha,\beta} s_{n,(\pi/\ell)^2}^{(-2)}(\ell) \sdd\mu_\Omega
    \le |\Omega|^2. \]
\label{l:cbk} \end{lemma}

\subsection{Mirrors and multiple images}
\label{s:multiple}

In this section, we establish the geometric inequalities needed
for Theorem~\ref{th:multiple}.   Let $M$ be a Riemannian manifold
with boundary $\dM$ (although $M$ might not be compact), and consider
geodesics that reflect from $\dM$ with equal angle of incidence and angle
of reflection. We then have an extension of the exponential map at any
point $x\in M$ beyond the time of reflection.  This extended exponential
map has a natural Jacobian; thus $M$ has an extended candle function
$j_M(\gamma,r)$ as treated in Section~\ref{s:candles}.  Then we say that $M$
is $\Candle(\kappa)$ in the sense of reflecting geodesics if this Jacobian
satisfies the $\Candle(\kappa)$ comparison.

Let $\hat{G}$ be the space of these geodesics, for simplicity considering
only those geodesics that are never tangent to $\dM$.  Then the results
of Section~\ref{s:etendue} still apply, with only slight modifications.
In particular $M$ might have a compactification $\Omega$ with $\dM = W
\subseteq \dOmega$.  Then \eqref{e:etendue} applies if we replace $\dOmega$
by $\dOmega \setminus W$; Lemma \ref{l:etendcand} holds; etc.

If $\Omega$ has a mirror $W$ as part of its boundary, then some pairs of
points have at least two connecting, reflecting geodesics.   We can suppose
in general that every two points in $(\Omega,W)$ are connected by at
most $m$ geodesics (which is also interesting even if $W$ is empty),
and we can suppose that $(\Omega,W)$ is $\Candle(\kappa)$ in the sense
of reflecting geodesics.  In this case it is straightforward to generalize
Lemma~\ref{l:core}.  The generalization will yield the linear programming
model for Theorem~\ref{th:multiple}.

\begin{lemma} If $(\Omega,W)$ is $\Candle(\kappa)$ and has at most $m$
reflecting geodesics between any pair of points, then:
\begin{align}
\int_{\alpha,\beta,\ell} \frac{s_{n,\kappa}(\ell)}{\cos(\alpha)\cos(\beta)}
    \sdd\mu_\Omega & \le m|\dOmega|^2  \label{e:mcroke} \\
\int_{\alpha,\beta,\ell} \frac{s_{n,\kappa}^{(-1)}(\ell)}{\cos(\alpha)}
    \sdd\mu_\Omega & \le m|\dOmega| |\Omega| \label{e:mprince} \\
\int_{\alpha,\beta,\ell} s_{n,\kappa}^{(-2)}(\ell) \sdd\mu_\Omega
    & \le m|\Omega|^2. \label{e:mteufel}
\end{align}
\eatline \label{l:multiple} \end{lemma}

\begin{proof} The proof is nearly identical to that of Lemma~\ref{l:core}.
In this case
\[ \tau: \Omega \times \Omega \to G \]
is not a partial map, but rather a multivalued correspondence which is
at most $1$ to $m$ everywhere.  We can define a pushforward measure such as
$\tau_*(\nu_{\Omega \times \Omega})$ by counting multiplicities.

By construction:
\begin{align*}
||\tau_*(\nu_{\dOmega \times \dOmega})|| &\le m|\dOmega|^2 \\ 
||\tau_*(\nu_{\dOmega \times \Omega})|| &\le m|\dOmega| |\Omega| \\
||\tau_*(\nu_{\Omega \times \Omega})|| &\le m|\Omega|^2.
\end{align*}
On the other hand,
\begin{align*}
||\tau_*(\nu_{\dOmega \times \dOmega})|| &= \int_G
    \frac{j(\gamma,\ell)}{\cos(\alpha)\cos(\beta)} \sdd\mu_G(\gamma) \\
||\tau_*(\nu_{\dOmega \times \Omega})|| &= \int_G \int_0^\ell
    \frac{j(\gamma,r)}{\cos(\alpha)} \sdd r \sdd\mu_G(\gamma) \\
||\tau_*(\nu_{\Omega \times \Omega})|| &= \int_G \int_0^\ell \int_0^t
    j(\gamma,r,t) \sdd r \sdd t \sdd\mu_G(\gamma).
\end{align*}
Using the $\Candle(\kappa)$ hypothesis, we obtain the desired inequalities.
\end{proof}

Finally, the following generalization of G\"unther's inequality
\cite{Gunther:volume,BC:book} shows that the $\Candle(\kappa)$ condition
is actually useful for reflecting geodesics.  

\begin{proposition} Let $M$ be a Riemannian manifold with $K \le \kappa$
for some $\kappa \in \R$, and suppose that $\dM$ is concave relative
to the interior.  If $\kappa > 0$, suppose also that $\chord(M) <
\pi/\sqrt{\kappa}$.  Then $M$ is $\LCD(\kappa)$ with respect to geodesics
that reflect from $\dM$.
\label{p:mgunther} \end{proposition}

Proposition~\ref{p:mgunther} generalizes Lemma 3.2 of Choe
\cite{Choe:double}, which claims $\Candle(0)$ using (in the proof) the
same hypotheses when $\kappa = 0$.  However, the argument given there
omits many details about reflection from a convex surface.  (Which is
thus concave from the other side as we describe it.)  We will prove
Proposition~\ref{p:mgunther} in Section~\ref{s:mgunther}.

\section{Linear programming and optimal transport}
\label{s:linear}

\subsection{A linear model for isoperimetric problems}
\label{s:linmod}

\subsubsection{Feasibility}

In this section we abstract the results of Section~\ref{s:etendue} and
\ref{s:core} into a linear programming model.

Assume that $\Omega$ is an $n$-manifold with boundary, with unique geodesics,
and with curvature at most $\kappa$.  Let $V=|\Omega|$ and $A=|\dOmega|$.
Then equations~\eqref{e:etendue}, \eqref{e:margin} \eqref{e:croke},
\eqref{e:prince}, and \eqref{e:teufel} show (after symmetrization in $\alpha$
and $\beta$) that $\mu=\mu_\Omega$ is a solution to the following infinite
linear programming problem.

\begin{linprog} Given $n$, $\kappa$, $A$, and $V$, let
\[ z(\theta) = \frac{\omega_{n-2}\sin(\theta)^{n-1}}{n-1}. \]
Is there a positive measure $\mu(\ell,\alpha,\beta)$ on $\R_{\ge 0} \times
[0,\pi/2)^2$, which is symmetric in $\alpha$ and $\beta$, and such that
\begin{align}
\alpha_*(\mu) = \int_{\ell,\beta} \back \dd \mu
    &= A\sdd z(\alpha) \label{e:lpmargin} \\
\int_{\ell,\alpha,\beta} \back s_{n,\kappa}(\ell)\sec(\alpha)\sec(\beta)
    \sdd\mu &\le A^2 \label{e:lpcroke} \\
\int_{\ell,\alpha,\beta} \back s_{n,\kappa}^{(-1)}(\ell)
    \big(\sec(\alpha)+\sec(\beta)\big) \sdd\mu &\le 2AV \label{e:lpprince} \\
\int_{\ell,\alpha,\beta} \back s_{n,\kappa}^{(-2)}(\ell) \sdd\mu
    &\le V^2 \label{e:lpteufel} \\
\int_{\ell,\alpha,\beta} \back \ell \sdd\mu &= \omega_{n-1} V? \label{e:lpvol} 
\end{align}
\eatline \label{lp:basic} \end{linprog}

(We could have written Problem~\ref{lp:basic} without symmetrization in
$\alpha$ and $\beta$.  It would have been equivalent, but more complicated.)

Since our ultimate goal is to prove a lower bound for $|\dOmega|$, we
want to show that given $n$, $\kappa$, and $V$, Problem \ref{lp:basic} is
infeasible for values of $A$ that are too low.  As usual in linear
programming, we will profit from stating a dual problem.

\begin{linprog}[Dual to Problem~\ref{lp:basic}] Given $n$, $\kappa$,
$A$, and $V$, are there numbers $a,b,c \ge 0$ and $d \in \R$ and a
continuous function $f:[0,\pi/2) \to \R$ such that
\begin{gather}
\begin{multlined}[b]
\back a s_{n,\kappa}(\ell)\sec(\alpha)\sec(\beta)
+ b s_{n,\kappa}^{(-1)}(\ell)\big(\sec(\alpha) + \sec(\beta)\big) \\
+ c s_{n,\kappa}^{(-2)}(\ell) - d\ell + f(\alpha) + f(\beta) \ge 0
\end{multlined} \label{e:lpdual1} \\
aA^2 + 2bAV + cV^2 - d\omega_{n-1}V + 2A\!\!\int_0^{\pi/2}
    \back \back f(\alpha) \sdd z(\alpha) < 0 \label{e:lpdual2}
\end{gather}
for all $(\alpha,\beta,\ell) \in [0,\pi/2)^2 \times \R_{\ge 0}$?
\label{lp:dual} \end{linprog}

(Note that the constant $d$ can have either sign.  We subtract it so that
it will be positive in actual usage.)

We will discuss in what sense Problem~\ref{lp:dual} is dual to
Problem~\ref{lp:basic}, and the consequences of this duality, in
Section~\ref{s:general}.   For now, we will concentrate on sufficient
criteria to prove our main theorems.  Problem~\ref{lp:dual} is strong
enough to prove Theorem~\ref{th:positive}.  In Sections~\ref{s:negative}
and \ref{s:other}, we will state other linear programming problems to
handle our other results stated in Section~\ref{s:intro}.

In the rest of this section (Section~\ref{s:linmod}), including in the
statements of the lemmas, we fix $V$, $n$, and $\kappa$, but not $A$.

\begin{lemma} Let $a, b, c \ge 0$, let $d \in \R$, and let
\begin{eq}{e:cost} E(\ell,\alpha,\beta)
    = a s_{n,\kappa}(\ell)\sec(\alpha)\sec(\beta)
    + b s_{n,\kappa}^{(-1)}(\ell)\big(\sec(\alpha)+\sec(\beta)\big)
    + c s_{n,\kappa}^{(-2)}(\ell) - d\ell. \end{eq}
Let $f:[0,\pi/2) \to \R$ be a continuous function such that $\int_0^{\pi/2}
f(\alpha) \sdd z(\alpha)$ is absolutely convergent.  If
\begin{eq}{e:confirm}
F(\ell,\alpha,\beta) \defeq E(\ell,\alpha,\beta) + f(\alpha) + f(\beta) \ge 0,
\end{eq}
then Problem~\ref{lp:basic} is infeasible for those $A \ge 0$ such that
\[ P(A) \defeq aA^2 + 2bAV + cV^2 - d\omega_{n-1}V + 2A\!\!\int_0^{\pi/2}
    \back f(\alpha) \sdd z(\alpha) < 0. \]
Explicitly, if $P(A)$ has two real roots, let the roots be $A_0 < A_1$;
if $P(A)$ is linear, let $A_1$ be its root and let $A_0 = -\infty$.
Then $A \in (A_0,A_1) \cap [0,\infty)$ is infeasible.
\label{l:dual} \end{lemma}

We introduce some terminology which will be justified
in Section~\ref{s:optrans}.  $E(\ell,\alpha,\beta)$ is a \emph{cost
function}, $f(\alpha)$ is a \emph{potential}, and $F(\ell,\alpha,\beta)$
is an \emph{adjusted cost function}.

\begin{proof} Except for a change of variables, the proposition
is the assertion that if Problem~\ref{lp:dual} is feasible, then
Problem~\ref{lp:basic} is infeasible.  More explicitly:  For 
any $a, b, c \ge 0$, $d \in \R$, and suitable $f:[0,\pi/2) \to \R$,
we can combine the relations in Problem~\ref{lp:basic} to produce
a formula of the form
\begin{eq}{e:FP}
\int_{\ell,\alpha,\beta} \back \back F(\ell,\alpha,\beta) d\mu \le P(A).
\end{eq}
If the integrand $F(\ell,\alpha,\beta)$ is non-negative while the upper
bound $P(A)$ is strictly negative, then the measure $\mu$ cannot exist.
\end{proof}

At this point, there is a potential difference between solving a geometric
isoperimetric problem and solving a linear programming model for one.
(Recall that Theorem~\ref{th:subsume} promises the latter.)  Geometrically,
the set of possible values of $A$ in each of our isoperimetric problems must
be an open or closed ray in $\R_+$.  Thus, if we apply Lemma~\ref{l:dual}
to exclude $A \in (A_0,A_1)$, then all values of $A \in [0,A_1)$ are
geometrically impossible even if $A_0 \ge 0$.  Problem~\ref{lp:basic}
has the same property because we can think of $A = A(\mu)$ as a function
of $\mu$, and we can increase $A$ by adding measure to $\mu$ at $\ell=0$.
But this is less clear for Problem~\ref{lp:extend} which we will need later;
in any case, a single use of Lemma~\ref{l:dual} need not satisfy $A_0 < 0$.
The simplest way to establish Theorem~\ref{th:subsume} is to obtain a
negative value of $A_0$ and a sharp value of $A_1$ in Lemma~\ref{l:dual}.
This is the case if and only if
\[ P(0) = cV^2 - d\omega_{n-1}V < 0, \]
which simplifies to
\begin{eq}{e:clean}
cV < d\omega_{n-1}.
\end{eq}
We will attain the condition \eqref{e:clean} for Problem~\ref{lp:basic}.
In our treatment of Problem~\ref{lp:extend}, we will take a slightly more
complicated approach.

\subsubsection{Optimality}
\label{s:optimal}

Suppose that for some value $A_1$, we find a solution $\mu$ to
Problem~\ref{lp:basic}, and we find $(a,b,c,d,f)$ in Lemma~\ref{l:dual}
with $P(A_1) = 0$ and $P'(A_1) > 0$; and suppose that \eqref{e:clean} also
holds.  Then the two solutions are an \emph{optimal pair}.  The existence
of $(a,b,c,d,f)$ shows that $A = A_1$ is the smallest feasible value in
Problem~\ref{lp:basic}; the existence of $\mu$ shows that $A = A_1$ is
the smallest infeasible value in Problem~\ref{lp:dual}.

There is a simple test of whether $\mu$ and $(a,b,c,d,f)$ are an optimal
pair.  If they are, then the conditions \eqref{e:FP}, $F \ge 0$, and $P(A)
= 0$ tell us that $\mu$ is supported on the zero locus of the adjusted cost
$F$.  On the other hand, if $(a,b,c,d,f)$ satisfies both Lemma~\ref{l:dual}
and \eqref{e:clean}, if $\mu$ is a solution to Problem~\ref{lp:basic}
that is supported on the zero locus of $F$, and if we happen to know that
all inequalities in Problem~\ref{lp:basic} are equalities, then we can
calculate that $P(A) = 0$.

If in addition $\mu = \mu_\Omega$ for an admissible domain $\Omega$,
then $A = |\dOmega|$ is the sharp isoperimetric value.  Recall that we
plan to prove that $\Omega = B_{n,\kappa}$ is an isoperimetric minimizer.
First, since this $\Omega$ is convex and has constant curvature $\kappa$,
Lemma~\ref{l:core} tells us that all inequalities in Problem~\ref{lp:basic}
are indeed equalities.   Second, the geodesics of this $\Omega$ have the
property that $\alpha = \beta$ and that $\ell = h(\alpha)$ is a function
of $\alpha$.  If we combine these properties with the assumption that
$\mu_\Omega$ is part of an optimal pair and is thus supported on the zero
locus of $F$, then we can solve for $f$, once we know $(a,b,c,d)$.
Together with the rest of the discussion in this section, we obtain the
following sufficient criterion.

\begin{lemma} Suppose that $a,b,c \ge 0$ and $d \in \R$ are numbers
that satisfy \eqref{e:clean}, and that $h:[0,\pi/2) \to \R_{\ge 0}$ is
a continuous function.  Let $\mu$ be the unique measure that satisfies
\eqref{e:lpmargin} for some $A = A(\mu)$ and that is supported on the set
$(h(\alpha),\alpha,\alpha)$, and let
\begin{eq}{e:potential}
f(\alpha) = -\frac{E(h(\alpha),\alpha,\alpha)}2.
\end{eq}
If $\mu$ is a solution to Problem~\ref{lp:basic}, and if $f$ satisfies
\eqref{e:confirm} for the same $A$, then $\mu$ and $(a,b,c,d,f)$ are an
optimal pair.  If in addition $\mu = \mu_\Omega$ for an admissible domain
$\Omega$, then $A = |\dOmega|$ is the sharp isoperimetric value.
\label{l:sharp} \end{lemma}

Lemma~\ref{l:sharp} is the basis for our proof of Theorem~\ref{th:positive}
and the corresponding part of Theorem~\ref{th:subsume}.   We will use
similar reasoning to prove Theorems~\ref{th:negative} and \ref{th:multiple}.
The calculations will be organized as follows.  We temporarily assume the
conclusion, that $\mu_\Omega$ is optimal when $\Omega = B_{n,\kappa}(r)$.
This yields the dependence $\ell = h(\alpha)$.  If our construction were
to work, the adjusted cost $F(\ell,\alpha,\beta)$ would attain a minimum
of $0$ at $(h(\alpha),\alpha,\alpha)$.  Thus, we can solve for $a$, $b$,
$c$, and $d$ by applying a derivative test to the cost $E$ or the adjusted
cost $F$, namely,
\begin{eq}{e:consist}
\dbyd{F}{\ell}(\ell,\alpha,\alpha) = \dbyd{E}{\ell}(\ell,\alpha,\alpha) = 0
\end{eq}
when $\ell = h(\alpha)$ and $0 \le \ell \le 2r$.

Having calculated $a$, $b$, $c$, and $d$, which determine
$E(\ell,\alpha,\beta)$, \eqref{e:potential} tells us $f(\alpha)$.   The
remaining hard part of the proof is then to confirm \eqref{e:confirm}.
We will carry out these calculations in Section~\ref{s:main}.

Our approach to solving Problem~\ref{lp:basic} as outlined in this section
may seem both lucky and creative.   It is indeed lucky, in the sense that
\eqref{e:consist} is an equality of functions used to solve for four numbers;
it only has solutions when $n \in \{2,4\}$.  In Section~\ref{s:general},
we will argue that solving Problem~\ref{lp:basic} follows the precepts of
linear programming and optimal transport with fairly little creativity.

\subsection{Generalities}
\label{s:general}

As explained at the end of Section~\ref{s:optimal}, this section
only provides context and is not needed for the proofs of our
results.

\subsubsection{Linear programming}
\label{s:linprog}

The genesis of linear programming is a structure theorem for finite
systems of linear equalities and inequalities due to Farkas and Minkowski
\cite{Farkas:einfachen,Minkowski:zahlen,Kjeldsen:goals}.

\begin{theorem}[Farkas-Minkowski] Let $x = \{x_i\}$ be a finite list of real
variables, and a finite system $L$ of linear inequalities and equalities,
given by two matrices $A$ and $B$:
\[ \sum_i A_{j,i} x_i \le a_j \qquad \sum_i B_{k,i} x_i = b_k. \]
Then:
\begin{enumerate}
\item The system $L$ is infeasible if and only if some linear combination
of the form
\begin{eq}{e:farkas}
\sum_{i,j} y_j A_{j,i} x_i + \sum_{i,k} z_k B_{k,i} x_i
    \le \sum_j y_j a_j + \sum_k z_k b_k \qquad \forall j, y_j \ge 0
\end{eq}
simplifies to the falsehood $0 \le -1$.  (Or $0 \le c$ for some constant
$c < 0$.)
\item A linear bound $\sum_i c_i x_i \le c$ holds for solutions to
$L$ if and only if it is expressible in the form \eqref{e:farkas}.
\item If $L$ is feasible and $\sum_i c_i x_i$ is bounded on its solution
set, then it has a maximum $c$, which is also the minimum of the right side of
\eqref{e:farkas} subject to the constraint that the left side simplifies
to $\sum_i c_i x_i$.
\end{enumerate}
\label{th:farkas} \end{theorem}

The coefficients $\{y_j\}$ and $\{z_k\}$, subject to the constraints in
one of the cases of Theorem~\ref{th:farkas}, is then a \emph{dual system}
$L^*$ to $L$.  If $\{x_i\}$ is feasible for $L$ and attains a value of $c$
for the objective $\sum_i c_i x_i$, and if $\{y_j\}$ and $\{z_k\}$ are
feasible for $L^*$ and attain the same $c$, then they are an \emph{optimal
pair}; each half of the pair proves that the other half is optimal.  Thus,
case 3 of Theorem~\ref{th:farkas} says that every maximization problem in
finite linear programming with a finite maximum can be solved by finding
an optimal pair.

We cannot directly apply Theorem~\ref{th:farkas} to Problem~\ref{lp:basic}
because it is an infinite{\hyp}dimensional problem.  The theorem still holds
in infinite dimensions, or in finite dimensions with infinitely many
inequalities, with an extra hypothesis such as compactness.  We do not
know a simple way to make Problem~\ref{lp:basic} compact, but we will find
optimal pairs anyway.

The standard notion of an optimal pair from Theorem~\ref{th:farkas}
is not exactly the same as that in Section~\ref{s:optimal}, because
Problem~\ref{lp:basic} is nonlinear in the variable $A$.  However,
the two concepts are analogous.  Indeed, we can change Problem~\ref{lp:basic}
slightly to make $A$ a linear variable, as follows.  First, 
we switch to the measure $\hmu \defeq \mu/A$ and divide through by
$A$.   Then the first three relations, \eqref{e:lpmargin}, \eqref{e:lpcroke},
and \eqref{e:lpprince}, are all linear in the variables $\hmu$ and $A$.
The last two relations, \eqref{e:lpteufel} and \eqref{e:lpvol},
now have a factor of $1/A$, but we can convert them to these two equations:
\[ \int_{\ell,\alpha,\beta} \back (\omega_{n-1}s_{n,\kappa}^{(-2)}(\ell)-V\ell)
    \sdd\hmu \le 0  \qquad \int_{\ell,\alpha,\beta} \back \ell
    \sdd\hmu \ge \frac{\omega_{n-1} V}A. \]
These relations are moderately weaker than \eqref{e:lpteufel} and
\eqref{e:lpvol}, but switching to them is nearly equivalent to assuming
the condition \eqref{e:clean}.  Although the second equation is still
nonlinear, it is a convex relation between $A$ and $\hmu$; it can be
expressed by a family of linear relations.

Finally, linear programming in infinite dimensions usually involves
topological vector spaces.  For instance, the measure $\mu$ in
Problem~\ref{lp:basic} lies in a space of Borel measures.  This raises the
question of the appropriate regularity of the dual variable $f(\alpha)$.
Because of \eqref{e:lpmargin}, the function $f(\alpha)$ could in principle
be integrable rather continuous; $f(\alpha)\sdd z(\alpha)$ could even
be replaced a Borel measure.  However, Proposition~\ref{p:convex} from
optimal transport theory tells us that an optimal $f(\alpha)$ is continuous.

\subsubsection{Optimal transport}
\label{s:optrans}

We can interpret Problem~\ref{lp:basic} as an optimal transport problem.
See Villani~\cite[Ch.3-5]{Villani:oldnew} for background material
on optimal transport.  Following Villani, we assume that $A\sdd
z(\alpha)$ is a distribution of boulangeries and $A \sdd z(\beta)$ is
a distribution of caf\'es.  Moreover, for each boulangerie $\alpha$ and
caf\'e $\beta$, there are a range of possible roads parametrized by $\ell$.
By \eqref{e:lpmargin}, $\mu$ is a transport of baguettes\footnote{Even though
in Section~\ref{s:etendue}, we transported photons.} from the boulangeries
to the caf\'es.  Problem~\ref{lp:basic} then asks whether the transport
is feasible given the constraints that we must pay separate road tolls
in Polish zlotys \eqref{e:lpcroke}, Czech korunas \eqref{e:lpprince}, and
Hungarian forints \eqref{e:lpteufel}; and given an exact labor requirement
\eqref{e:lpvol} (neither more nor less).  Strictly speaking, this is
a feasible transport problem rather than an optimal transport problem,
but we can convert it to optimal transport.

The function $E(\ell,\alpha,\beta)$ defined in equation \eqref{e:cost}
is a natural reduction of all four resource limits into one combined
cost function, which we can then optimize to test feasibility.  In the
economics interpretation, the coefficients $a$, $b$, and $c$ are currency
conversions, while $d$ is a wage rate. The last term $d\ell$ is naturally
subtracted if employment is the goal of the program and thus a negative cost.
Certainly if any choice of $a,b,c,d$ yields a cost function $E$ such that
\begin{eq}{e:feastrans}
\int_{\ell,\alpha,\beta} E(\ell,\alpha,\beta) d\mu
    \le aA^2 + 2bAV + cV^2 - d\omega_{n-1}V
\end{eq}
is infeasible, then the original multi-resource transport problem is
also infeasible.  We won't try to prove the converse for all $n$
and $\kappa$: that if Problem~\ref{lp:basic} is infeasible, then there
exist $(a,b,c,d)$ such that \eqref{e:feastrans} is also infeasible.

Even so, once $a,b,c,d$ are chosen, Problem~\ref{lp:basic} reduces
to just \eqref{e:lpmargin} and \eqref{e:feastrans}, which is a nearly
standard optimal transport problem.  The two differences are:
\begin{enumerate}
\item We have a choice of ``roads" parametrized by $\ell$.  Given a
scalar cost, we can convert it to a standard optimal transport problem if
we choose the most efficient road for each pair $(\alpha,\beta)$ and let
the cost be
\[ E(\alpha,\beta) = \min_\ell E(\ell,\alpha,\beta). \]
\item The transport $\mu$ does not usually have to be symmetric
in $\alpha$ and $\beta$.   We can live without this constraint because
Problem~\ref{lp:basic} is itself symmetric in $\alpha$ and $\beta$, if we
add the relation $\beta_*(\mu) = A\sdd z(\beta)$, which is the other half
of \eqref{e:margin}.  We can symmetrize any solution using
\[ \hmu(\ell,\alpha,\beta) \defeq \frac{\mu(\ell,\alpha,\beta)
    + \mu(\ell,\beta,\alpha)}2. \]
\end{enumerate}

Having fixed $a,b,c,d$, the remaining dual variable in Problem~\ref{lp:dual}
is $f(\alpha)$.   Its sole constraint is \eqref{e:confirm}.
In optimal transport terminology, $f(\alpha)$ is known as a \emph{Kantorovich
potential}.  We will call the left side, $F(\ell,\alpha,\beta)$, the
\emph{adjusted cost function}.  In standard optimal transport, we would
have two potentials $f(\alpha)$ and $g(\beta)$ satisfying the equation
\begin{eq}{e:fg}
E(\ell,\alpha,\beta) + f(\alpha) + g(\beta) \ge 0.
\end{eq}
But, just as symmetry is optional in Problem~\ref{lp:basic}, it is
also optional in Problem~\ref{lp:dual}; we can symmetrize a solution
to make $f = g$.

\begin{proposition} An optimal potential $f(\alpha)$ in Problem~\ref{lp:dual}
is a convex function of $\sec(\alpha)$ and therefore continuous.
\label{p:convex} \end{proposition}

Proposition~\ref{p:convex} is a standard type of result in optimal
transport theory.  A potential that satisfies an equation such
as \eqref{e:costconv} below is called \emph{cost convex}.

\begin{proof} We assume two potentials $f(\alpha)$ and $g(\beta)$.
In an asymmetric variation of Problem~\ref{lp:dual}, they are chosen
to minimize
\[ \int_0^{\pi/2} \back \back f(\alpha) \sdd z(\alpha) +
    \int_0^{\pi/2} \back \back g(\beta) \sdd z(\beta). \]
For each fixed $g(\beta)$, we can minimize this integral subject to the
constraint \eqref{e:fg} by choosing
\begin{eq}{e:costconv}
f(\alpha) = \sup_{\ell,\beta} \big[ - E(\ell,\alpha,\beta) - g(\beta) \big].
\end{eq}
For each fixed value of $\ell$ and $\beta$, the supremized function on
the right side is linear in $\sec(\alpha)$ by \eqref{e:cost}.  It follows
that $f(\alpha)$ is convex in $\sec(\alpha)$ and thus continuous; the same
is true of $g(\beta)$.  If this asymmetric optimization yields $f \ne g$,
then their average $(f+g)/2$ has all of the desired properties.
\end{proof}

\section{Proofs of the main results}
\label{s:main}

In this section, we will complete the proofs Theorems~\ref{th:positive} and
\ref{th:negative}, picking up from Section~\ref{s:linmod}.  Up to rescaling,
we can assume that $\kappa \in \{-1,0,1\}$.  We first explicate the condition
\eqref{e:consist}, which we will use to check whether Problem~\ref{lp:basic}
has any hope of producing a sharp isoperimetric inequality, and to calculate
the parameters in Lemma~\ref{l:dual}.

If $\Omega = B_{n,\kappa}(r)$, then the length of a geodesic chord that makes
an angle of $\alpha$ from the normal to $\dOmega$ is given by the relation
\begin{eq}{e:chord}
\cos(\alpha) = T_{\kappa,r}(\ell) \defeq \begin{cases}
    \dfrac{\tan(\ell/2)}{\tan(r)} &\mbox{ if $\kappa=1$} \\[2ex]
    \dfrac{\ell}{2r} &\mbox{ if $\kappa=0$.} \\[2ex]
    \dfrac{\tanh(\ell/2)}{\tanh(r)} &\mbox{ if $\kappa=-1$}
\end{cases}
\end{eq}
Equation \eqref{e:chord} thus gives us the function $\ell = h(\alpha)$
in the statement of Lemma~\ref{l:sharp}.

We combine equations \eqref{e:cost} and \eqref{e:consist} to obtain
\[ \dbyd{E}{\ell}(\ell,\alpha,\alpha) 
    = a\frac{s'_{n,\kappa}(\ell)}{\cos(\alpha)^2} 
    + 2b\frac{s_{n,\kappa}(\ell)}{\cos(\alpha)}
    + cs_{n,\kappa}^{(-1)}(\ell) - d = 0. \]
Combining with \eqref{e:chord}, we obtain
\begin{eq}{e:check}
a\frac{s_{n,\kappa}'(\ell)}{T_{\kappa,r}(\ell)^2} +
    2b\frac{s_{n,\kappa}(\ell)}{T_{\kappa,r}(\ell)}
    + cs_{n,\kappa}^{(-1)}(\ell) - d = 0.
\end{eq}
Again, \eqref{e:check} is an equation for the coefficients $a,b,c,d$
that should hold for $0 \le \ell \le 2r$.  In each case, the coefficients
will be unique up to rescaling by a positive real number.  Note that the
factor of $\tan(r)$, $r$, or $\tanh(r)$ that appears in $T_{\kappa,r}$
factors of out of the question of whether there is a solution, since this
factor can be absorbed into the constants $a$ and $b$.

Our proofs in this section follow a set pattern:
\begin{enumerate}
\item Working either from Problem~\ref{lp:basic} or
Problem~\ref{lp:extend}, and their dual problems, calculate $(a,b,c,d)$
using \eqref{e:check}.

\item Change variables from $\alpha$ and $\beta$ to $x$ and $y$
using \eqref{e:zxy}, \eqref{e:pxy}, or \eqref{e:nxy}.  Calculate the cost
$E(\ell,x,y)$, the potential $f(x)$, and the adjusted cost $F(\ell,x,y)$
in the new variables.

\item Using calculus methods, establish that the adjusted cost
$F(\ell,x,y)$ is non-negative, according to \eqref{e:confirm}.  This will
fulfill the hypotheses of Lemma~\ref{l:dual} when $\kappa \ge 0$, or its
equivalent when $\kappa < 0$, and finish the proof.  In the hardest two
cases ($n=4$ and $\kappa \ne 0$), this step depends crucially on symbolic
algebra software.
\end{enumerate}
Also, we abbreviate $s = s_{n,\kappa}$ throughout.

Here are two general remarks about dimension $n=2$.   First, for every
value of $\kappa$, there is a separation \eqref{e:sep} in this dimension.
This means that we could have proved the results with a simpler measure
$\mu(\ell,\alpha)$ that depends on only one angle, in the spirit of
Section~\ref{s:illum}.  Second, $a=0$ when $n=2$, so we can immediately
accept $A$ as a linear variable in Problem~\ref{lp:basic} or \ref{lp:extend}.

\subsection{Weil's and Croke's theorems}
\label{s:zero}

This case is a warm-up to the more difficult cases with $\kappa \ne 0$.
We introduce the change of variables
\begin{eq}{e:zxy}
(x,y) \defeq \Big(\frac{\sec(\alpha)}{r},
    \frac{\sec(\beta)}{r}\Big)
\end{eq}
in place of $\alpha$ and $\beta$.  We will give them the range $x,y \in
\R_{\ge 0}$.  By abuse of notation, we can change variables without changing
the names of functions; for example, we can write
\[ E(\ell,\alpha,\beta) = E(\ell,\alpha(x),\beta(y)) = E(\ell,x,y). \]

If $\kappa = 0$, then
\begin{align*}
s(\ell) &= \ell^{n-1}, & s'(\ell) &= (n-1)\ell^{n-2}, \\
s^{(-1)}(\ell) &= \frac{\ell^{n}}n, & s^{(-2)}(\ell)
    &= \frac{\ell^{n+1}}{n(n+1)}.
\end{align*}
Equation \eqref{e:check} becomes
\[ 4(n-1)r^2 a \ell^{n-4} + 4r b \ell^{n-2} + \frac{c\ell^n}n - d = 0 .\]
Obviously this has solutions if $n \in \{2,4\}$ and not otherwise;
this point was known to Croke (personal communication).

When $n=2$, the solution is
\[ a = 0, \qquad b = \frac{1}r, \qquad c = 0, \qquad d = 4. \]
From \eqref{e:cost}, we thus obtain
\[ E(\ell,\alpha,\beta) = \frac{\ell^2(\sec(\alpha) + \sec(\beta))}{2r}
    - 4\ell \]
Then \eqref{e:potential} and \eqref{e:chord} give us the potential
\[f(\alpha) = -\frac{E(2r\cos(\alpha),\alpha,\alpha)}2 = 2r\cos(\alpha) .\]
Then the adjusted cost \eqref{e:confirm} separates as
\begin{eq}{e:sep}
F(\ell,\alpha,\beta) = G(\ell,\alpha) + G(\ell,\beta) 
\end{eq}
with 
\[ G(\ell,\alpha) = \frac{\ell^2\sec(\alpha)}{2r} - 2\ell + 2r\cos(\alpha)
    = \frac{2r\cos(\alpha) - \ell}{2r\cos(\alpha)} \ge 0. \]
Thus $F \ge 0$, which establishes Weil's theorem.

When $n=4$, the solution to \eqref{e:check} is
\[ a = \frac{1}{r^2}, \qquad b = 0, \qquad c = 0, \qquad d = 12. \]
These coefficients plainly satisfy condition \eqref{e:clean}.  The cost
function is
\[ E(\ell,\alpha,\beta) = \frac{\ell^3\sec(\alpha)\sec(\beta)}{r^2}
    - 12\ell, \]
the potential is
\[f(\alpha) = -\frac{E(2r\cos(\alpha),\alpha,\alpha)}2 = 8r\cos(\alpha), \]
and their sum is
\[ F(\ell,\alpha,\beta) = \frac{\ell^3\sec(\alpha)\sec(\beta)}{r^2} - 12\ell
    + 8r(\cos(\alpha) + \cos(\beta)). \]
Using the change of variables \eqref{e:zxy},
\[ F(\ell,x,y) = \ell^3xy - 12\ell + \frac8x + \frac8y. \]
We want to show that $F \ge 0$.  For each fixed value of $xy$, $F$ is
minimized when $x=y$.  We can then calculate
\[ F(\ell,x,x) = \ell^3x^2 - 12\ell + \frac{16}x
    = \frac{(\ell x + 4)(\ell x - 2)^2}x \ge 0. \]
This establishes Croke's theorem.

Following the comments after the proof of Theorem~\ref{th:point} in
Section~\ref{s:illum}, our proof of Croke's theorem is only superficially
different from Croke's proof.  The extra point here is that Croke's
theorem (and Weil's theorem along with it) hold in Model~\ref{lp:basic},
which establishes part of Theorem~\ref{th:subsume}.

\subsection{The positive case}
\label{s:positive}

In this section we will establish Theorem~\ref{th:positive}.  We will
let $\kappa = 1$, but before we do that, we note that $\kappa = 0$ is a
limiting case of $\kappa > 0$.   Section~\ref{s:zero} established that a
sharp result in the case $\kappa = 0$ is only possible when $n \in \{2,4\}$,
this justifies the same restriction in Theorem~\ref{th:positive}.

We use the change of variables
\begin{eq}{e:pxy}
(x,y) \defeq \Big(\frac{\sec(\alpha)}{\tan(r)},
    \frac{\sec(\beta)}{\tan(r)}\Big)
\end{eq}
with the range $x,y \in \R_{\ge 0}$.  Note that equation \eqref{e:chord}
simplifies to
\begin{eq}{e:xchord}
\tan(\frac{\ell}2) = \frac1x.
\end{eq}

\subsubsection{Dimension 2}

In dimension $n=2$,
\begin{align*}
s(\ell) &= \sin(\ell), & s'(\ell) &= \cos(\ell), \\
s^{(-1)}(\ell) &= 1-\cos(\ell), & s^{(-2)}(\ell) &= \ell-\sin(\ell) 
\end{align*}
when $\ell < \pi$, and
\[ s(\ell) = 0, \qquad s^{(-1)}(\ell) = 2, \qquad s^{(-2)} = 2\ell-\pi \]
for $\ell \ge \pi$.  Equation~\eqref{e:check}, with \eqref{e:chord}, becomes
\[ \frac{a\tan(r)^2\cos(\ell)}{\tan(\ell/2)^2}
    + \frac{2b\tan(r)\sin(\ell)}{\tan(\ell/2)} + c(1-\cos(\ell)) - d = 0. \]
The solution is
\[ a = 0, \qquad b = \frac{1}{\tan(r)}, \qquad c = 2, \qquad d = 4. \]
In the variables \eqref{e:pxy}, the cost function \eqref{e:cost} is
\[ E(\ell,x,y) = (1-\cos(\ell))(x+y) - 2\sin(\ell) - 2\ell, \]
for $\ell \le \pi$, and is constant in $\ell$ for $\ell \ge \pi$:
\begin{eq}{e:lconst}
E(\ell,x,y) = E(\pi,x,y) \qquad \forall \ell \ge \pi.
\end{eq}
Using \eqref{e:xchord}, the potential \eqref{e:potential} becomes
\[ f(x) = \ell = 2\arctan(\frac{1}x). \]
The adjusted cost \eqref{e:confirm} again separates according to
\eqref{e:sep}, where this time
\[ G(\ell,x) = (1-\cos(\ell))x - \sin(\ell) - \ell + 2\arctan(\frac{1}x). \]
We can minimize $G$ with the derivative test either in $\ell$ or in $x$.
The latter is slightly simpler and gives us
\[ \dbyd{G}{x}(\ell,x) = \frac{x^2(1-\cos(\ell))-(\cos(\ell)+1)}{x^2+1}. \]
We learn that $\del G/\del x$ crosses $0$ exactly once, when $x$ and $\ell$
satisfy \eqref{e:xchord}; this is therefore the minimum of $G$ for each
fixed $\ell$.  Since the relation \eqref{e:xchord} is used to define the
potential $f(x)$, it is automatic that this minimum value is $0$; the
substitution $x = 1/\tan(\ell/2)$ also establishes it.  Thus $G(\ell,x)
\ge 0$, which confirms \eqref{e:confirm} and establishes the $n=2$ case
of Theorem~\ref{th:positive}.

\subsubsection{Dimension 4}

In dimension $n=4$,
\begin{align*}
s(\ell) &= \sin(\ell)^3, \\
s'(\ell) &= 3\sin(\ell)^2\cos(\ell), \\
s^{(-1)}(\ell) &= \frac{\cos(\ell)^3-3\cos(\ell)+2}3, \\
s^{(-2)}(\ell) &= \frac{6\ell-\sin(\ell)^3-6\sin(\ell)}{9}
\end{align*}
when $\ell < \pi$, and
\[ s(\ell) = 0,\quad s^{(-1)}(\ell) = \frac43,
    \quad s^{(-2)}(\ell) = \frac{4\ell-2\pi}3 \]
when $\ell \ge \pi$.  Equation~\eqref{e:check} becomes
\[ \frac{3a\tan(r)^2\cos(\ell)\sin(\ell)^2}{\tan(\ell/2)^2}
    + \frac{2b\tan(r)\sin(\ell)^3}{\tan(\ell/2)}
    + \frac{c(\cos(\ell)^3-3\cos(\ell)+2)}{3} - d = 0. \]
The solution is
\[ a = \frac{1}{\tan(r)^2}, \qquad b = \frac{3}{\tan(r)},
    \qquad c = 9, \qquad d = 12. \]
The clean optimality condition \eqref{e:clean} becomes
\[ 9V < 12\omega_3 = 24\pi^2. \]
Since $V$ is at most the volume of a hemisphere, we have
\[ 9V < \frac{9\omega_4}2 = 12\pi^2. \]
Thus \eqref{e:clean} holds.

\begin{figure}[htb]
\begin{center}
\includegraphics[width=.7\columnwidth]{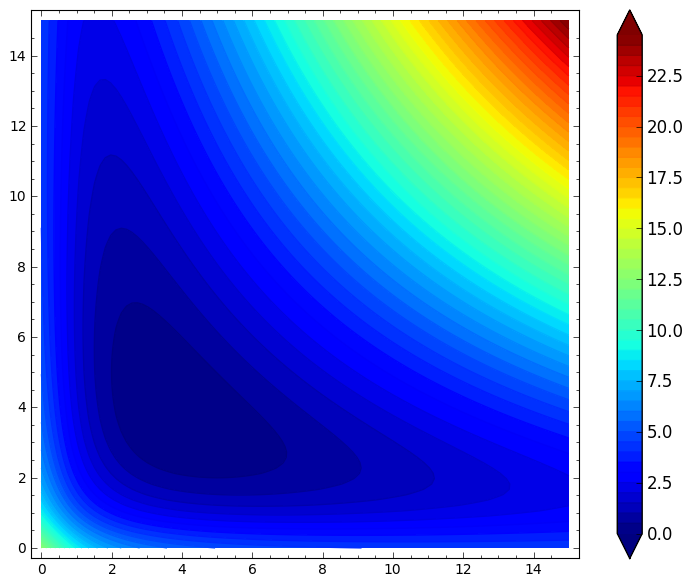}
\end{center}
\caption{The slice $F(\pi/6,x,y)$ in the case $\kappa = 1$.}
\label{f:Fpi6} \end{figure}

The cost function \eqref{e:cost} is
\[ E(\ell,x,y) = \sin(\ell)^3xy + (\cos(\ell)^3-3\cos(\ell)+2)(x+y)
    -\sin(\ell)^3-6\sin(\ell) - 6\ell. \]
for $\ell \le \pi$, while once again $E$ is constant in $\ell$ for $\ell
\ge \pi$, as in \eqref{e:lconst}.  The potential from \eqref{e:potential}
and \eqref{e:xchord} is
\[ f(x) = 6\arctan(\frac1x) + \frac{2x}{x^2+1}. \]
We will include the values $x=0$ and $y=0$ in our calculations, so it is
helpful to recall that
\[ \arctan(\frac1x) = \frac{\pi}2 - \arctan(x). \]
The adjusted cost \eqref{e:confirm} is
\begin{multline}
F(\ell,x,y) = \sin(\ell)^3xy + (\cos(\ell)^3-3\cos(\ell)+2)(x+y)
    - \sin(\ell)^3 - 6\sin(\ell) - 6\ell \\ + 6\pi - 6\arctan(x)
    + \frac{2x}{x^2+1} - 6\arctan(y) + \frac{2y}{y^2+1} .
\label{e:Fp4} \end{multline}
The remainder of the proof of Theorem~\ref{th:positive} is given by the
following lemma.  Although the lemma is evident from contour plots (\eg,
Figure~\ref{f:Fpi6}), the authors found it surprisingly tricky to find
a rigorous proof.

\begin{lemma} The function $F(\ell,x,y)$ on $[0,\pi] \times \R_{\ge 0}^2$
given by \eqref{e:Fp4} is non-negative, and vanishes only when
\[ x = y = \frac{1}{\tan(\ell/2)}. \]
\label{l:tech4p} \end{lemma}

\begin{proof}
We will use these immediate properties of the potential $f(x)$:
\[ f(0) = 3\pi, \qquad f(x) > 0. \]

We first check the non-compact direction of the domain of $F$.  There exists
a constant $k > 0$ such that
\[ s^{(-1)}(\ell) \ge k\ell^4. \]
(Because $\ell^4/s^{(-1)}(\ell)$ is continuous on $[0,\pi]$ and therefore
bounded.  In fact
\[ k = \frac{s^{(-1)}(\pi)}{\pi^4} = \frac{4}{3\pi^4} \]
works.)  Thus
\begin{align*}
F(\ell,x,y) &= s(\ell)xy + 3s^{(-1)}(\ell)(x+y) + 9s^{(-2)}(\ell) - 12\ell 
    + f(x) + f(y) \\ 
    &\ge 3k(x+y)\ell^4 - 12\ell
\end{align*}
by discarding positive terms and simplifying $s^{(-1)}(\ell)$.  Thus
\begin{align*}
\liminf_{x+y \to \infty} \big( \min_\ell F(\ell,x,y) \big)
    &\ge \liminf_{x+y \to \infty} \Big(\min_{\ell \ge 0}
    \big(3k(x+y)\ell^4 - 12\ell\big)\Big) \\
    &= \liminf_{x+y \to \infty} \frac{-9}{\sqrt[3]{k(x+y)}} = 0.
\end{align*}
The inequality comes from discarding positive terms, while the equality
follows just from the properties of $s^{(-1)}(\ell)$ that it is continuous,
and that it is positive for $\ell > 0$.

Having confined the locus of $F(\ell,x,y) \le -\eps$ to a compact region
for every $\eps > 0$, we will calculate derivatives and boundary values to
show that this region cannot have a local minimum and must thus be empty.
First, taking $\ell = 0$, we get
\[ F(0,x,y) = f(x) + f(y) > 0. \]
Second, taking $\ell = \pi$, we get
\[ F(\pi,x,y) = 4(x+y)-6\pi+f(x)+f(y). \]
Here we check that 
\[ \dbyd{F}{x}(\pi,x,y) = \frac{4x^4}{(x^2+1)^2} \ge 0, \qquad
    \dbyd{F}{y}(\pi,x,y) = \frac{4y^4}{(y^2+1)^2} \ge 0, \qquad 
    F(\pi,0,0) = 0. \]
Fourth, taking $x = y = 0$, we obtain
\[ F(\ell,0,0) = 9s^{(-2)}(\ell) - 12\ell + 6\pi.\]
We check in this case that
\[ F(\pi,0,0) = 0, \qquad \dbyd{F}{\ell}(\ell,0,0)
    = 9s^{(-1)}(\ell)-12 \le 0. \]

The fifth case is the case $y=0$ with $x$ and $\ell$ interior, which by
symmetry is equivalent to the case $x=0$ with $y$ and $\ell$ interior.
The sixth and final case is the interior for all three coordinates.
We will handle the fifth and sixth cases together.  Using the final change
of variables
\[ t \defeq \tan\big(\frac{\ell}2\big), \]
and with the help of Sage, we learn that
\begin{align*}
\dbyd{F}{\ell}(t,x,y) &= -12\frac{(t^4-t^2)xy - 2t^3(x+y) + 3t^2+1}
    {(t^2+1)^3} \\
\dbyd{F}{x}(t,x,y) &= 4\frac{(2t^3y-3t^2-1)(x^2 + 1)^2 + x^4(t^2+1)^3}
    {(t^2+1)^3(x^2 + 1)^2}.
\end{align*}
Note that the partial derivatives are rational functions in $x$, $y$,
and $t$.  We can rigorously determine the common zeroes of their numerators
by finding their associated prime ideals in the ring $\Q[x,y,t]$ using
the ``\texttt{associated{\_}primes}" function in Sage\footnote{See the
attached Sage files in the source file of the arXiv version of this paper.}
(In other words, we use the Lasker-Noether factorization theorem converted
to an algorithm by the Gr\"oebner basis method.)  The solution set is
characterized by five prime ideals:
\begin{align*}
I_1 &= (x - y, yt - 1) \\ 
I_2 &= (x + t, t^2 + 1) \\
I_3 &= (y + t, t^2 + 1) \\
I_4 &= (x, 2yt^3 - 3t^2 - 1) \\
I_5 &= \begin{multlined}[t]
    (2x^2y + 3x^2t + xyt + x + y, x^2t^3 + xyt^3 - xt^2 - yt^2 - 2x + 2t, \\
    xy^2t^2 + 2xyt^3 + y^2t^3 - 2xyt - 3xt^2 - yt^2 - y - t, \\
    y^2t^4 + y^2t^2 + xt^3 + 3yt^3 + 2xy + 3xt - 3yt - 7t^2 - 1, \\
    xyt^4 - xyt^2 - 2xt^3 - 2yt^3 + 3t^2 + 1). \end{multlined}
\end{align*}
Four of these ideals cannot vanish when $x > 0$ and $y,\ell \ge 0$: $I_2$
and $I_3$ contain $t^2+1$, $I_4$ contains $x$, and $I_5$ contains
\[ 2x^2y + 3x^2t + xyt+x+y > 0. \]
The ideal $I_1$ yields the desired locus $x = y = 1/t$.

A careful examination of the equality cases shows that $x = y = 1/t$
is the only possibility for the minimum value $F = 0$.
\end{proof}

\subsection{The negative case}
\label{s:negative}

In this section we will establish Theorem~\ref{th:negative}.  As in
Section~\ref{s:positive}, we let $\kappa = -1$ and we must take $n \in
\{2,4\}$.  We cannot use Problem~\ref{lp:basic}, because in both dimensions,
one of the dual coefficients turns out to be negative.  Instead we the
use the following model, which is provided by Lemma~\ref{l:extend}.

\begin{linprog} Given $n$, $A$, $V$, and $L$, let
\[ q = \frac{n-1}{\tanh(L)}. \]
Is there a symmetric, positive measure $\mu(\ell,\alpha,\beta)$ such that
\begin{align*}
\alpha_*(\mu) = \int_{\ell,\beta} \back \dd \mu &= A\sdd z(\alpha) \\
\int_{\ell,\alpha,\beta} \!\! \Big(s(\ell)\sec(\beta) -
    qs^{(-1)}(\ell)\Big)\sec(\alpha) \sdd \mu_\Omega &\le A^2-qAV. \\
\int_{\ell,\alpha,\beta} \! \big(s^{(-1)}(\ell)\sec(\alpha) -
    qs^{(-2)}(\ell)\big) \sdd \mu_\Omega &\le AV - qV^2 \\
\int_{\ell,\alpha,\beta} \back s^{(-2)}(\ell) \sdd\mu &\le V^2 \\
\int_{\ell,\alpha,\beta} \back \ell \sdd\mu &= \omega_{n-1} V?
\end{align*}
\eatline \label{lp:extend} \end{linprog}

We will need the dual problem, which we can state without changing variables.

\begin{linprog}[Dual to Problem~\ref{lp:extend}] Given $n$, $A$, $V$,
and $L$, let
\[ q = \frac{n-1}{\tanh(L)}. \]
Are there numbers $a,b,c, d \in \R$ and a continuous function $f:[0,\pi/2)
\to \R$ such that
\begin{gather}
a \ge 0 \qquad 2b+qa \ge 0 \qquad c+q(2b+qa) \ge 0 \label{e:slack} \\
\begin{multlined}
a s_{n,-1}(\ell)\sec(\alpha)\sec(\beta)
+ b s_{n,-1}^{(-1)}(\ell)\big(\sec(\alpha) + \sec(\beta)\big) \\
+ c s_{n,-1}^{(-2)}(\ell) - d\ell + f(\alpha) + f(\beta) \ge 0
\end{multlined} \nonumber \\
aA^2 + 2bAV + cV^2 - d\omega_{n-1}V
+ 2A\!\!\int_0^{\pi/2} \back \back f(\alpha) \sdd z(\alpha) < 0?
\nonumber \end{gather}
\eatline \label{lp:edual} \end{linprog}

We will use the change of variables
\begin{eq}{e:nxy}
(x,y) \defeq \Big(\frac{\sec(\alpha)}{\tanh(r)},
    \frac{\sec(\beta)}{\tanh(r)}\Big)
\end{eq}
with the range $x,y \in (1,\infty)$.  Equation \eqref{e:chord} simplifies to
\begin{eq}{e:xchordh}
\tanh(\frac{\ell}2) = \frac1x.
\end{eq}

\subsubsection{Dimension $2$}

In dimension $n=2$,
\begin{align*}
s(\ell) &= \sinh(\ell), & s'(\ell) &= \cosh(\ell), \\
s^{(-1)}(\ell) &= \cosh(\ell)-1, & s^{(-2)}(\ell) &= \sinh(\ell) - \ell.
\end{align*}
Equation~\eqref{e:check}, with \eqref{e:chord}, becomes
\[ \frac{a\tanh(r)^2\cosh(\ell)}{\tanh(\ell/2)^2}
    + \frac{2b\tanh(r)\sinh(\ell)}{\tanh(\ell/2)}
    + c(\cosh(\ell)-1) - d = 0. \]
The solution is
\[ a = 0, \qquad b = \frac{1}{\tanh(r)}, \qquad c = -2, \qquad d = 4. \]
We need to check the third case of condition \eqref{e:slack}, which
reduces to
\[ c + 2qb = -2 + \frac{2}{\tanh(r)\tanh(L)} \ge 0. \]
Since the $\tanh$ function is bounded above by $1$, this is immediate.

In the variables \eqref{e:nxy}, the cost function \eqref{e:cost} is
\[ E(\ell,x,y) = (\cosh(\ell)-1)(x+y) - \sinh(\ell) -2\ell. \]
The potential \eqref{e:potential} is
\[ f(x) = 2\arctanh(\frac{1}x). \]
The adjusted cost \eqref{e:confirm} separates according to
\eqref{e:sep} with
\[ G(\ell,x) = (\cosh(\ell)-1)x - \sinh(\ell) - \ell + 2\arctanh(\frac{1}x). \]
We minimize $G$ using the derivative test in $x$ to obtain
\[ \dbyd{G}{x}(\ell,x) = \frac{x^2(\cosh(\ell)-1)-(\cosh(\ell)+1)}{x^2-1}. \]
We learn that the minimum of $G$ in $x$ for each fixed $\ell$ occurs when
$x$ and $\ell$ satisfy \eqref{e:xchordh} and it is easy to confirm that
the value is 0.  Thus $G(\ell,x) \ge 0$, which confirms \eqref{e:confirm}
and establishes the $n=2$ case of Theorem~\ref{th:negative}.

\subsubsection{Dimension $4$}

In dimension $n=4$,
\begin{align*}
s(\ell) &= \sinh(\ell)^3, \\
s'(\ell) &= 3\cosh(\ell)\sinh(\ell)^2, \\
s^{(-1)}(\ell) &= \frac{\cosh(\ell)^3 - 3\cosh(\ell) + 2}3, \\
s^{(-2)}(\ell) &= \frac{\sinh(\ell)^3 - 6\sinh(\ell) + 6\ell}9.
\end{align*}

Equation~\eqref{e:check} becomes
\[ \frac{3a\tanh(r)^2\cosh(\ell)\sinh(\ell)^2}{\tanh(\ell/2)^2}
    + \frac{2b\tanh(r)\sinh(\ell)^3}{\tanh(\ell/2)}
    + \frac{c(\cosh(\ell)^3-3\cosh(\ell)+2)}{3} - d = 0. \]
The solution is
\[ a = \frac{1}{\tanh(r)^2}, \qquad b = -\frac{3}{\tanh(r)}, \qquad c = 9,
    \qquad d = 12. \]    
We need to check the second case of condition \eqref{e:slack}:
\[ 2b + qa = -\frac{6}{\tanh(r)} + \frac{3}{\tanh(L)\tanh(r)^2} \ge 0. \]
This condition is equivalent to the smallness hypothesis \eqref{e:small}.

Following the notation of Lemma~\ref{l:dual}, we claim that our choices
for $a,b,c,d$ and a suitable choice of $f$ produce a polynomial $P(A)$
with two real roots $A_0 < A_1$, where $A_1$ is the sharp isoperimetric
value, that together are a solution to Problem~\ref{lp:edual} when $A \in
(A_0,A_1)$.   In other words, we want to confirm \eqref{e:confirm}, but this
is the hardest part of the proof and we save it for the end of the section.
In the meantime, we settle a different and simpler difficulty.  Recall that
\[ V = |B_{4,-1}(r)| = \omega_3 s^{(-1)}(r), \qquad
    A_1 = |\del B_{4,-1}(r)|  = \omega_3 s(r), \]
and that $A_0 < 0$ is equivalent to the clean optimality condition
\eqref{e:clean}.   Since $s^{(-1)}(r)$ is unbounded, \eqref{e:clean}
does not hold for all $V$.  This criterion is not needed to prove
Theorem~\ref{th:negative}, since the set of geometrically feasible
$A$ is the closed ray $[A_1,\infty)$; but it is important for
Theorem~\ref{th:subsume}, which asserts that every $A \in [0,A_1)$
is infeasible for Problem~\ref{lp:extend}.  To prove this part of
Theorem~\ref{th:subsume}, we claim a second set of values
$a_2,b_2,c_2,d_2,f_2$ that are feasible for Problem~\ref{lp:edual} and
that produce $P_2(A)$ with $P_2(0) < 0$ and $P_2(A_1) < 0$.

Let $A_2 = 3\tanh(r)V$.  On the one hand, the proof of Theorem~\ref{th:yau}
in Section~\ref{s:yaulin} produces the feasible values
\[ a_2 = 0, \qquad b_2 = 1, \qquad c_2 = -3,\qquad d_2 = 0, \qquad f_2 = 0, \]
and a linear $P_2(A)$ that vanishes at $A_2$.  Thus, every $A \in [0,A_2)$
is infeasible.  On the other hand, we claim that $P(A_2) < 0$.  We use
$P(A_1) = 0$ to calculate that
\[ \int_0^{\pi/2} \back \back f(\alpha) \sdd z(\alpha)
    = \frac{d\omega_3 V - aA_1^2 -  2bA_1V - cV^2}{2A_1}
    = \frac{4\pi^2(\cosh(r)-1)}{\sinh(r)} \]
and
\[ P(A_2) = \Big(\frac{A_2}{\tanh(r)} - 3V\Big)^2 - 12\omega_3V
    + 2A_2\!\!\int_0^{\pi/2} \back \back f(\alpha) \sdd z(\alpha)
    = -\frac{24\pi^2}{\cosh(r)} < 0, \]
so that $A_2 \in (A_0,A_1)$, as desired.

It remains to confirm \eqref{e:confirm} for the given values of $a,b,c,d$.
The cost function \eqref{e:cost} is
\[ E(\ell,x,y) = \sinh(\ell)^3 xy - (\cosh(\ell)^3-3\cosh(\ell)+2)(x+y)
    + \sinh(\ell)^3 - 6\sinh(\ell) - 6\ell. \]
The potential from \eqref{e:potential} and \eqref{e:xchordh} is
\[ f(x) = 6\arctanh(\frac1x) + \frac{2x}{x^2-1}. \]
The adjusted cost \eqref{e:confirm} is
\begin{multline}
F(\ell,x,y) = \sinh(\ell)^3xy - (\cosh(\ell)^3-3\cosh(\ell)+2)(x+y)
    + \sinh(\ell)^3 - 6\sinh(\ell) - 6\ell \\ + 6\arctanh(\frac1x)
    + \frac{2x}{x^2-1} + 6\arctanh(\frac1y) + \frac{2y}{y^2-1}.
\label{e:Fn4} \end{multline}
We conclude the proof of Theorem~\ref{th:negative} with the following
lemma.  The lemma is also numerically evident but surprisingly tricky
(for the authors).

\begin{lemma} The function $F(\ell,x,y)$ on $\R_{\ge 0} \times (1,\infty)^2$
given by \eqref{e:Fn4} is non-negative, and vanishes only when
\[ x = y = \frac{1}{\tanh(\ell/2)}. \]
\end{lemma}

\begin{proof} The proof is analogous to that of Lemma~\ref{l:tech4p}, but
differs in its technical details.  Throughout the proof, we will fix $y$
and minimize $F(\ell,x,y)$ with respect to $x$ and $\ell$.

To check the non-compact limits of $x$ and $\ell$, we re-express $F$ as:
\[ F(\ell,x,y) = \sinh(\ell)^3(x-1)(y-1) + h(\ell)(x+y)
    - 6\sinh(\ell) - 6\ell + f(x) + f(y), \]
where 
\[ h(\ell) = (\sinh(\ell)^3 - \cosh(\ell)^3+3\cosh(\ell)-2)
    = \frac{(3e^\ell + 1)(1 - e^{-\ell})^3}4 > 0. \]
We also have 
\[ f(x) = \frac1{x-1} + \frac1{x+1} + \arctanh(\frac1x) > \frac1{x-1} \]
and the elementary relation $\sinh(\ell) \ge \ell$.  We combine these
comparisons to obtain the bound
\[ \hF(\ell,x,y) \defeq \sinh(\ell)^3(x-1)(y-1) + \frac1{x-1}
    - 12\sinh(\ell) < F(\ell,x,y). \]
The function $\hF$ is useful for minimizing with respect to either $\ell$
or $x$, leaving the other variables fixed.  It is a bit simpler to use
the variables
\[ (x_1,y_1) \defeq (x-1,y-1), \]
which we will need anyway later in the proof.  We obtain
\begin{align*}
\hF(\ell,x_1,y_1) &= \sinh(\ell)^3x_1y_1 + \frac1{x_1} - 12\sinh(\ell) \\
\min_\ell \hF(\ell,x_1,y_1) &= \frac{-16}{\sqrt{x_1y_1}} + \frac1{x_1} \\
\min_{x_1} \hF(\ell,x_1,y_1) &= 2\sqrt{\sinh(\ell)^3y_1} - 12\sinh(\ell).
\end{align*}
We obtain these uniform lim infs:
\begin{align*}
\liminf_{x \to \infty} \big(\inf_\ell F(\ell,x,y)\big)
    &\ge \lim_{x_1 \to \infty} \big(\min_\ell \hF(\ell,x_1,y_1)\big) = 0 \\
\liminf_{x \to 1} \big(\inf_\ell F(\ell,x,y)\big)
    &\ge \lim_{x_1 \to 0} \big(\min_\ell \hF(\ell,x_1,y_1)\big) = \infty \\
\liminf_{\ell \to \infty} \big(\inf_x F(\ell,x,y)\big)
    &\ge \lim_{\ell\to\infty} \big(\min_{x_1}\hF(\ell,x_1,y_1)\big) = \infty \\
\liminf_{\ell \to 0} \big(\inf_x F(\ell,x,y)\big)
    &\ge \lim_{\ell \to 0} \big(\min_{x_1} \hF(\ell,x_1,y_1)\big) = 0. \\
\end{align*}
Once we control $x$, we can also check the last case more directly by
calculating that
\[ F(0,x,y) = f(x) + f(y) > 0. \]
Either way, this establishes that we can use the derivative test for each
fixed $y$ to confirm that $F(\ell,x,y) \ge 0$.

We use the final change of variables
\[ t \defeq \tanh(\frac{\ell}2). \]
Sage tells us that
\begin{align*}
\dbyd{F}{\ell}(t,x,y) &= -12\frac{(t^4+t^2)xy - 2t^3(x+y) + 3t^2-1}
    {(t^2-1)^3}, \\
\dbyd{F}{x}(t,x,y) &= -4\frac{(2t^3y-3t^2+1)(x^2 - 1)^2 + x^4(t^2-1)^3}
    {(t^2-1)^3(x^2 - 1)^2}.
\end{align*}
We again rigorously determine the common zeroes of their numerators by
finding their associated prime ideals in the ring $\Q[x,y,t]$ using Sage.
The solution set in this case is characterized by 7 prime ideals:
\begin{align*}
I_1 &= (x - y, yt - 1) \\
I_2 &= (t + 1, x + 1) \\
I_3 &= (t + 1, y + 1) \\
I_4 &= (t - 1, x - 1) \\
I_5 &= (t - 1, y - 1) \\
I_6 &= (x, 2yt^3 - 3t^2 + 1) \\
I_7 &= \begin{multlined}[t]
    (2x^2y - 3x^2t - xyt - x - y, x^2t^3 + xyt^3 - xt^2 - yt^2 + 2x + 2t, \\
    xy^2t^2 - 2xyt^3 - y^2t^3 - 2xyt + 3xt^2 + yt^2 - y + t, \\
    y^2t^4 - y^2t^2 + xt^3 + 3yt^3 + 2xy - 3xt + 3yt - 7t^2 + 1, \\
    xyt^4 + xyt^2 - 2xt^3 - 2yt^3 + 3t^2 - 1).
    \end{multlined}
\end{align*}
The ideal $I_1$ yields the desired locus $x = y = 1/t$, while the other
six do not vanish when $0 \le t < 1$ and $x,y > 1$.  Five of these cases
are easy:  The ideals $I_2$ and $I_3$ contain $t+1$, the ideals $I_4$
and $I_5$ contain $t-1$, and the ideal $I_6$ contains $x$.

The ideal $I_7$ is obviously more complicated.  Setting the first
two generators to zero, we obtain:
\begin{align*}
2x^2y - 3x^2t - xyt - x - y &= 0 \\
x^2t^3 + xyt^3 - xt^2 - yt^2 + 2x + 2t &= 0.
\end{align*}
Since the first equation is linear in $t$, we can eliminate it by
substitution, and then clear the denominator and eliminate a factor of $x+y$.
The resulting equation in $x$ and $y$ is
\[ 4x^6y^3 - 12x^5y^2 - 8x^4y^3 + 27x^5 + 27x^4y + 17x^3y^2 + 5x^2y^3
    - 11x^3 - 11x^2y - 5xy^2 - y^3 = 0. \]
We can express this in the variables $x_1$ and $y_1$ as
\begin{multline*}
4x_1^6y_1^3 + 12x_1^6y_1^2 + 24x_1^5y_1^3 + 12x_1^6y_1 + 60x_1^5y_1^2
    + 52x_1^4y_1^3 + 4x_1^6 + 48x_1^5y_1 + 96x_1^4y_1^2 \\
    + 48x_1^3y_1^3 + 39x_1^5 + 63x_1^4y_1 + 41x_1^3y_1^2 + 8x_1^2y_1^3
    + 154x_1^4 + 46x_1^3y_1 + 312x_1^3 + 55x_1^2y_1 \\
    + 336x_1^2 + 56x_1y_1 + 176x_1 + 16y_1 + 32
    + 9x_1^2y_1(y_1-1)^2+2x_1y_1(y_1-2)^2 = 0.    
\end{multline*}
The left side is manifestly a sum of positive terms when $x_1,y_1 > 0$
and thus cannot vanish.  Thus $I_7$ cannot vanish when $x,y > 1$, which
completes the derivative test for $F(\ell,x,y)$.
\end{proof}

\section{Proofs of other results}
\label{s:other}

\subsection{Uniqueness}
\label{s:uniq}

Problems~\ref{lp:basic} and \ref{lp:extend} both place strong restrictions
on $\mu$ and therefore on $\Omega$ in the sharp case.   First, all of the
inequalities in Problem~\ref{lp:basic} become equalities when $\kappa >
0$; all of the inequalities in Problem~\ref{lp:extend} become equalities
when $\kappa < 0$.  In particular, equation \eqref{e:teufel} becomes
an equality, which implies that $\Omega$ is convex and that the candle
comparison is an equality at short distances.  That in turn implies that
$\Omega$ satisfies $\Ric \ge (n-1)\kappa g$ and that it is the equality case
of Bishop's inequality \cite[Sec. 11.10]{BC:book}, which implies that it
has constant curvature $K = \kappa$.

The case $\kappa = 0$ does not use \eqref{e:teufel}, but it does use
\eqref{e:croke}.  This again implies that $\Omega$ is convex.  The stronger
assumption that $\Omega$ is $\RRic$ class $0$ together with equality in
\eqref{e:croke} tells us again that $\Omega$ has constant curvature $K = 0$.

Second, sharpness tells us that $\mu_\Omega$ is concentrated on the locus
given by equation \eqref{e:chord}.  In other words, every chord in $\Omega$
has the same length and incident angles as if $\Omega$ were a round ball
$B_{n,\kappa}(r)$.  If $\Omega$ is convex with constant curvature, this
implies that $\Omega$ is isometric to $B_{n,\kappa}(r)$.

\subsection{G\"unther's inequality with reflections}
\label{s:mgunther}

In this section, we will prove Proposition \ref{p:mgunther}.

If $\gamma(t)$ is a smooth curve in $M$ with $t \in [0,r]$, then it is
a constant-speed geodesic if and only if it is a critical point of the
energy functional
\[ E(\gamma) = \int_0^r \frac{\braket{\gamma'(t),\gamma'(t)}}2 \,\dd t\]
assuming Dirichlet boundary conditions (\ie, that we fixed the endpoints
of $\gamma$).   Let $\gamma$ be such a geodesic with unit speed, and let
$y(t)$ be a smooth, infinitesimal normal displacement.  Then we can define
a relative energy
\[ E(y) \defeq E(\gamma+y) - E(\gamma) + O(||y||^3), \]
which is just the second variational derivative of the curve energy,
equivalently half of the second variation of the curve length.  We can
identify the normal bundle to $\gamma(t)$ with $\R^{n-1}$ using parallel
transport, thus view $y$ as a function with values $y(t) \in \R^{n-1}$.
If $\gamma$ is an ordinary geodesic without reflections, then by a standard
calculation,
\[ E(y) = \int_0^r \big[\braket{y'(t),y'(t)} - \braket{y(t),R(t)y(t)}\big]
    \sdd t, \]
where
\[ R(t) = R(\cdot,\gamma'(t),\cdot,\gamma'(t)) \]
is the Riemann curvature tensor specialized at the unit tangent $\gamma'$.
This leads to the differential equation
\begin{eq}{e:jacobi} y''(t) = -R(t)y(t), \end{eq}
which is satisfied by $y$ when it is a Jacobi field, \ie, a geodesic
displacement of $\gamma$.

\begin{figure}[htb]
\begin{center}
\begin{tikzpicture}[scale=1.25]
\draw[thick,darkblue,->] (-3,0) -- (-1.5,0);
\draw[thick,darkblue] (-1.5,0) -- (0,0);
\draw[thick,darkblue,dashed] (0,0) -- (3,0);
\draw[thick,darkblue,->] (0,0) -- (0,-1);
\fill[black!30!white] (.795,-1.035) arc (30:60:5)
    -- (-.795,1.035) arc (210:240:5) -- cycle;
\draw[thick] (.795,-1.035) arc (30:60:5);
\draw[dashed,thick] (-.795,1.035) arc (210:240:5);
\draw[darkgreen] (-3,.560) .. controls (-2,.26) and
    (-1.5,.76) .. (-.668,.560);
\draw[very thick,darkred] (-.668,.560) -- (-.560,.668);
\draw[darkgreen,dashed] (-.560,.668) .. controls (.46,.868) and (2,.268)
    .. (3,.668);
\draw (-3,0) node[anchor=east] {$\gamma(t)$};
\draw (-3,.560) node[anchor=east] {$(\gamma+y)(t)$};
\draw (.795,-1.035) node[anchor=north east] {$\dM$};
\draw[gray] (1.035,-.795) node[anchor=south west] {$\dM$};
\end{tikzpicture}
\end{center}
\caption{Diagram of a vector field $y$ that displaces a geodesic $\gamma$
    in a curved surface (non-geodesically), and a continuation if $\gamma$
    were straight.  The short red segment is the length variation of
    $\gamma+y$ due to the reflection.}
\label{f:reflect}\end{figure}
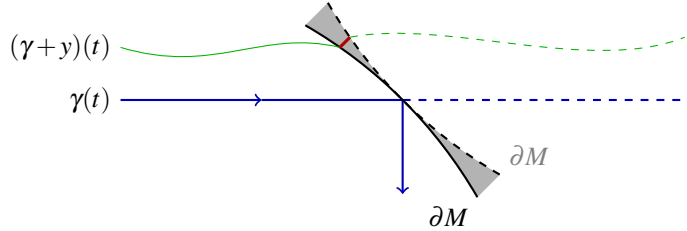

If $\gamma$ reflects from $\dM$, then the energy has extra terms.  We will
derive the energy \eqref{e:energy} and a modified Jacobi field equation
\eqref{e:distrib}.  Although these equation are not really new \cite[Sec.
2]{Innami:billiards}, we give a geometric argument that we have not seen
elsewhere. To understand the extra terms in $E(y)$ due to the reflections,
suppose that $\gamma$ reflects from $\dM$ at a point $p = \gamma(t)$,
and let $Q = Q(p)$ be shape operator $\dM$ relative to the inward unit
normal $w = w(p)$, \ie,
\[ Qu= -\nabla_u w. \]
If we give $\gamma$ a ghost extension as in Figure~\ref{f:reflect}, then
the displacement $\gamma+y$ has a gap when $\dM$ is curved.  (The figure
shows the convex case with a positive gap; the gap can also have negative
width.) We first assume the simplest case in which $\gamma$ is normal
to $\dM$.  The quadratic form $\langle\cdot,Q\cdot\rangle/2$ osculates
$\dM$, so that the width of the gap, and thus the negative of the change
in length, is $\braket{y,Qy}$.  If the angle of incidence of $\gamma$ is
$\theta \ne 0$, then this answer is subject to two corrections.   First,
the gap is at an angle of $\theta$ from $\gamma$, so the length saved is
$\cos(\theta)\braket{\cdot,Q\cdot}$.  Second, $y$ no longer represents the
position that $\gamma+y$ meets $T_p\dM$, again because the surface is angled.

To derive where $\gamma+y$ meets $T_p\dM$, we call $T_p(\dM)$ the
\emph{tangent hyperplane}, the normal $N_p(\gamma)$ to $\gamma(t)$ the
\emph{coronal hyperplane}, and the 2-dimensional plane spanned by $w(p)$
and $\gamma'(t)$ the \emph{sagittal plane}\footnote{This terminology
is borrowed from human anatomy.}.  Let $P$ be the orthogonal projection
from the tangent hyperplane to the coronal hyperplane.   If we choose an
orthonormal coronal basis $e_1,\ldots,e_{n-1}$ such that $e_1$ is in the
sagittal plane, and a matching tangent basis, then
\[ P = \begin{pmatrix} \cos(\theta) & 0 & \cdots & 0 \\ 0 & 1 & & 0 \\
\vdots & & \ddots & \vdots \\ 0 & 0 & \cdots & 1 \end{pmatrix}. \]
Then the change in length, and therefore the extra energy term, is
\[ -\cos(\theta)\braket{P^{-1}y,QP^{-1}y}. \]
(This formula still works when $\theta = 0$ if we take $P$ to be the
identity matrix.)  If $\gamma$ reflects from a sequence of boundary
points $\{p_k\}$ at times $\{t_k\}$, then we have the same change in
length using angles $\theta_k$ and symmetric matrices $P_k$ and $Q_k$,
and we can abbreviate the result by letting
\[ A_k \defeq \cos(\theta_k)P_k^{-1}Q_kP_k^{-1}. \]
Then energy of the normal field $y$ is
\begin{eq}{e:energy}
E(y) = \int_0^r \big[\braket{y'(t),y'(t)}
    - \braket{y(t),R(t)y(t)}\big] \sdd t
    - \sum_k \braket{y(t_k),A_k y(t_k)}.
\end{eq}
Thus, if $y$ is a (reflecting) Jacobi field, it satisfies the distributional
differential equation
\begin{eq}{e:distrib}
y''(t) = -R(t)y(t) - \sum_k A_k y(t_k) \delta_{t_k}(t),
\end{eq}
where $\delta_t$ is a Dirac delta measure on $\R$ concentrated at $t$.
Note that if $\dM$ is concave, then $Q_k$ is negative semidefinite and
therefore so is $A_k$.

We now follow a standard proof of G\"unther's inequality
\cite[Thm. 3.101]{GHL:riemannian}.  First, we will need that the energy
\eqref{e:energy} is positive definite, so that if $y$ Jacobi field,
it is an energy minimum (assuming Dirichlet boundary conditions) and
not just a critical point.  This is standard in the proof of G\"unther's
inequality without the $A_k$ terms, with the aid of the length restriction
when $\kappa > 0$.  It is still true with the $A_k$ terms, since each such
term is positive semidefinite.

Second, we consider a matrix solution $Y$ to \eqref{e:distrib} with $Y(0)
= 0$ and $Y(r) = I$.   Then the candle function of $\gamma$ satisfies
\[ j(\gamma,0,\ell) = \frac{\det Y(\ell)}{\det Y(r)}, \]
and the logarithmic derivative at $r$ is given by
\[ \dbyd{}{t}\Big|_{t=r} \back\log(j(\gamma,0,t))
    = (\det Y)'(r) = \Tr(Y)(r). \]
We generalize the energy \eqref{e:energy} to the matrix argument $Y$,
and we interpret it as a function of $Y$, $R(t)$, and each $A_k$:
\begin{eq}{e:menergy}
E(Y,R,A) = \int_0^r \big[\braket{Y'(t),Y'(t)} -
    \braket{Y(t),R(t)Y(t)}\big] \sdd t - \sum_k \braket{Y(t_k),A_k Y(t_k)},
\end{eq}
using the Hilbert-Schmidt inner product
\[ \braket{X,Y} = \Tr(X^TY). \]
If $Y$ is a solution to \eqref{e:distrib}, then integration by parts yields
the remarkable equality
\[ \Tr(Y)(r) = E(Y,R,A) .\]
If we minimize $E$ with respect to all three arguments $Y$, $R$, and $A$,
then we both solve \eqref{e:distrib} and minimize the logarithmic derivative
of $\gamma$.  If we fix $Y$, then it is immediate from \eqref{e:menergy}
and from the constraints that we should take $R = \kappa I$ and $A_k =
0$, \ie, maximum curvature and flat mirrors.

\subsection{Multiple images}
\label{s:mulproof}

Lemma~\ref{l:multiple} yields the following model.

\begin{linprog} Given $n$, $\kappa$, $A$, $V$, and $m$, is there a symmetric,
positive measure $\mu(\ell,\alpha,\beta)$ such that
\begin{align*}
\alpha_*(\mu) = \int_{\ell,\beta} \back \dd \mu
    &= A\sdd z(\alpha) \\
\int_{\ell,\alpha,\beta} \back s_{n,\kappa}(\ell)\sec(\alpha)\sec(\beta)
    \sdd\mu &\le mA^2 \\
\int_{\ell,\alpha,\beta} \back s_{n,\kappa}^{(-1)}(\ell)
    \big(\sec(\alpha)+\sec(\beta)\big) \sdd\mu &\le 2mAV \\
\int_{\ell,\alpha,\beta} \back s_{n,\kappa}^{(-2)}(\ell) \sdd\mu
    &\le mV^2 \\
\int_{\ell,\alpha,\beta} \back \ell \sdd\mu &= \omega_{n-1} V? 
\end{align*}
\eatline \label{lp:multiple} \end{linprog}

Theorem~\ref{th:multiple} now follows as a porism\footnote{A corollary of
proof.} of Theorem~\ref{th:positive}.  If we apply the transformation
\[ \tilde{V} = mV, \qquad \tilde{A} = mA, \qquad \tilde{\mu} = m\mu, \]
then Problem~\ref{lp:multiple} becomes Problem~\ref{lp:basic}.

\subsection{Alternative functionals}
\label{s:alldim}

In this section we prove Theorem~\ref{th:alldim}.  The proof is
almost the same as the proof of Croke's theorems in Section~\ref{s:zero}.

Given $L = L(\Omega)$, we consider the following linear programming problem
based on equation~\ref{e:croke}.

\begin{linprog}
Given $n$, $A$, and $L$, is there a symmetric
positive measure $\mu(\ell,\alpha,\beta)$ such that
\begin{align*}
\alpha_*(\mu) = \int_{\ell,\beta} \dd \mu
    &= A\sdd z(\alpha) \\
\int_{\ell,\alpha,\beta} \back \ell^{n-1}\sec(\alpha)\sec(\beta)
    \sdd\mu &\le A^2 \\
\int_{\ell,\alpha,\beta} \back \ell^{n-3} \sdd\mu &= L?
\end{align*}
\label{lp:alldim} \end{linprog}

We can apply a version of Lemma~\ref{l:sharp}
to establish Theorem~\ref{th:alldim} as a sharp inequality.

Given $a \ge 0$ and $d \in \R$, we consider the cost function
\[ E(\ell,\alpha,\beta) = a\ell^{n-1}\sec(\alpha)\sec(\beta) - d\ell^{n-3}. \]
By design, given a radius $r > 0$, there are values of $a,d>0$ such that
$E(\ell,\alpha,\alpha)$ is minimized in $\ell$ when
\[ \ell = 2r\cos(\alpha), \]
which thus satisfies \eqref{e:chord}.  We can take
\[ a = \frac{n-3}{r^2}, \qquad d = 4(n-1). \]
(Note that we need $n \ge 4$.  If $n < 3$, then $a$ would be negative.
If $n=3$, then $L(\Omega) \propto |\dOmega|$ and Theorem~\ref{th:alldim}
is vacuous.)

We define the potential
\[ f(\alpha) = -\frac{E(2r\cos(\alpha),\alpha,\alpha)}2 = 
    2^{n-1}(r\cos(\alpha))^{n-3}. \]
Applying the change of variables \eqref{e:zxy}, the adjusted cost function is
\[ F(\ell,x,y) = (n-3)\ell^{n-1}xy - 4(n-1)\ell^{n-3}
    + 2^{n-1}(x^{3-n} + y^{3-n}). \]
We want to show that $F \ge 0$.  For any fixed value of $xy$,
$F(\ell,x,y)$ is minimized when $x = y$.  Then
\begin{align*}
F(\ell,x,x) &= (n-3)\ell^{n-1}x^2 - 4(n-1)\ell^{n-3} + 2^nx^{3-n} \\
    &= \big((n-3)(\ell x)^{n-1} - 4(n-1)(\ell x)^{n-3} + 2^n\big)x^{3-n}.
\end{align*}
The first factor is a polynomial in $\ell x$ that, by univariate calculus,
decreases to $0$ at $\ell x = 2$ and then increases again.  This completes
the proof of Theorem~\ref{th:alldim}.

\subsection{Old wine in new decanters}
\label{s:oldwine}

In this section, we complete the proof of Theorem~\ref{th:subsume}.  The rest
of this paper has covered all cases except Theorem~\ref{th:yau}, due to Yau,
and Theorem~\ref{th:croke2}, due to Croke.  The arguments given here are
equivalent to the original proofs, only restated in linear programming form.

\subsubsection{Yau's linear isoperimetric inequality}
\label{s:yaulin}

If $\Omega$ is $n${\hyp}dimensional and $\LCD(-1)$, then
Problem~\ref{lp:extend} yields
\[ \int_{\ell,\alpha,\beta} \back \back \big(s^{(-1)}(\ell)\sec(\alpha) -
    (n-1)s^{(-2)}(\ell)\big) \sdd \mu_\Omega \le AV - (n-1)V^2 \]
since $q > n-1$.  The integrand is positive, since it is the second
antiderivative of
\[ s'(\ell)\sec(\alpha)-(n-1)s(\ell) = (n-1)\sinh(\ell)^{n-2}
    \big(\cosh(\ell)\sec(\alpha) - \sinh(\ell)\big) > 0. \]
Thus the right side is positive, and Theorem~\ref{th:yau} follows.  In terms
of optimal transport, the result follows if we define a cost function
\[ E(\ell,\alpha,\beta) = s^{(-1)}(\ell)\sec(\alpha) - (n-1)s^{(-2)}(\ell), \]
and then a vanishing potential $f(\alpha) = 0$.

\subsubsection{Croke's curvature-free inequality} 
\label{s:croke2}

For simplicity, we take $\rho = 1$.

Suppose that $\Omega$ is $n$-dimensional with unique geodesics.
Lemma~\ref{l:cbk} produces the following simple model independent of
$\kappa$, and that can be combined with Problem~\ref{lp:basic}.

\begin{linprog} Given $n$, $A$, and $V$, is there a symmetric,
positive measure $\mu(\ell,\alpha,\beta)$ such that
\begin{align*}
\alpha_*(\mu) = \int_{\ell,\beta} \back \dd \mu &= A\sdd z(\alpha) \\
\int_{\ell,\alpha,\beta} s_{n,(\pi/\ell)^2}^{(-2)}(\ell) \sdd\mu
    &\le V^2 \\
\int_{\ell,\alpha,\beta} \ell \sdd\mu &= \omega_{n-1} V?
\end{align*}
\eatline \label{lp:croke2} \end{linprog}

To analyze this model, we simplify it in two respects.  First,
we can integrate away $\alpha$ and $\beta$, because none of the integrals
explicitly depend on them.  We call the resulting measure $\mu(\ell)$.
Second, we can explicitly evaluate the integrand that arises
from Lemma~\ref{l:cbk}:
\[ s_{n,(\pi/\ell)^2}^{(-2)}(\ell) = \Big(\frac{\ell}{\pi}\Big)^{n+1}
    s_{n,1}^{(-2)}(\pi) = \frac{\ell^{n+1} \omega_n}{2\pi^n \omega_{n-1}}. \]
The first equality is just rescaling by $\ell/\pi$.  The second equality
is a tricky but standard integral; the answer can also be inferred from
the optimal case of a hemisphere $Y_{n,1}$.  The simplified model is then
as follows.

\begin{linprog} Given $n$, $A$, and $V$, is there a positive measure
$\mu(\ell)$ on $\R_{\ge 0}$ such that
\begin{align*}
\int_\ell \dd \mu &= \frac{\omega_{n-2}}{n-1}A \\
\int_\ell \frac{\ell^{n+1} \omega_n}{2\pi^n \omega_{n-1}}
    \sdd\mu_\Omega &\le V^2 \\
\int_\ell \ell \sdd\mu_\Omega &= \omega_{n-1} V
\end{align*}
\eatline \label{lp:crokesimp} \end{linprog}

As usual, we state the dual of Problem~\ref{lp:crokesimp}.

\begin{linprog} Given $n$, $A$, and $V$, are there constants
$c \ge 0$ and $f,d \in \R$ such that
\begin{align}
f + c \frac{\ell^{n+1} \omega_n}{2\pi^n \omega_{n-1}} - d\ell &\ge 0
    \nonumber \\
f \frac{\omega_{n-2}}{n-1}A + cV^2 - d\omega_{n-1}V &< 0? \label{e:crokedip}
\end{align}
\eatline \end{linprog}

In the optimal case of $Y_{n,1}$, we have $\ell = \pi$ everywhere and $V
= \omega_n/2$.  We can solve for the constants $f$, $c$, and $d$ assuming
that the left side of equation \eqref{e:crokedip} reaches 0 there and is
non-negative for other values of $\ell$.   We obtain
\[ c = 2\omega_{n-1}, \qquad d = (n+1)\omega_n, \qquad f = n\pi\omega_n. \]
Assuming that $A$ is feasible for Problem~\ref{lp:crokesimp},
\eqref{e:crokedip} then gives us the inequality
\[ A \ge \frac{(n-1)\omega_n\omega_{n-1}}{2\pi \omega_{n-2}} = \omega_{n-1}. \]
This establishes Theorem~\ref{th:croke2}.

\begin{remark} It may seem wrong that $Y_{n,1}$ does not itself have unique
geodesics.  But it is a limit of manifolds that do, which is good enough.
In any case the proof of Theorem~\ref{th:croke2} only really uses that
$\Omega$ has unique geodesics in its interior.
\end{remark}

\section{Closing questions}
\label{s:questions}

Of course, we want Theorem~\ref{th:negative} without the smallness
condition \eqref{e:small}.   It would suffice to prove a stronger version
of Lemma~\ref{l:extend}.  This would be implied by the $n=4$ case of the
following conjecture.

\begin{conjecture}  Let $j(r,t) = j_M(\gamma,r,t)$ be the candle function
of a geodesic in $\gamma$ in an $n$-manifold $M$ with curvature $K \le -1$.
Then
\begin{eq}{e:gunther2} \Big[(n-1)^2j - (n-1)\dbyd{j}{t} + (n-1)\dbyd{j}{r}
    - \frac{\del^2 j}{\del r \del t}\Big](r,t) \end{eq}
is minimized when $M$ has constant curvature $K = -1$.  
\label{c:candle} \end{conjecture}

As in the proof of Lemma~\ref{l:extend}, we would use Conjecture~\ref{c:candle}
to obtain the inequality
\begin{multline*}
\int_{\ell,\alpha,\beta} \Big(\frac{s_{n,-1}(\ell)}{\cos(\alpha)\cos(\beta)} -
    (n-1)s_{n,-1}^{(-1)}(\ell)\big(\frac1{\cos(\alpha)}
    + \frac1{\cos(\beta)}\big)
    + (n-1)^2 s_{n,-1}^{(-2)}(\ell) \Big) \sdd \mu_\Omega \\
    \le |\dOmega|^2- 2(n-1)|\dOmega| |\Omega| + (n-1)^2 |\Omega|^2,
\end{multline*}
and thus sharpen Problem~\ref{lp:extend}, by integrating over $(\gamma
\cap \Omega) \times (\gamma \cap \Omega)$ for a general geodesic $\gamma$.
If we integrate over a connected interval $[0,\ell]$, which suffices when
$\Omega$ is convex, then Conjecture~\ref{c:candle} implies that
\[ j(0,\ell) - (n-1) \int_0^\ell j(0,t) \sdd t
    - (n-1) \int_0^\ell j(s,\ell) \sdd s 
    + (n-1)^2 \int_0^\ell\int_0^t j(s,t) \sdd s \sdd t \]
is minimized when $K = -1$.  Note that even this relation is not true
under the weaker hypothesis $\LCD(-1)$.  For example, it does not hold
when $\ell$ is large enough if $M$ is the complex hyperbolic plane $\CH^2$,
normalized to be $(-9/4,-9/16)$-pinched.

The following relaxation of Kleiner's theorem is open even though, as
explained in Section~\ref{s:nonsharp}, it is close to true.  The motivation
is that the even strongest form holds in dimension $n=4$ following the
proof of Croke's theorem.

\begin{question} Suppose that $\Omega$ is a compact 3-manifold with boundary,
and with unique geodesics, non-positive curvature, and fixed volume $V =
|\Omega|$.  Then is its surface area $|\dOmega|$ minimized when $\Omega$
is a round, Euclidean ball?   What if non-positive curvature is replaced
by $\Candle(0)$?   What if $\Candle(0)$ is only required for pairs of
boundary points?
\label{q:kleiner}
\end{question}

Question~\ref{q:kleiner} could also be asked in dimension $n \ge 5$ and
for other curvature bounds $\kappa \ne 0$.

Finally, the following conjecture would give a more robust proof 
of Theorem~\ref{th:equality}, with a weaker hypothesis as well
when $\kappa = 0$.

\begin{conjecture} Suppose that $\Omega$ is a convex, compact Riemannian
$n${\hyp}manifold with boundary with unique geodesics.   Suppose that for
some constants $\kappa$ and $r$, all chords in $\Omega$ satisfy equation
\eqref{e:chord}.   Then $\Omega$ is isometric to $B_{n,\kappa}(r)$.
\end{conjecture}

\bibliographystyle{hamsalpha}
\bibliography{dg,mg,web,me}
\end{document}